\documentclass[aos,preprint]{imsart}

\RequirePackage{amsthm,amsmath,amsfonts,amssymb}
\RequirePackage[numbers,sort&compress]{natbib}
\RequirePackage[colorlinks,citecolor=blue,urlcolor=blue,linkcolor=blue]{hyperref}
\RequirePackage{graphicx}
\usepackage{mathtools,enumitem,booktabs,tabularx,array}
\usepackage{bibunits}

\usepackage{xcolor}

\defaultbibliographystyle{spr-ims-number}
\defaultbibliography{reference}
\startlocaldefs
\numberwithin{equation}{section}

\theoremstyle{plain}
\newtheorem{theorem}{Theorem}[section]
\newtheorem{proposition}[theorem]{Proposition}
\newtheorem{lemma}[theorem]{Lemma}
\newtheorem{corollary}[theorem]{Corollary}
\theoremstyle{definition}
\newtheorem{definition}[theorem]{Definition}
\newtheorem{remark}[theorem]{Remark}
\newtheorem{example}[theorem]{Example}
\newcommand{\R}{\mathbb R}
\newcommand{\Pp}{\mathcal P}
\newcommand{\Mm}{\mathcal M}

\newcommand{\E}{\mathbb E}
\newcommand{\Z}{\mathcal Z}
\newcommand{\K}{\mathcal K}
\newcommand{\F}{\mathcal F}
\newcommand{\dd}{\,\mathrm d}
\newcommand{\supp}{\operatorname{supp}}
\newcommand{\Lip}{\operatorname{Lip}}
\newcommand{\osc}{\operatorname{osc}}
\newcommand{\diam}{\operatorname{diam}}

\newcommand{\TV}{\operatorname{TV}}
\newcommand{\dist}{\operatorname{dist}}
\newcommand{\rank}{\operatorname{rank}}

\newcommand{\cA}{\mathcal A}
\newcommand{\ip}[2]{\left\langle #1,#2\right\rangle}
\newcommand{\norm}[1]{\left\|#1\right\|}
\newcommand{\abs}[1]{\left|#1\right|}

\newcommand{\bx}{\boldsymbol x}
\newcommand{\by}{\boldsymbol y}
\newcommand{\bm}{\boldsymbol m}
\newcommand{\bu}{\boldsymbol u}
\newcommand{\bv}{\boldsymbol v}
\newcommand{\ba}{\boldsymbol a}
\newcommand{\be}{\boldsymbol e}
\newcommand{\bxi}{\boldsymbol\xi}
\newcommand{\bz}{\boldsymbol z}
\newcommand{\bk}{\boldsymbol k}
\newcommand{\btheta}{\boldsymbol\theta}

\DeclareMathOperator*{\argmin}{arg\,min}
\endlocaldefs

  \definecolor{darkspringgreen}{rgb}{0.09, 0.45, 0.27} 
 \definecolor{darkgray}{rgb}{0.66, 0.66, 0.66}

\begin{document}

\begin{frontmatter}
\title{Random Basket Quotes: Minimax Distributional Recovery and Price Estimation
}
\runtitle{Random Basket Quotes}
\begin{aug}
\author[A]{\fnms{Tongseok}~\snm{Lim}\ead[label=e1]{lim336@purdue.edu}}
\author[A]{\fnms{Yichen}~\snm{Zhang}\ead[label=e2]{yichen@purdue.edu}}
\address[A]{Department of Quantitative Methods, Purdue University\printead[presep={,\ }]{e1,e2}}
\end{aug}

\begin{abstract}
We study the recovery of a compactly supported multivariate
probability law from random, possibly noisy basket-call prices.
Exact quotes define a convex set of admissible laws, whose
Wasserstein diameter equals the largest remaining price interval
over Lipschitz claims. After a finite-dimensional moment
correction, the $L^2$ pricing discrepancy is proportional to
a homogeneous Sobolev norm of order $-(d+3)/2$. For densities in an $H^\alpha$ ball, we establish the sharp
minimax expected $W_1$ risk
\(
 n^{-(\alpha+1)/d}
 +(\varsigma^2/n)^{(\alpha+1)/(2\alpha+2d+3)},
\)
uniformly over bounded noise levels, including exact quotes.
A bounded systematic quote error of size $\epsilon$ adds
the term $\epsilon^{(\alpha+1)/(\alpha+(d+3)/2)}$.
The upper bound is attained by a probability-valued estimator
combining local reconstruction of the pricing surface with
residual-corrected Fourier moments. For a given payoff, we characterize minimax price estimation
at fixed positive noise through an $L^2$ modulus of continuity.
Its convex dual is a penalized static replication problem,
which yields an affine estimator attaining the minimax rate. Finally, long-only basket quotes identify the joint law and estimate every moment of the spread $x_1-x_2$ at the parametric rate, yet the exchange option $(x_1-x_2)_+$ has minimax risk of order $(\log n)^{-(\alpha+3/2)}$, uniformly down to zero noise. 
\end{abstract}


\end{frontmatter}
\begin{bibunit}

\section{Introduction}\label{sec:intro}
Many statistical problems observe an unknown probability measure $\mu$ only through random linear functionals. A design point $Z_i$ is drawn, and the value $\int\phi_{Z_i}\dd\mu$ of a known test function is recorded, exactly or with noise. Exact observations determine a convex set of measures consistent with the data, and noisy observations define a statistical inverse problem. The difficulty of recovering $\mu$, or a single functional $\int g\dd\mu$, depends on the geometry of the test functions, the design, the probability constraints on the unknown, and the loss. We study these effects for a family of test functions arising in multivariate pricing.

We consider basket calls at one fixed maturity. A direction $\bu\in S^{d-1}$ and strike $k$ define the payoff $(\bu\cdot\bx-k)_+$, whose price under a law $\mu\in\Pp(\R^d)$ is $C_\mu(\bu,k)=\int(\bu\cdot\bx-k)_+\,\mu(\mathrm d\bx)$. We observe the direction, strike, and price of each sampled basket. A finite collection of such prices imposes linear moment constraints on the pricing law, leaving an interval of admissible prices for an unquoted payoff. We treat recovery of the whole law and estimation of a single claim as two uses of this common information structure.

For distributional recovery, we measure error in $W_1$, the 1-Wasserstein distance. This metric controls all Lipschitz payoffs, and the $W_1$ diameter of a quote fiber equals the largest width of their admissible price intervals. The pricing operator has smoothing order $s=(d+3)/2$: after correcting for the first moments, its $L^2$ pricing discrepancy is proportional to a homogeneous negative Sobolev norm of order $-s$. We use this identity to establish sharp inverse moduli for arbitrary probability laws and Sobolev density classes. These results quantify how discrepancies between pricing surfaces translate into uncertainty about the underlying law.

Our first main result determines the minimax expected risk over a Sobolev density class $\F_\alpha(L)\subset H^\alpha(\R^d)$ supported strictly inside the observed strike range. Under uniform i.i.d.\ quote locations and independent errors of standard deviation $\varsigma$, Theorem~\ref{thm:nu-smooth} establishes the sharp minimax rate
\[
 \mathfrak R_n\bigl(\F_\alpha(L),\varsigma\bigr)
 \asymp n^{-\frac{\alpha+1}{d}}
 +\Big(\frac{\varsigma^2}{n}\Big)^{\frac{\alpha+1}{2\alpha+2d+3}},
 \qquad 0\le\varsigma\le\overline\varsigma.
\]
The first term is the information limit imposed by finitely many exact random quotes; the second is the noise-dependent contribution at fixed positive noise. They have the same order at $\varsigma\asymp n^{-(2\alpha+d+3)/(2d)}$. The theory also accommodates bounded systematic quote errors: an error of size $\epsilon$, such as a mid-quote bias, adds $\epsilon^{(\alpha+1)/(\alpha+(d+3)/2)}$ to the rate. Unlike independent noise, this contribution does not vanish as the number of quotes grows; see Corollary~\ref{cor:nu-marginals}.

A probability-valued estimator combining local polynomial reconstruction of the pricing surface $C_\mu$ with residual correction of inverse Fourier features attains the rate uniformly over the noise level. For the lower bound, Theorem~\ref{thm:exact-information} establishes a general limitation of exogenous random linear information, while localized alternatives capture the additional effects of noise and systematic bias.

Our second main result concerns the price $\int g\dd\mu$ of a given payoff $g$. For a convex class $\mathcal C$ of laws, let $\omega_g(\epsilon)$ be the largest price difference between laws in $\mathcal C$ whose pricing surfaces differ by at most $\epsilon$ in the $L^2$ norm of the design. Theorem~\ref{thm:minimax-functional} shows that, at fixed positive noise, the minimax risk is comparable to $\omega_g(\varsigma/\sqrt n)$. Moreover, the modulus has the dual representation
\[
 \omega_g(\epsilon)
 =\inf_w\big[\operatorname{osc}_{\mathcal C}
 \bigl(g-\cA^*(wQ)\bigr)+\epsilon\norm{w}_{L^2(Q)}\big],
\]
where $\cA^*(wQ)$ is a static basket portfolio and $\operatorname{osc}_{\mathcal C}$ measures the remaining price uncertainty over $\mathcal C$. This duality yields a near-minimax affine estimator whose weights are determined by the replicating portfolio.

Our third main result is that long-only quotes identify the joint law, yet the exchange option $(x_1-x_2)_+$ has minimax risk of exact order $(\log n)^{-(\alpha+3/2)}$, uniformly over bounded noise levels, including exact quotes, as shown in Theorem~\ref{thm:exchange-uniform}. The exchange option is a kink $(\be\cdot\bx-k)_+$ in the spread direction, which is unquoted. The same logarithmic rate holds for every kink in an unquoted direction, hence for spread calls and for the best-of and worst-of payoffs, and the law therefore admits no polynomial rate in $W_1$. The reason is that long-only quotes determine the $i$th moment of the spread within the quote error times a factor that is geometric in $i$, whereas the price of a kink depends on moments of every degree. A polynomial claim depends on finitely many moments, and the same quotes therefore estimate every polynomial claim, including every moment of the spread, at the parametric rate $n^{-1/2}$. A kink is Lipschitz stable once its direction is quoted, as shown in Theorem~\ref{thm:log-modulus}, and quoting the spread direction therefore makes all these payoffs Lipschitz stable.

\subsection{Related Literature}
\par\noindent\textbf{Random information and exact quotes.} The term $n^{-(\alpha+1)/d}$ is the error of recovering an $H^\alpha$ density in $W_1$ from $n$ random linear functionals, a question of information-based complexity \cite{novak2008,hinrichs2020}. For point evaluations of functions in the Sobolev space $W^s_p$, i.i.d. points attain the optimal rate of $L_q$-approximation when $q<p$ \cite{krieg2024,krieg2022}. We establish the corresponding result in Corollary~\ref{cor:nu-general-design} for basket quotes, where each observation is the integral of a payoff $(\bu\cdot\bx-k)_+$ against the pricing law and the error is measured in $W_1$. The estimator attains the rate of optimal information from i.i.d. quotes without logarithmic loss, and it does not use the design density. The lower bound in Theorem~\ref{thm:exact-information} holds for any exogenous random information and any randomized estimator, with no linear or convex structure on the class or the loss. Bayes risks under a design-independent prior also give the rates of \cite{pathak2024} for noisy random linear observations under elliptical constraints. For deterministic information, the argument reduces to the Bernstein-width bounds of \cite{novak2008}. Rates that are uniform in the noise level are obtained in \cite{devore2025} for point values of a Besov function at locations chosen by the statistician, which leaves open the case of random locations and indirect observations treated in Theorem~\ref{thm:nu-smooth}.

\par\noindent\textbf{Tomography.} The pricing operator integrates the Radon transform twice in the strike variable, making basket quotes tomographic data with the basket direction specifying the projection angle. The two integrations raise the smoothing order from $(d-1)/2$ to $s=(d+3)/2$, which enters the noise and bias exponents of Theorem~\ref{thm:nu-smooth}. Statistical tomography studies density recovery from random projection data \cite{johnstone1990,korostelev1993,cavalier2000}, with integrated or pointwise squared error. Its rates correspond to the second term in Theorem~\ref{thm:nu-smooth}. The first term concerns exact quotes and has no counterpart in this literature. For exact quotes the bound for Fourier inversion by sample averages in Proposition~\ref{prop:noisy-upper} remains of the fixed-noise order, since the quote locations are random. Related work treats unobserved projection angles \cite{panaretos2009} and inverse problems in density and white-noise models \cite{klemela2010}. Statistical inverse learning studies regularized estimation under source conditions \cite{blanchard2018}. 
Four features distinguish our setting. The unknown is a probability law, and recovery is measured in $W_1$. The admissible laws form a convex set whose diameter equals the largest price uncertainty over Lipschitz claims. The pricing operator admits a negative-Sobolev identity that yields sharp inverse moduli, as in Theorems~\ref{thm:operator-identity} and~\ref{thm:modulus}. Random sampling of quote locations imposes an information limit even when prices are observed exactly. Wasserstein losses appear in direct density estimation through Besov duality \cite{nilesweed2022,weed2019} and in deconvolution \cite{caillerie2011,dedecker2013}. Unlike these settings, the basket transform is neither direct nor a convolution.

\par\noindent\textbf{Single prices and static replication.} For a linear functional over a convex class, the minimax risk is determined by a modulus of continuity and attained by an affine rule. This is the optimal-recovery theory of \cite{micchelli1977} in the statistical form given by \cite{donoho1994}, with further work on convergence rates, nonconvex classes, indirect observations and convex programming in \cite{donoho1991}, \cite{cai2004}, \cite{goldenshluger2000} and \cite{juditsky2009}. We establish this modulus characterization in Theorem~\ref{thm:minimax-functional} for random-design regression with arbitrary bounded continuous test functions and convex weakly closed classes of laws. Under basket quotes the affine rule is a static portfolio of quoted calls and the modulus is a penalized replication problem. Spanning of multi-asset payoffs by basket calls is treated in \cite{bossu2024}, and the one-dimensional recovery problem from single-asset options in \cite{aitsahalia1998,aitsahalia2003}. Static-arbitrage bounds for basket options are obtained by convex optimization in \cite{bertsimas2002} and by super- and subreplicating portfolios in \cite{hobson2005,aspremont2006,tankov2011,neufeld2023}. The penalized replication problem is the statistical form of this duality.

\par\noindent\textbf{Missing directions.} Derivative families that determine the joint law are studied in \cite{molchanov2014}. Long-only quotes identify the law by the homogeneous-polynomial criterion of \cite{cuesta2007}, a refinement of the Cram\'er--Wold theorem \cite{cramer1936,belisle1997}. Weak convergence from positive-measure sets of projections is studied in \cite{cholaquidis2026} under moment assumptions, and compact support permits the direct argument in Theorem~\ref{thm:direction-identification}. Section~\ref{sec:directional-pricing} quantifies the instability of this identification. The missing cone around the spread direction corresponds to the missing wedge of limited-angle tomography \cite{quinto1993,frikel2013,koch2021}, where the unknown image has singularities that the data do not determine stably. In our problem the singularity is the kink of the payoff, and the logarithmic rate in Theorem~\ref{thm:exchange-uniform} is the order of the instability for its price, with a lower bound that holds for exact quotes.

\par\noindent\textbf{Organization.} Section~\ref{sec:model} specifies the observation experiments and establishes identification and consistency. Section~\ref{sec:geometry} develops the operator estimates, and Section~\ref{sec:exact} treats exact sampling and information lower bounds. Section~\ref{sec:noise-uniform} proves the minimax rate uniformly over bounded noise levels. Section~\ref{sec:stable-pricing} relates admissible sets to price intervals and characterizes minimax estimation of a single price, and Section~\ref{sec:directional-pricing} analyzes directional information and the exchange option. Section~\ref{sec:numerical} illustrates the rates numerically. Section~\ref{sec:discussion} concludes and discusses open problems.

\section{Observation Model and Identification}\label{sec:model}
We specify the probability classes and observation experiments, then characterize directional identification and the almost-sure limit of the admissible sets.

\subsection{Quotes and Admissible Laws}
We absorb discounting into the quoted prices and use a single terminal pricing law. Fix $d\ge2$ and $R>0$. Write $B_R=\{\bx\in\R^d:\abs{\bx}\le R\}$ and $\Z=S^{d-1}\times[-R,R]$. For $\mu\in\Pp(B_R)$, define
\begin{equation}\label{eq:call}
 C_\mu(\bu,k)=\int_{B_R}(\bu\cdot\bx-k)_+\,\mu(\mathrm d\bx).
\end{equation}
We observe $n$ independent quotes $(Z_i,V_i)$, where
\begin{equation}\label{eq:exact-model}
 Z_i=(U_i,K_i)\overset{\mathrm{i.i.d.}}{\sim} Q \in \Pp(\Z),
 \qquad V_i=C_\mu(Z_i),\qquad i=1,\ldots,n.
\end{equation}

We use the product metric
$
\rho\left((\bu,k),(\bv,l)\right)
   =\abs{\bu-\bv}+\frac{\abs{k-l}}{R}.
$
For every $\mu,\nu\in\Pp(B_R)$,
\begin{equation}\label{eq:lipschitz}
 0\le C_\mu\le2R,\qquad
 \Lip_\rho(C_\mu)\le R,\qquad
 \norm{C_\mu-C_\nu}_\infty\le W_1(\mu,\nu).
\end{equation}
For a deterministic design $z=(z_1,\ldots,z_n)$ and a feasible price vector $v=(v_1,\ldots,v_n)$, define \begin{align} \K_n(z,v) &=\left\{\nu\in\Pp(B_R): C_\nu(z_i)=v_i,\ 1\le i\le n\right\}, \label{eq:fiber}\\ \K_n(\mu) &=\left\{\nu\in\Pp(B_R): C_\nu(Z_i)=C_\mu(Z_i),\ 1\le i\le n\right\}. \label{eq:random_fiber} \end{align} We call $\K_n(z,v)$ the quote fiber of the price vector $v$, and $\K_n(\mu)$ the random quote fiber of $\mu$. Both sets are convex and compact in the $W_1$ topology whenever they are nonempty.

\subsection{Regularity Classes and Statistical Risk}

Let $\sigma$ be normalized surface measure on $S^{d-1}$, write $\omega_{d-1}=\abs{S^{d-1}}$, and define
\begin{equation}\label{eq:q0}
 Q_0(\mathrm d\bu\,\mathrm dk)=\sigma(\mathrm d\bu)\frac{\mathrm dk}{2R},
 \qquad s=\frac{d+3}{2}.
\end{equation}
For a signed measure $\eta$ supported on $B_R$, let $\cA\eta=C_\eta$ denote the linear extension of \eqref{eq:call}. We use the Fourier convention
\[
 \widehat\eta(\bxi)=\int e^{-\mathrm i\bxi\cdot\bx}\,\eta(\mathrm d\bx),
 \qquad
 \norm{\eta}_{\dot H^{-s}}^2=(2\pi)^{-d}
    \int_{\R^d}\abs{\widehat\eta(\bxi)}^2\abs{\bxi}^{-2s}\dd\bxi.
\]
We also use the inhomogeneous Sobolev norm $\displaystyle \norm{f}_{H^\beta(\R^d)}^2 =(2\pi)^{-d}\int_{\R^d} (1+\abs{\bxi}^2)^\beta \abs{\widehat f(\bxi)}^2\dd\bxi$, $\beta \in \R$. The observation model is
\begin{equation}\label{eq:nu-model}
 Z_i\overset{\mathrm{i.i.d.}}{\sim}Q_0,\qquad
 Y_i=C_\mu(Z_i)+\varsigma\xi_i,\qquad
 0\le\varsigma\le\overline\varsigma<\infty,
\end{equation}
with errors $\xi_i$ that are i.i.d., independent of the design, with $\E\xi_1=0$ and $\E\xi_1^2=1$. Positive-noise lower bounds assume an absolutely continuous error density $p_\xi$ with finite Fisher information $I_\xi=\int(p_\xi')^2/p_\xi$. Gaussian errors have $I_\xi=1$.
For the recovery results, fix $0<r<R$ and assume that the unknown law $\mu$ is supported in $B_r$. The estimator may use the noise level $\varsigma$. The gap between $r$ and $R$ permits bounded extension in the strike coordinate. We define
\begin{equation}\label{eq:smooth-class}
 \F_\alpha(L)=\left\{f\ge0:\int f\dd\bx=1, \quad
       \supp f\subseteq B_r,\quad
       \norm{f}_{H^\alpha(\R^d)}\le L\right\},\quad \alpha\ge0.
\end{equation}
The zero extension of $f$ outside $B_r$ belongs to $H^\alpha(\R^d)$.
For a class $\mathcal C\subseteq\Pp(B_r)$, define
\begin{equation}\label{eq:minimax-risk}
 \mathfrak R_n(\mathcal C,\varsigma)
 =\inf_{\widehat\mu}\sup_{\mu\in\mathcal C}
        \E_{\mu,\varsigma}W_1(\widehat\mu,\mu),
\end{equation}
where the infimum is over measurable estimators based on the observations $(Z_i,Y_i)_{i=1}^n$ and taking values in $\Pp(B_r)$.

\subsection{Directional Identification and Set Consistency}\label{sec:directions}
\begin{definition}
The design distribution $Q$ identifies $\Pp(B_R)$ if equality of $C_\mu$ and $C_\nu$ on $\supp Q$ implies $\mu=\nu$ for all $\mu,\nu\in\Pp(B_R)$.
\end{definition}
Full support on $\Z$ is sufficient, since equality of call functions for all strikes determines the projected law in each direction. Equality in all directions then determines the characteristic function on $\R^d$.

Let $E\subseteq S^{d-1}$ be the closed support of the direction distribution, and assume the strike distribution has full support on $[-R,R]$. For a product design, $\supp Q=E\times[-R,R]$.

\begin{theorem}[Directional identification]\label{thm:direction-identification}
The following statements are equivalent.
\begin{enumerate}[label=(\roman*)]
\item The pricing surfaces on $E\times[-R,R]$ determine every law in $\Pp(B_R)$ uniquely.
\item No nonzero homogeneous real polynomial on $\R^d$ vanishes on every direction in $E$.
\end{enumerate}
If (ii) fails, there are distinct smooth probability densities supported strictly inside $B_R$ with identical call prices for every direction in $E$ and every strike. Every nonempty relatively open set of directions satisfies (ii), a finite set of directions fails (ii), and in dimension two every infinite set of directions satisfies (ii).
\end{theorem}

\begin{theorem}[Set consistency]\label{thm:set-consistency}
For every $\mu\in\Pp(B_R)$, almost surely the decreasing random compact sets $\K_n(\mu)$ converge in Hausdorff distance induced by $W_1$ to
\begin{equation}\label{eq:identified-set}
 \K_\infty(\mu;Q)=\left\{\nu\in\Pp(B_R):
       C_\nu=C_\mu\ \text{on }\supp Q\right\}.
\end{equation}
If $Q$ identifies $\Pp(B_R)$, then $\sup_{\nu\in\K_n(\mu)}W_1(\nu,\mu)\longrightarrow0$ almost surely. Hence, every feasible selection is strongly consistent.
\end{theorem}

\section{Geometry of the Pricing Operator}\label{sec:geometry}

How close must two probability laws be if their basket-call prices are close? We address this question through the pricing operator $\cA$ defined in Section~\ref{sec:model}, for which $\cA(\mu-\nu)=C_\mu-C_\nu$. The results below relate pricing discrepancies to Sobolev norms and Wasserstein distance, and provide representations of smooth functionals used in the estimation procedures. We first consider differences of laws with the same mean. Such differences have zero total mass and zero first moments, and the following theorem expresses their $L^2(Q_0)$ pricing discrepancy exactly in Fourier and spatial form.

\begin{theorem}[Exact operator identity]\label{thm:operator-identity}
Let $\eta$ be a finite signed measure supported on $B_R$ with
$\eta(B_R)=0$ and $\int\bx\,\eta(\mathrm d\bx)=0$. Then
\begin{equation}\label{eq:fourier-identity}
 \norm{\cA\eta}_{L^2(Q_0)}^2
 =\frac{1}{2\pi R\omega_{d-1}}
     \int_{\R^d}\abs{\widehat\eta(\bxi)}^2
                  \abs{\bxi}^{-d-3}\dd\bxi
 =\frac{(2\pi)^{d-1}}{R\omega_{d-1}}
     \norm{\eta}_{\dot H^{-s}}^2.
\end{equation}
There is also the spatial identity
\begin{equation}\label{eq:cubic-identity}
 \norm{\cA\eta}_{L^2(Q_0)}^2
 =\frac{m_{3,d}}{24R}
   \iint\abs{\bx-\by}^3\,\eta(\mathrm d\bx)\eta(\mathrm d\by),
 \quad
 m_{3,d}=\frac{\Gamma(d/2)}
                 {\sqrt\pi\,\Gamma((d+3)/2)}.
\end{equation}
\end{theorem}

The equal-mean restriction is not needed for the stability bounds below: Lemma~\ref{lem:mean} controls the difference of means by the pricing discrepancy, allowing a smooth correction before applying Theorem~\ref{thm:operator-identity}.

\begin{theorem}[Inverse moduli]\label{thm:modulus}
There is $C=C(d,R)>0$ so that for all
$\mu,\nu\in\Pp(B_R)$,
\begin{equation}\label{eq:universal-modulus}
 W_1(\mu,\nu)\le C
 \norm{C_\mu-C_\nu}_{L^2(Q_0)}^{\,2/(d+3)}.
\end{equation}
If $\mu=f\dd\bx$, $\nu=g\dd\bx$, and
$\norm{f}_{H^\alpha},\norm{g}_{H^\alpha}\le L$
for $\alpha\ge0$, then, with
$D=\norm{C_\mu-C_\nu}_{L^2(Q_0)}$,
\begin{equation}\label{eq:smooth-modulus}
 W_1(\mu,\nu)
 \le C_{d,R,\alpha} \big[
 (1+L)^{(s-1)/(\alpha+s)}
 D^{(\alpha+1)/(\alpha+s)}+D \big].
\end{equation}
The exponent in \eqref{eq:universal-modulus} is optimal.
For every $\alpha\ge0$ and every sufficiently large $L$,
the exponent $(\alpha+1)/(\alpha+s)$ in
\eqref{eq:smooth-modulus} is also optimal over the indicated
density class.
\end{theorem}

Recovery from sampled quotes also requires Sobolev bounds on the pricing surface. We use a fixed Hilbert norm on the standard $L^2$ Sobolev space $H^t(\Z)$.

\begin{lemma}[Forward and inverse Sobolev bounds]
\label{lem:sobolev-map}
Let $\eta\in H^\alpha(\R^d)$ be supported on $B_r$ and satisfy
$\int\eta\dd\bx=0$. Then
\begin{equation}\label{eq:forward-sobolev}
 \norm{\cA\eta}_{H^{\alpha+s}(\Z)}
 \le C_{d,r,R,\alpha}\norm{\eta}_{H^\alpha(\R^d)}.
\end{equation}
For any two probability densities $f,g$ supported on $B_r$
for which the right side is finite,
\begin{equation}\label{eq:inverse-sobolev}
 W_1(f\dd\bx,g\dd\bx)
 \le C_{d,r,R}\norm{\cA(f-g)}_{H^{s-1}(\Z)}.
\end{equation}
\end{lemma}

Finally, smooth linear functionals can be written as weighted integrals of the pricing surface, plus a constant. The $L^2(Q_0)$ norm of the weight function bounds the estimator variance.

\begin{lemma}[Inverse features]\label{lem:dual-features}
For every real $\psi\in C_c^\infty(\R^d)$ there are a constant
$a_\psi$ and a function $\widetilde G_\psi\in L^2(Q_0)$ such
that, for all $\mu\in\Pp(B_r)$,
\begin{equation}\label{eq:dual-functional}
 \int\psi\dd\mu
 =a_\psi+
   \int_\Z\widetilde G_\psi(\bu,k)C_\mu(\bu,k)\dd Q_0,
 \qquad
 \abs{a_\psi}+ \|\widetilde G_\psi\|_{L^2(Q_0)}
 \le C\norm{\psi}_{H^s(\R^d)}.
\end{equation}
The construction extends to complex-valued $\psi$ by linearity.
\end{lemma}

\section{Exact Recovery and Random Information}\label{sec:exact}

How much uncertainty remains about the law after finitely many
exact basket quotes have been observed? We first bound the
Wasserstein diameter of the quote fiber in terms of the coverage
of the observed locations, without imposing smoothness on the
law. We then obtain sharper bounds on Sobolev density classes
and construct an estimator that does not require the class
radius. Finally, we establish lower bounds showing the limits
of recovery from finitely many exact linear observations.
Extensions of the coverage bounds to stopping times and
adaptive quote arrivals are given in Appendix~\ref{AppendixC}.
Throughout this paper, we set $\diam_{W_1}\varnothing=0$.

\subsection{Recovery Without Smoothness Assumptions}
\label{sec:exact-general}

For $z=(z_i)_{i=1}^n$, define its fill distance
\begin{equation}\label{eq:fill}
 h(z)=\sup_{w\in\Z}\min_{1\le i\le n}\rho(w,z_i),
 \qquad h_n=h(Z_1,\ldots,Z_n).
\end{equation}
$h(z)$ is the largest distance from a possible quote
location to the nearest observed location. As 
$\Z=S^{d-1}\times[-R,R]$ has dimension $d$,
for $0<t<1$, it admits a $t$-net of cardinality at most
$C_dt^{-d}$, and every $\rho$-ball of radius $t$ centered
in $\Z$ has $Q_0$-probability at least $c_dt^d$. 

If two laws agree at the observed locations, the Lipschitz
bound in \eqref{eq:lipschitz} controls the difference of their
prices elsewhere by the fill distance. The inverse modulus
in Theorem~\ref{thm:modulus} then bounds their Wasserstein
distance.

\begin{theorem}[Distribution-free exact recovery bound]
\label{thm:exact-general}
For any deterministic design $z$,
\begin{equation}\label{eq:general-diameter}
 \sup_{v\in\R^n} \diam_{W_1}\K_n(z,v)
 \le C_{d,R}h(z)^{2/(d+3)}.
\end{equation}
If $Q\ge q_*Q_0$ for some $q_*>0$, then with probability
at least $1-\delta$,
\begin{equation}\label{eq:general-random-rate}
 \sup_{v\in\R^n} \diam_{W_1}\K_n(Z,v)
 \le C_{d,R,q_*}
 \Big(\frac{\log(en/\delta)}{n} \Big)^{\frac{2}{d(d+3)}}.
\end{equation}
\end{theorem}

If instead
$\abs{C_\mu(z_i)-v_i}\le\epsilon$ and
$\abs{C_\nu(z_i)-v_i}\le\epsilon$ for every $i$, the same argument gives
\[
 W_1(\mu,\nu)
 \le C_{d,R}(\epsilon+h(z))^{2/(d+3)}.
\]
Proposition~\ref{prop:bounded-errors} below gives the
corresponding bound on Sobolev classes.

The quantitative coverage assumption $Q\ge q_*Q_0$ cannot
be replaced by full support alone if a rate uniform over
design distributions is desired. For an arbitrary
full-support $Q$, the coverage inequality
\eqref{eq:fill-tail} in Appendix~\ref{AppendixC} holds with
$cq_*t^d$ replaced by
$
 \beta_Q(t)=\inf_{w\in\Z}Q(B_\rho(w,t)).
$
This quantity is positive for each $t>0$, by compactness
and a finite covering by smaller balls, but it need not
have a uniform power-law lower bound. The following
example shows the resulting limitation.

\begin{proposition}[No rate uniform over full-support designs]
\label{prop:no-uniform-rate}
Fix $z_0\in\Z$ and $0<\epsilon<1$, and let
$Q_\epsilon=(1-\epsilon)\delta_{z_0}+\epsilon Q_0$,
which has full support. For exact quotes at locations
$Z_i$ i.i.d.\ $Q_\epsilon$ and estimators $\widehat\mu$
based on $(Z_i,V_i)_{i=1}^n$ with values in $\Pp(B_r)$,
\[
 \inf_{\widehat\mu}\sup_{\mu\in\Pp(B_r)}
 \E_{\mu,Q_\epsilon}W_1(\widehat\mu,\mu)
 \ge r(1-\epsilon)^n.
\]
Hence the supremum of this minimax risk over full-support
designs is at least $r$ for every $n$.
\end{proposition}

\subsection{Recovery on Sobolev Classes}
\label{sec:exact-smooth}

Sobolev regularity yields a sharper bound on the uncertainty
remaining after exact quotes. The difference of two feasible
pricing surfaces vanishes at every observed location.
A Sobolev sampling inequality, combined with the forward
and inverse bounds of Lemma~\ref{lem:sobolev-map}, converts
these equalities into a bound of order $h(z)^{\alpha+1}$
on the Wasserstein distance. For $\F_\alpha(L)$ in \eqref{eq:smooth-class}, let
\[
 \K_{n,\alpha,L}(z,v)
 =\left\{f\dd\bx:
 f\in\F_\alpha(L),\
 C_f(z_i)=v_i,\ 1\le i\le n\right\}.
\]

\begin{theorem}[Exact recovery on Sobolev classes]
\label{thm:exact-smooth}
There are constants $C,h_*>0$, depending only on
$d,r,R,\alpha$, such that every design with $h(z)\le h_*$
satisfies
\begin{equation}\label{eq:smooth-exact-diameter}
 \sup_{v\in\R^n}
 \diam_{W_1}\K_{n,\alpha,L}(z,v)
 \le CLh(z)^{\alpha+1}.
\end{equation}
If $Q\ge q_*Q_0$ for some $q_*>0$, then with probability at least $1-\delta$,
\begin{equation}\label{eq:smooth-exact-random}
 \sup_{v\in\R^n}
 \diam_{W_1}\K_{n,\alpha,L}(Z,v)
 \le C_{q_*}L
 \Big(\frac{\log(en/\delta)}{n}\Big)^{\frac{\alpha+1}{d}}.
\end{equation}
\end{theorem}

We also allow each quoted price to be known only up to a
bounded error. For $\epsilon\ge0$, let
\[
 \K^\epsilon_{n,\alpha,L}(z,v)
 =\left\{f\dd\bx:
 f\in\F_\alpha(L),\
 \abs{C_f(z_i)-v_i}\le\epsilon,\ 1\le i\le n\right\}.
\]

\begin{proposition}[Bounded quote errors]
\label{prop:bounded-errors}
Let $\vartheta=(\alpha+1)/(\alpha+s)$.
There are $C,h_*>0$, depending only on $d,r,R,\alpha$,
such that every design with $h(z)\le h_*$ satisfies
\begin{equation}\label{eq:bounded-errors}
 \sup_{v\in\R^n}
 \diam_{W_1}\K^\epsilon_{n,\alpha,L}(z,v)
 \le C\left[
 Lh(z)^{\alpha+1}
 +L^{1-\vartheta}\epsilon^\vartheta
 \right].
\end{equation}
If $Q\ge q_*Q_0$ for some $q_*>0$, then with probability
at least $1-\delta$, simultaneously for all $v$ and
$\epsilon$, the right side can be replaced by
$
 C_{q_*}\big[
 L\big(\frac{\log(en/\delta)}{n}\big)^{(\alpha+1)/d}
 +L^{1-\vartheta}\epsilon^\vartheta\big].
$
\end{proposition}

The two terms in the upper bound represent uncertainty
due to finite sampling and bounded quote errors,
respectively. They are of the same order when
$\epsilon\asymp Lh(z)^{\alpha+s}$; above this scale,
the quote-error term is larger.
Section~\ref{sec:noise-uniform} studies the corresponding
minimax expected risk when bounded systematic errors
are combined with random quote noise.

For a specified smoothness index $\alpha$, define the
minimum-norm interpolant
\begin{equation}\label{eq:minnorm-estimator}
 \widehat f_n
 =\argmin_f\left\{
 \norm{f}_{H^\alpha}^2:
 f\ge0,\ \supp f\subseteq B_r,\
 \int f\dd\bx=1, \ C_f(Z_i)=V_i\ (1\le i\le n)
 \right\}. \nonumber
\end{equation}
If the exact quotes are generated by a density
$f_0\in\F_\alpha(L)$, then minimality gives
$\|\widehat f_n \|_{H^\alpha}
 \le \|f_0\|_{H^\alpha}\le L$.
Hence $\widehat f_n\in\F_\alpha(L)$, although the
optimization requires the smoothness index $\alpha$
but not the radius $L$. Existence, uniqueness, and
measurability are shown in Appendix~\ref{AppendixC}.
On the event in \eqref{eq:smooth-exact-random},
$
 \sup_{f_0\in\F_\alpha(L)}
 W_1(\widehat f_n\dd\bx,f_0\dd\bx)
 \le CL
 \big(\frac{\log(en/\delta)}{n}\big)^{(\alpha+1)/d},
$
where, for each $f_0$, the estimator $\widehat f_n$
is computed from the exact quotes generated by $f_0$. 

For i.i.d. locations with $Q\ge q_*Q_0$, the proof of
Theorem~\ref{thm:random-exact-lower} gives the
worst-case expected error bound
$
 \sup_{f_0\in\F_\alpha(L)}
 \E_{f_0}W_1(\widehat f_n\dd\bx,f_0\dd\bx)
 \le CL
 \big(\frac{\log(en)}{n}\big)^{(\alpha+1)/d}.
$
Corollary~\ref{cor:random-exact-lower} below shows that
this estimator is minimax optimal up to a factor
$(\log(en))^{(\alpha+1)/d}$.

\subsection{Lower Bounds for Exact Linear Information}
\label{sec:exact-lower-bounds}

We now ask whether a different choice of quotes or
estimator could improve the recovery rate.
The next theorem shows that, for any fixed collection
of $n$ continuous payoffs, two densities in $\F_\alpha(L)$
have identical prices but remain separated in $W_1$
by at least $cn^{-(\alpha+1)/d}$.
Thus the polynomial exponent in
\eqref{eq:smooth-exact-random} cannot be improved
for uniform recovery from a fixed design.

\begin{theorem}[Deterministic uniform lower bound]
\label{thm:exact-lower}
For every $\alpha\ge0$, there is $L_0<\infty$ such that
for each $L>L_0$ there is $c>0$ with the following
property. For every $n\ge1$ and every collection of
continuous payoffs $\phi_1,\ldots,\phi_n$ on $B_r$,
there are $f^+,f^-\in\F_\alpha(L)$ satisfying
\begin{equation}\label{eq:exact-lower}
 \int\phi_i f^+\dd\bx
 =\int\phi_i f^-\dd\bx\quad(1\le i\le n),
 \qquad
 W_1(f^+\dd\bx,f^-\dd\bx)
 \ge cn^{-(\alpha+1)/d}.
\end{equation}
Hence for a deterministic quote design, the supremum of fiber diameters is at least $n^{-(\alpha+1)/d}$.
\end{theorem}

The densities $f^+,f^-$ in Theorem~\ref{thm:exact-lower} may depend on the measurements, whereas minimax expected risk fixes the density before averaging over the design. To obtain a lower bound for random measurements, we use a design-independent prior whose average risk bounds the worst-case risk below.  Unlike \eqref{eq:smooth-exact-random}, which holds for all feasible price vectors on a common design event, this expected-risk criterion requires no logarithmic loss for exact quotes; see Theorem~\ref{thm:nu-smooth}. The next theorem supplies the corresponding information lower bound.

\begin{theorem}[Lower bound for exact linear observations]
\label{thm:exact-information}
Let $\mathsf d$ be a pseudometric on a set $\mathcal X$,
and let
$\{f_{\btheta}:\btheta\in[-1,1]^M\}\subseteq\mathcal X$
be continuous in $\btheta$ for $\mathsf d$ and satisfy
\[
 \mathsf d(f_{\btheta},f_{\btheta'})
 \ge\delta\norm{\btheta-\btheta'}_1,
 \qquad
 \btheta,\btheta'\in[-1,1]^M,
\]
for some $\delta>0$.
Let the data consist of $n$ observations
$V_i=\mathsf L_i(f)+a_i$, where the linear functionals
$\mathsf L_1,\ldots,\mathsf L_n$ and the shift
$a\in\R^n$ have a joint distribution that does not
depend on $f$, and
$\btheta\mapsto\mathsf L_i(f_{\btheta})$
is affine for every realization.
Then every measurable estimator $\widehat f\in\mathcal X$ based on $(V_i,\mathsf L_i,a_i)_{i=1}^n$ and independent randomization satisfies
\begin{equation}\label{eq:exact-information}
 \sup_{\btheta\in[-1,1]^M}
 \E_{f_{\btheta}}\,
 \mathsf d(\widehat f,f_{\btheta})
 \ge\frac{\delta(M-n)}{4\pi^2}.
\end{equation}
\end{theorem}
The proof uses a fixed product prior to convert coefficient uncertainty left by a rank $\le n$ observation matrix into expected loss. This Bakhvalov-type argument is related to Bernstein-width bounds \cite{novak2008} and Bayesian bounds for random linear recovery \cite{pathak2024}. It allows random measurements and randomized estimators without requiring $\mathcal X$ to be linear or convex. Theorem~\ref{thm:exact-information}, applied to the family in Lemma~\ref{lem:bumps}, gives a lower bound for random basket quotes.

For a design distribution $\mathsf Q_n$ on $\Z^n$ independent of the unknown law, let
\begin{equation}\label{eq:exact-minimax-definition}
 \mathfrak R_n^{\mathrm{ex}}(\mathsf Q_n,\alpha,L)
 =\inf_{\widehat\mu}\sup_{f\in\F_\alpha(L)}
 \E_{f,\mathsf Q_n}W_1(\widehat\mu,f\dd\bx),
\end{equation}
where the observations are $Z\sim\mathsf Q_n$ and
$V_i=C_f(Z_i)$, and estimators take values in $\Pp(B_r)$.

\begin{corollary}[Exact random-design lower bound]
\label{cor:random-exact-lower}
Fix $\alpha\ge0$ and $L>L_0$ as in
Theorem~\ref{thm:exact-lower}.
There is $c>0$, independent of $\mathsf Q_n$, such that
\begin{equation}\label{eq:random-exact-lower}
 \mathfrak R_n^{\mathrm{ex}}(\mathsf Q_n,\alpha,L)
 \ge cn^{-(\alpha+1)/d}.
\end{equation}
The bound  remains valid for any $n$ exogenously chosen continuous linear payoffs.
\end{corollary}

Thus the rate $n^{-(\alpha+1)/d}$ is a lower limit
for worst-case expected recovery error even when
the observations are exact and the measurement
distribution can be chosen freely, provided that
it is independent of the unknown law.
Section~\ref{sec:noise-uniform} shows that i.i.d.\
basket quotes with adequate coverage attain this
rate without a logarithmic loss.

\section{Minimax Recovery Across Noise Levels}\label{sec:noise-uniform}
How does recovery improve as quotes become more accurate, and what limits remain when their number increases? We distinguish independent quote noise from a bounded systematic bias. The former can be reduced by averaging, whereas the latter need not disappear even with arbitrarily many observations.

Specifically, let $e:\Z\to\R$ be a Borel function with $\sup_{z\in\Z}|e(z)|\le\epsilon$, and observe
\begin{equation}\label{eq:systematic-model}
 Y_i=C_\mu(Z_i)+e(Z_i)+\varsigma\xi_i,\qquad i=1,\ldots,n,
\end{equation}
where $Z_i$ are i.i.d.\ with law $Q_0$ and the errors $\xi_i$ satisfy the assumptions of \eqref{eq:nu-model}. The function $e$ specifies the systematic bias at each direction--strike pair. It is fixed before sampling and may depend on the unknown law, but not on the realized locations or noise. No smoothness or mean-zero condition is imposed on $e$.
For $\mathcal C\subseteq\Pp(B_r)$, define
\[
 \mathfrak R_n(\mathcal C,\varsigma,\epsilon)
 =\inf_{\widehat\mu}\sup_{\mu\in\mathcal C}
   \sup_{|e|\le\epsilon}
       \E_{\mu,e,\varsigma}W_1(\widehat\mu,\mu),
\]
where estimators take values in $\Pp(B_r)$. As before, we identify a density class with its associated probability laws. In particular,
$\mathfrak R_n(\mathcal C,\varsigma,0)=\mathfrak R_n(\mathcal C,\varsigma)$.
Throughout this section, write
\[
 t=\alpha+s,\qquad D_\alpha=2t+d=2\alpha+2d+3,
 \qquad \vartheta=(\alpha+1)/(\alpha+s).
\]

\begin{theorem}[Sharp rate uniformly in the noise level]\label{thm:nu-smooth}
Fix $\alpha\ge0$, $L>L_0$ as in Theorem~\ref{thm:exact-lower}, and finite nonnegative bounds $\overline\varsigma$ and $\overline\epsilon$. There are constants $0<c<C<\infty$, independent of $n,\varsigma,\epsilon$, such that, for every $n\ge1$, $0\le\varsigma\le\overline\varsigma$, and $0\le\epsilon\le\overline\epsilon$,
\begin{equation}\label{eq:nu-smooth-rate}
 c\Big[n^{-\frac{\alpha+1}{d}}
       +\Big(\frac{\varsigma^2}{n}\Big)^{\frac{\alpha+1}{D_\alpha}}+\epsilon^\vartheta\Big]
 \le \mathfrak R_n(\F_\alpha(L),\varsigma,\epsilon)
 \le C\Big[n^{-\frac{\alpha+1}{d}}
       + \Big(\frac{\varsigma^2}{n} \Big)^{\frac{\alpha+1}{D_\alpha}}+\epsilon^\vartheta\Big].
\end{equation}
The upper bound holds for independent, centered, unit-variance errors, without a density or tail assumption. 
The constants depend on $d,r,R,\alpha,L,\overline\varsigma,\overline\epsilon$, and, for the lower bound, on $I_\xi$.
\end{theorem}

The three terms describe distinct limitations. Even exact quotes leave uncertainty of order $n^{-(\alpha+1)/d}$ because only finitely many linear measurements are observed. Independent noise adds the second term, which decreases with both the sample size and the noise level. The last term is the worst-case uncertainty caused by systematic bias. Since the sum is comparable to its largest term, the largest of these three contributions determines the recovery rate.

\begin{corollary}[Transition scales]\label{cor:nu-transition}
Under the assumptions of Theorem~\ref{thm:nu-smooth}, the exact-information term and the independent-noise term have the same order when
\(
 \varsigma\asymp n^{-\frac{2\alpha+d+3}{2d}}.
\)
The exact-information term and the systematic-bias term have the same order when
\(
 \epsilon\asymp n^{-\frac{2\alpha+d+3}{2d}}.
\)
If $\epsilon=0$ and $\varsigma_n=n^{-\gamma}$ with $\gamma\ge0$, then the minimax risk has order $n^{-\beta}$, where
\(
 \beta=\min\big\{\frac{\alpha+1}{d},
           \frac{(2\gamma+1)(\alpha+1)}{2\alpha+2d+3}\big\}.
\)
For example, when $d=2$ and $\alpha=1$, the risk is comparable to
$n^{-1}+(\varsigma^2/n)^{2/9}+\epsilon^{4/7}$, and both transition scales above are $n^{-7/4}$.
\end{corollary}

\begin{remark}[An irreducible error from systematic bias]\label{rem:systematic}
For a fixed $\epsilon>0$, increasing the number of quotes cannot make the worst-case recovery error tend to zero. Two different laws can have pricing surfaces within $2\epsilon$ of one another, so that admissible biases make their observed mean prices identical. The term $\epsilon^\vartheta$ quantifies the resulting ambiguity on $\F_\alpha(L)$. 
\end{remark}

\subsection{Direct Fourier Reconstruction}\label{sec:noise}
We first describe a direct estimator for unbiased quotes, that is, $e=0$. This provides a baseline and introduces the Fourier reconstruction used below. Lemma~\ref{lem:dual-features} expresses each smooth moment as an integral of the pricing surface. Replacing this integral by a sample average gives, for $\psi\in C_c^\infty(\R^d)$,
\begin{equation}\label{eq:moment-estimator}
 \widehat\theta_\psi=a_\psi+
       \frac1n\sum_{i=1}^n\widetilde G_\psi(Z_i)Y_i.
\end{equation}
This estimator is unbiased for $\int\psi\dd\mu$, and
\begin{equation}\label{eq:moment-variance}
 \E_{\mu,\varsigma}\left|\widehat\theta_\psi-\int\psi\dd\mu\right|^2
 \le\frac{4R^2+\varsigma^2}{n}
       \norm{\widetilde G_\psi}_{L^2(Q_0)}^2
 \le\frac{C}{n}\norm{\psi}_{H^s}^2.
\end{equation}
Indeed, $\E(Y_i^2\mid Z_i)=C_\mu(Z_i)^2+\varsigma^2\le4R^2+\varsigma^2$.

To reconstruct the law from these moments, work on the torus
$\mathbb T_T^d=\R^d/(T\mathbb Z^d)$ with $T>4R$, represented by the cube $[-T/2,T/2)^d$ with opposite faces identified. The ball $B_R$ lies strictly inside this cube. Moreover, torus geodesic distance agrees with Euclidean distance between points of $B_r$.
Choose a real $\chi\in C_c^\infty((-T/2,T/2)^d)$ with $\chi=1$ on $B_r$, and set
\(
 \psi_{\bk}(\bx)=\chi(\bx)e^{-2\pi\mathrm i\bk\cdot\bx/T},\)
\( \bk\in\mathbb Z^d.
\)
Then $\norm{\psi_{\bk}}_{H^s}\le C(1+|\bk|)^s$.
Write $a_{\bk}=a_{\psi_{\bk}}$ and
$\widetilde G_{\bk}=\widetilde G_{\psi_{\bk}}$.
For $\bk\ne0$, let $\widehat\theta_{\bk}$ be \eqref{eq:moment-estimator} with $\psi=\psi_{\bk}$; it estimates
\(
 \theta_{\bk}(\mu)=\int e^{-2\pi\mathrm i\bk\cdot\bx/T}\dd\mu(\bx).
\)
Set $\widehat\theta_{\boldsymbol0}=1$. Choose the features consistently with complex conjugation, so that
$\widehat\theta_{-\bk}=\overline{\widehat\theta_{\bk}}$.

Fix a real even $w\in C_c^\infty(\R^d)$ equal to one on the unit ball and zero outside the ball of radius two. For a frequency cutoff $M\ge1$, define
\begin{equation}\label{eq:fourier-estimator}
 \widetilde f_M(\bx)=\frac1{T^d}\sum_{\bk\in\mathbb Z^d}
       w(\bk/M)\widehat\theta_{\bk}e^{2\pi\mathrm i\bk\cdot\bx/T}.
\end{equation}
$\widetilde f_M$ is finite and real. Orthogonality of the Fourier modes gives
\(
 \int_{\mathbb T_T^d}\widetilde f_M(\bx)\dd\bx
       =w(0)\widehat\theta_{\boldsymbol0}=1.
\)
However, $\widetilde f_M$ need not be nonnegative, nor need it be supported in $B_r$. It therefore defines a signed measure of mass one, not necessarily a probability law. We enforce the probability and support constraints by projection. For a finite signed measure $\zeta$ of mass zero on $\mathbb T_T^d$, put
\[
 \norm{\zeta}_{\Lip^*}
 =\sup_{\Lip(F)\le1}\left|\int F\dd\zeta\right|,
\]
where Lipschitz constants use torus geodesic distance. For probabilities supported on $B_r$, this norm of their difference equals their Euclidean $W_1$ distance. Choose $\widehat\mu_M\in\Pp(B_r)$ satisfying
\begin{equation}\label{eq:kr-projection}
 \norm{\widehat\mu_M-\widetilde f_M\dd\bx}_{\Lip^*}
 \le\inf_{\nu\in\Pp(B_r)}
       \norm{\nu-\widetilde f_M\dd\bx}_{\Lip^*}+\frac1n.
\end{equation}

\begin{proposition}[Direct reconstruction bounds]\label{prop:noisy-upper}
Under \eqref{eq:nu-model}, for every $M\ge1$ and $n\ge1$, the estimator in \eqref{eq:kr-projection} satisfies, with constants uniform over $0\le\varsigma\le\overline\varsigma$,
\begin{align}\label{eq:smooth-bias-variance}
 \sup_{f\in\F_\alpha(L)}\E_{f,\varsigma}W_1(\widehat\mu_M,f\dd\bx)
 &\le C\left[LM^{-\alpha-1}+n^{-1/2}M^{d+1/2}\right]+\frac1n,\\
 \sup_{\mu\in\Pp(B_r)}\E_{\mu,\varsigma}W_1(\widehat\mu_M,\mu)
 &\le C\left[M^{-1}+n^{-1/2}M^{d+1/2}\right]+\frac1n.
 \label{eq:general-bias-variance}
\end{align}
\end{proposition}

Taking $M\asymp n^{1/D_\alpha}$ in \eqref{eq:smooth-bias-variance} gives the rate $n^{-(\alpha+1)/D_\alpha}$. The next result shows that this rate is optimal when the noise level is fixed and positive and there is no systematic bias.
\begin{proposition}[Lower bound at fixed positive noise]\label{prop:noisy-lower}
Fix $\varsigma>0$, $\alpha\ge0$, and $L>L_0$ as in Theorem~\ref{thm:exact-lower}. Under the finite-Fisher-information assumption on the error density, there is $c>0$ such that
\(
 \mathfrak R_n(\F_\alpha(L),\varsigma)
       \ge c n^{-(\alpha+1)/(2\alpha+2d+3)}.
\)
\end{proposition}

\subsection{Estimating the Pricing Surface Locally}
The direct estimator does not fully exploit highly accurate quotes: even when $\varsigma=0$, its sample average fluctuates because the quote locations are random. To reduce this variation, we first estimate the entire pricing surface from part of the data. We call this preliminary estimate a \emph{pilot}; the remaining observations will correct its error rather than estimate the original surface from scratch. The pilot fits polynomials to nearby direction--strike observations. Its bandwidth $h$ is the nominal side length of the local fitting regions: smaller $h$ resolves finer variation, whereas larger $h$ pools more observations. Where a region contains too few suitably placed observations, the construction enlarges it until a stable fit is possible. The Appendix gives the precise construction.

\begin{lemma}[Pricing-surface pilots]\label{lem:nu-pilot}
Under \eqref{eq:nu-model}, for $m\ge1$ and $\varsigma\ge0$, there are measurable estimators $\widehat C_m$ taking values in $[0,2R]$ with the following properties. For $\alpha\ge0$ and $t=\alpha+s$, a local polynomial construction of fixed coordinatewise degree at least $\lceil t\rceil$ satisfies
\begin{equation}\label{eq:nu-smooth-pilot}
 \sup_{f\in\F_\alpha(L)}
 \E_{f,\varsigma}\norm{\widehat C_m-C_f}_{L^2(Q_0)}^2
 \le C\left[L^2m^{-2t/d}+\varsigma^2\right].
\end{equation}
A construction of any fixed coordinatewise degree at least one satisfies
\begin{equation}\label{eq:nu-universal-pilot}
 \sup_{\mu\in\Pp(B_r)}
 \E_{\mu,\varsigma}\norm{\widehat C_m-C_\mu}_{L^2(Q_0)}^2
 \le C\left[m^{-3/d}+\varsigma^2\right].
\end{equation}
More generally, for any prescribed $h\in[m^{-1/d},1]$, the corresponding estimators $\widehat C_{m,h}$ satisfy
\begin{equation}\label{eq:nu-bandwidth-pilot}
\begin{split}
 \sup_{f\in\F_\alpha(L)}\E_{f,\varsigma}
       \norm{\widehat C_{m,h}-C_f}_{L^2(Q_0)}^2
 &\le C\Big[L^2h^{2t}+\frac{\varsigma^2}{mh^d}\Big],\\
 \sup_{\mu\in\Pp(B_r)}\E_{\mu,\varsigma}
       \norm{\widehat C_{m,h}-C_\mu}_{L^2(Q_0)}^2
 &\le C\Big[h^3+\frac{\varsigma^2}{mh^d}\Big].
\end{split}
\end{equation}
The two bounds concern the respective polynomial constructions above. Before pointwise truncation to $[0,2R]$, each estimator is linear in $(Y_i)_{i=1}^m$ conditional on $(Z_i)_{i=1}^m$. All assertions remain valid for i.i.d. locations with law $Q\ge q_*Q_0$, with constants depending on $q_*>0$; the estimators do not require knowledge of $Q$ or a density of $Q$. Under \eqref{eq:systematic-model}, the same bounds hold with an additional term $C\epsilon^2$ on each right-hand side.
\end{lemma}

In \eqref{eq:nu-bandwidth-pilot}, the first term is approximation error and the second is the contribution of independent noise. The choice $h\asymp m^{-1/d}$ gives \eqref{eq:nu-smooth-pilot} and \eqref{eq:nu-universal-pilot}. At this fine resolution the pilot alone need not average away the noise; its role is to make the subsequent residuals small when quotes are accurate.

\subsection{Residual Correction and Probability-Valued Recovery}
For $e=0$, the variance of a direct moment estimate has the exact decomposition
\[
 \operatorname{Var}(\widehat\theta_\psi)
 =\frac1n\operatorname{Var}\!\left(\widetilde G_\psi(Z)C_\mu(Z)\right)
  +\frac{\varsigma^2}{n}\norm{\widetilde G_\psi}_{L^2(Q_0)}^2,
\]
The first term is caused by random quote locations and remains nonzero for exact quotes. Thus the direct bound in Proposition~\ref{prop:noisy-upper} does not improve to the exact-information rate merely by setting $\varsigma=0$. Residual correction removes most of this location-induced variation. For $n \ge 2$, we split the sample into groups of sizes $n_0=\lfloor n/2\rfloor$ and $n_1=n-n_0$, and construct a pilot $\widehat C$ from the first group. For $\bk\ne0$, define
\begin{equation}\label{eq:nu-cv-moment}
 \widehat\theta_{\bk}^{\mathrm{cv}}
 ={}a_{\bk}+\int_\Z\widetilde G_{\bk}(z)\widehat C(z)\dd Q_0(z)
 +\frac1{n_1}\sum_{i=n_0+1}^n
       \widetilde G_{\bk}(Z_i)\left(Y_i-\widehat C(Z_i)\right),
 \quad \widehat\theta_{\boldsymbol0}^{\mathrm{cv}}=1.
\end{equation}
The integral accounts for the estimated surface using the known design $Q_0$; only the residual is averaged over the second group. If the pilot equaled $C_\mu$ and the quotes were exact and unbiased, these residuals would vanish. More generally, conditional on the first group, the corrected moment has mean
\(
 \theta_{\bk}(\mu)+\int_\Z\widetilde G_{\bk}(z)e(z)\dd Q_0(z).
\)
In particular, pilot approximation does not introduce bias into the corrected moment when $e=0$.

Use \eqref{eq:fourier-estimator} with these corrected moments to obtain $\widetilde f_M^{\mathrm{cv}}$, and choose $\widehat\mu_M^{\mathrm{cv}}\in\Pp(B_r)$ satisfying
\begin{equation}\label{eq:nu-projection}
 \norm{\widehat\mu_M^{\mathrm{cv}}-\widetilde f_M^{\mathrm{cv}}\dd\bx}_{\Lip^*}
 \le\inf_{\nu\in\Pp(B_r)}
       \norm{\nu-\widetilde f_M^{\mathrm{cv}}\dd\bx}_{\Lip^*}+e^{-n}.
\end{equation}
The corrected construction preserves complex conjugation and total mass. A measurable projection is obtained as in Proposition~\ref{prop:noisy-upper}. The smaller projection tolerance $e^{-n}$ avoids limiting the attainable rate when $\alpha$ is large.

\begin{proposition}[Pilot-to-law risk bound]\label{prop:nu-transfer}
Let the data follow \eqref{eq:systematic-model}. Suppose $\widehat C$ is based only on the first group and
$\E\|\widehat C-C_\mu\|_{L^2(Q_0)}^2\le v_n(\mu)$.
For every $M,n$ and $f\in\F_\alpha(L)$,
\begin{equation}\label{eq:nu-transfer-smooth}
\E_{f,e}W_1(\widehat\mu_M^{\mathrm{cv}},f\dd\bx)
 \le C\Big[L M^{-\alpha-1}+\epsilon M^{s-1}
       +\sqrt{\frac{\varsigma^2+\epsilon^2+v_n(f)}{n}}\,M^{d+\frac{1}{2}}\Big]+e^{-n}.
\end{equation}
For every $\mu\in\Pp(B_r)$,
\begin{equation}\label{eq:nu-transfer-all}
 \E_{\mu,e}W_1(\widehat\mu_M^{\mathrm{cv}},\mu)
 \le C\Big[M^{-1}+\epsilon M^{s-1}
       +\sqrt{\frac{\varsigma^2+\epsilon^2+v_n(\mu)}{n}}\,M^{d+\frac{1}{2}}\Big]+e^{-n}.
\end{equation}
\end{proposition}

The three contributions in \eqref{eq:nu-transfer-smooth} are Fourier approximation, amplified systematic bias, and random fluctuation. Unlike the direct bound, the last contribution is controlled by the pilot error and the noise level rather than by the size of the pricing surface itself.

\begin{proof}[Proof sketch of Theorem~\ref{thm:nu-smooth}]
Lemma~\ref{lem:nu-pilot}, with $m=n_0$ and $h\asymp n^{-1/d}$, gives
$v_n(f)\le C[L^2n^{-2t/d}+\varsigma^2+\epsilon^2]$ uniformly over the class and admissible biases. Substitution into \eqref{eq:nu-transfer-smooth} and the cutoff
\(
 M\asymp\min\left\{n^{1/d},
       (n/\varsigma^2)^{1/D_\alpha},
       (L/\epsilon)^{1/t}\right\}
\)
give the upper bound. Entries involving division by zero are interpreted as $+\infty$; the cutoff is rounded and truncated below at one, with bounded-parameter cases handled by $W_1\le2r$.

For the lower bound, the exact-information argument of Theorem~\ref{thm:exact-information} applies even after independent noise is added. The testing argument of Proposition~\ref{prop:noisy-lower}, with spatial scale $(\varsigma^2/n)^{1/D_\alpha}$, gives the noise-dependent term uniformly over positive bounded noise levels. Finally, Lemma~\ref{lem:curved-packets} provides two laws separated by order $\epsilon^\vartheta$ whose pricing surfaces differ uniformly by at most $2\epsilon$. Choosing opposite biases to make their mean prices coincide gives the systematic term. The maximum of the three lower bounds controls their sum up to a constant. Full details are in the Appendix.
\end{proof}

\subsection{General Designs and Additional Marginal Information}
The integral in \eqref{eq:nu-cv-moment} uses the known uniform design. The same recovery rates can nevertheless be attained without knowing the sampling law, provided that every region receives sufficiently many quotes. For this purpose, we smooth a local estimate of the pricing surface and project it onto the pricing surfaces generated by the density class.

For $h\in[n^{-1/d},1]$, construct $\widehat C_{n,h}$ from all observations as in Lemma~\ref{lem:nu-pilot}. Let $J_h$ be a fixed linear smoothing operator that removes frequencies above order $h^{-1}$ in coordinate charts, as constructed in the Appendix. Define
\(
 \widehat C^{\mathrm{proj}}_{n,h}
 =\operatorname*{argmin}_{G\in\{C_f:f\in\F_\alpha(L)\}}
       \|G-J_h\widehat C_{n,h}\|_{H^{s-1}(\Z)},
\)
and let $\widehat\mu_{n,h}$ be the law with this pricing surface. The projection is well defined and measurable: its target is a compact convex subset of the Hilbert space $H^{s-1}(\Z)$, and Lemma~\ref{lem:sobolev-map} makes the inverse map to laws continuous. Write $\widehat\mu_n=\widehat\mu_{n,h}$ at the selected bandwidth.

Under $Q\ge q_*Q_0$, the construction gives
\begin{equation}\label{eq:nu-direct-upper}
 \sup_{f\in\F_\alpha(L)}\sup_{|e|\le\epsilon}
 \E_{f,e,\varsigma}W_1(\widehat\mu_{n,h},f\dd\bx)
 \le C\left[Lh^{\alpha+1}
      +\varsigma n^{-1/2}h^{-d-1/2}
      +\epsilon h^{-(s-1)}\right].
\end{equation}
The terms represent approximation, independent noise, and systematic bias. A bandwidth of order \( \min\left\{1,\max\left\{n^{-1/d}, (\varsigma^2/n)^{1/D_\alpha}, (\epsilon/L)^{1/t}\right\}\right\} \) yields the upper rate in Theorem~\ref{thm:nu-smooth}. The local construction rounds $h$ upward to a power of two, changing it by at most a factor of two. This discretization does not affect the rate. The estimator uses the observed locations and prices and may use $\alpha,L,\varsigma,\epsilon$ and the support, but does not require knowledge of $Q$.

\begin{corollary}[General i.i.d. designs]\label{cor:nu-general-design}
Under the class assumptions of Theorem~\ref{thm:nu-smooth}, every i.i.d. design $Q\ge q_*Q_0$ with $q_*>0$ satisfies
\begin{equation}\label{eq:nu-general-design}
 \mathfrak R_n^{\mathrm{ex}}(Q^{\otimes n},\alpha,L)
       \asymp n^{-(\alpha+1)/d}.
\end{equation}
Under \eqref{eq:systematic-model} with $Q$ in place of $Q_0$, the upper bound in \eqref{eq:nu-smooth-rate} also holds, uniformly over $0\le\varsigma\le\overline\varsigma$ and $0\le\epsilon\le\overline\epsilon$. If additionally $Q\le q^*Q_0$ for some finite $q^*$, the matching lower bound holds as well.
\end{corollary}

Single-asset option prices provide additional information: observing their full strike curves determines each coordinate marginal law. One might therefore expect that knowing all these marginals would improve recovery of the joint law. The next result shows that, on a suitable fixed-marginal subclass, the remaining uncertainty about dependence is already sufficient to produce all three terms of the minimax rate.

\begin{corollary}[Fixed marginal laws]\label{cor:nu-marginals}
Let $f_0$ be the base density of Lemma~\ref{lem:bumps}. Let $\mathcal P_0$ consist of the laws $f\dd\bx$ with $f\in\F_\alpha(L)$ whose one-dimensional coordinate marginals equal those of $f_0\dd\bx$.
Under the assumptions of Theorem~\ref{thm:nu-smooth}, its lower bound remains valid with $\F_\alpha(L)$ replaced by $\mathcal P_0$. Its estimator attains the matching upper bound on $\mathcal P_0$.
\end{corollary}
Corollary~\ref{cor:nu-marginals} concerns exact knowledge of the full coordinate marginals, not merely finitely many single-asset quotes. Even this additional information does not improve the worst-case rate on the specified subclass: the alternatives used for all three lower bounds have identical coordinate marginals.

\subsection{Arbitrary Probability Laws}
The preceding sharp rate assumes a Sobolev density. The universal pilot in \eqref{eq:nu-universal-pilot} also applies to atomic and singular laws. When $e=0$, combining it with the same residual correction yields the following bounds without a density assumption.

\begin{theorem}[Bounds over all laws uniformly in the noise level]\label{thm:nu-all-laws}
Fix a finite $\overline\varsigma\ge0$. Under \eqref{eq:nu-model}, there are constants $c,C>0$, independent of $n\ge1$ and $0\le\varsigma\le\overline\varsigma$, such that
\begin{equation}\label{eq:nu-all-laws}
 c\Big[n^{-\frac{1}{d}}+\Big(\frac{\varsigma^2}{n}\Big)^{\frac{1}{2d+3}}\Big]
 \le\mathfrak R_n(\Pp(B_r),\varsigma)\le C\Big[n^{-\frac{d+3}{d(2d+3)}}
                  +\Big(\frac{\varsigma^2}{n}\Big)^{\frac{1}{2d+3}}\Big].
\end{equation}
For every fixed $c_0>0$, the two sides have the same order on the regime $\varsigma^2\ge c_0n^{-3/d}$.
\end{theorem}

Thus the noise-dominated rate is sharp, but a gap remains for exact or very low-noise quotes over arbitrary laws. For example, when $d=2$ and $\varsigma=0$,
\(
 c n^{-1/2}\le\mathfrak R_n(\Pp(B_r),0)\le C n^{-5/14}.
\)
Closing this gap requires a sharper recovery argument for laws that need not have a smooth density.

\section{Price Intervals and Minimax Estimation of a Price}\label{sec:stable-pricing}
Recovering the law controls the prices of all Lipschitz claims, but estimating the price of a specified payoff can be substantially easier. We first express the uncertainty left by exact quotes as a price interval. We then introduce two claim-specific moduli: one for uniformly bounded quote discrepancies, and one for discrepancies measured in the design norm. The latter characterizes minimax estimation at a fixed positive noise level. Its dual formulation yields both an affine estimator and confidence intervals through a penalized static replication problem.

\subsection{Price Intervals and Wasserstein Diameter}\label{sec:price-intervals}
For $\K\subseteq\Pp(B_R)$ and $g\in C(B_R)$, define
\[
 I_g(\K)=\left[\min_{\nu\in\K}\int g\dd\nu,\ \max_{\nu\in\K}\int g\dd\nu\right].
\]
If $\K$ is compact and convex, every value in this interval is the price of $g$ under some law in $\K$. In particular, $I_g(\K_n(z,v))$ is the range of prices consistent with the observed exact quotes and the support constraint. Its width measures the remaining uncertainty about the claim $g$.

\begin{proposition}[Uncertainty equals Wasserstein diameter]\label{prop:price-diameter}
For compact $\K\subseteq\Pp(B_R)$,
\begin{equation}\label{eq:price-diameter}
 \sup_{\Lip(g)\le1}
 \left(\max_{\nu\in\K}\int g\dd\nu-\min_{\nu\in\K}\int g\dd\nu\right)
 =\diam_{W_1}\K.
\end{equation}
For $\K=\K_n(\mu)$ and  $g\in C(B_R)$, the endpoints of $I_g(\K_n(\mu))$ converge a.s. to the corresponding endpoints of $I_g(\K_\infty(\mu;Q))$. If $Q$ identifies $\Pp(B_R)$, both endpoints converge to $\int g\dd\mu$.
\end{proposition}

Thus a Wasserstein diameter bound is also a bound on every $1$-Lipschitz claim's price interval. More explicitly, if $Q\ge q_*Q_0$, Theorem~\ref{thm:exact-general} gives, with probability at least $1-\delta$,
\[
 |I_g(\K_n(Z,v))|
 \le C_{d,R,q_*}\Big(\frac{\log(en/\delta)}{n}\Big)^{2/[d(d+3)]}
\]
for every price vector $v$ and $g$ with $\Lip(g)\le1$. $|I|$ denotes interval length. The event depends only on the quote locations: once these locations provide adequate coverage, the bound holds for all such claims and price vectors simultaneously. No smoothness restriction is imposed on the admissible laws. If instead the laws are restricted to $\F_\alpha(L)$, Theorem~\ref{thm:exact-smooth} gives
\[
 |I_g(\K_{n,\alpha,L}(Z,v))|
 \le CL\Big(\frac{\log(en/\delta)}{n}\Big)^{(\alpha+1)/d}
\]
on an event of probability at least $1-\delta$, again simultaneously over $v$ and $1$-Lipschitz $g$. These are worst-case bounds over claims; the uncertainty for a specified payoff may be smaller.

\subsection{Uniform Quote Discrepancies and Static Replication}\label{sec:modulus}
Let $E\subseteq S^{d-1}$ be a nonempty closed set of quoted directions, put $\Z_E=E\times[-R,R]$, and let $\Mm(\Z_E)$ denote the finite signed measures on $\Z_E$. Write $\phi_z(\bx)=(\bu\cdot\bx-k)_+$ for $z=(\bu,k)$. A signed measure $m\in\Mm(\Z_E)$ specifies a static portfolio of basket calls with payoff
\[
 \cA^*m(\bx)=\int_{\Z_E}\phi_z(\bx)\,m(\mathrm dz),\qquad \bx\in B_R.
\]
This payoff is $\norm m$-Lipschitz, where $\norm m$ is the total variation norm, and \(
 \int\cA^*m\dd\eta=\int_{\Z_E}C_\eta\dd m
\)
for every $\eta \in \Mm(B_R)$. The payoffs obtainable from cash and quoted calls form
\begin{equation}\label{def:V_E}
 V_E=\{c+\cA^*m:c\in\R,\ m\in\Mm(\Z_E)\}\subset C(B_R).
\end{equation}
If $E$ has positive surface measure, let $Q_E$ be the product of normalized surface measure on $E$ and the uniform strike distribution on $[-R,R]$. Define the uniform-quote pricing modulus
\[
 \varpi_g(\epsilon)=\sup\Big\{\int g\dd(\mu-\nu):\
 \mu,\nu\in\Pp(B_R),\ \sup_{\Z_E}|C_\mu-C_\nu|\le\epsilon\Big\} \ \ \text{for } g\in C(B_R), \ \epsilon\ge0.
\]
Let $\varpi_g^{\alpha,L}(\epsilon)$ denote the same supremum restricted to laws with densities in $\F_\alpha(L)$. 
The modulus is the largest price difference compatible with quote discrepancies of at most $\epsilon$ at every quoted direction and strike. Write $\osc h=\max_{B_R}h-\min_{B_R}h$.

\begin{lemma}[Replication bound]\label{lem:replication-bound}
For every $g\in C(B_R)$, $m\in\Mm(\Z_E)$, and $\epsilon\ge0$,
\begin{equation}\label{eq:replication-bound}
 \varpi_g(\epsilon)\le\osc(g-\cA^*m)+\epsilon\norm m.
\end{equation}
In particular, $\varpi_g(\epsilon)\le\epsilon\norm m$ whenever $g=c+\cA^*m\in V_E$.
\end{lemma}

The first term bounds the price uncertainty of the unreplicated payoff; the second bounds the effect of quote discrepancies on the portfolio. Cash does not affect either price differences or oscillation. Exact replication therefore gives a linear bound in the quote tolerance, while approximate replication trades residual uncertainty against the size of the portfolio. Section~\ref{sec:directional-pricing} uses this bound to distinguish claims that are stable under restricted directional information.

\subsection{The Design Modulus and Minimax Price Estimation}\label{sec:minimax-functional}
Fix $K\in\{B_r,B_R\}$ and a nonempty convex weakly closed class $\mathcal C\subseteq\Pp(K)$. The cases used below are $\Pp(B_R)$, $\Pp(B_r)$, and $\mathcal C_{\alpha,L}=\{f\dd\bx:f\in\F_\alpha(L)\}$. Let $Q$ be a known probability measure on $\Z_E$, and observe
\[
 Z_i\overset{\mathrm{i.i.d.}}\sim Q,\qquad
 Y_i=C_\mu(Z_i)+\varsigma\xi_i,\qquad i=1,\ldots,n,
\]
where $\varsigma>0$ is fixed and the errors are independent of the design, centered, and have unit variance, as in \eqref{eq:nu-model}. There is no systematic bias in this experiment. For $g\in C(K)$, put
\[
 \mathfrak R_n^Q(g;\mathcal C)
 =\inf_{\widehat p}\sup_{\mu\in\mathcal C}
 \E_\mu\left|\widehat p-\int g\dd\mu\right|,
\]
where the infimum is over real-valued estimators based on $(Z_i,Y_i)_{i=1}^n$. The relevant discrepancy for this experiment is measured in $L^2(Q)$ rather than uniformly on $\Z_E$. Define
\[
 \omega_g(\epsilon;\mathcal C)
 =\sup\Big\{\int g\dd(\mu-\nu):\
 \mu,\nu\in\mathcal C,\ \norm{C_\mu-C_\nu}_{L^2(Q)}\le\epsilon\Big\},
 \qquad \epsilon\ge0,
\]
and define its price-range width over $\mathcal C$ by
\(
 \osc_{\mathcal C}(h)
 =\max_{\mu\in\mathcal C}\int h\dd\mu
  -\min_{\mu\in\mathcal C}\int h\dd\mu
\)
for $h\in C(K)$. For $\mathcal C=\Pp(K)$, this width is $\max_Kh-\min_Kh$. For $w\in L^2(Q)$, the portfolio
\[
 \cA^*(wQ)(\bx)=\int_{\Z_E}\phi_z(\bx)w(z)\,Q(\mathrm dz)
\]
is bounded in absolute value by $2R\norm w_{L^1(Q)}$ and is $\norm w_{L^1(Q)}$-Lipschitz on $K$.

\begin{theorem}[Minimax estimation of a price]\label{thm:minimax-functional}
Let $g\in C(K)$.
\begin{enumerate}[label=(\roman*)]
\item For every $\epsilon\ge0$,
\begin{equation}\label{eq:l2-modulus-duality}
 \omega_g(\epsilon;\mathcal C)
 =\inf_{w\in L^2(Q)}\left\{
 \osc_{\mathcal C}(g-\cA^*(wQ))
 +\epsilon\norm w_{L^2(Q)}\right\}.
\end{equation}
\item Suppose additionally, for the lower bound, that $I_\xi < \infty$ as in Section~\ref{sec:model}. For every $n\ge1$,
\begin{equation}\label{eq:minimax-functional}
 \frac{1}{8\sqrt{1\vee I_\xi}}\,
 \omega_g\left(\frac{\varsigma}{\sqrt n};\mathcal C\right)
 \le\mathfrak R_n^Q(g;\mathcal C)
 \le\sqrt{1+\frac{4R^2}{\varsigma^2}}\,
 \omega_g\left(\frac{\varsigma}{\sqrt n};\mathcal C\right).
\end{equation}
The upper bound is attained up to an arbitrarily small additive error by
\(
 \widehat p_w=c_w+\frac1n\sum_{i=1}^nw(Z_i)Y_i,
\)
where $w$ is a deterministic near-minimizer in \eqref{eq:l2-modulus-duality} at $\epsilon=\varsigma/\sqrt n$, and
\(
 c_w=\frac12\left\{
 \max_{\mu\in\mathcal C}\int(g-\cA^*(wQ))\dd\mu
 +\min_{\mu\in\mathcal C}\int(g-\cA^*(wQ))\dd\mu\right\}.
\)

\end{enumerate}
\end{theorem}

The dual objective balances residual price uncertainty against the $L^2(Q)$ size of the portfolio weights. The midpoint $c_w$ centers the residual price range, giving worst-case absolute bias at most $\tfrac12\osc_{\mathcal C}(g-\cA^*(wQ))$. Moreover, for $\mu\in\mathcal C$, 
\(
 \operatorname{Var}_\mu(\widehat p_w)
 \le\frac{4R^2+\varsigma^2}{n}\norm w_{L^2(Q)}^2. 
 \)
This variance includes both quote noise and variation due to the random locations. Consequently, the comparison in \eqref{eq:minimax-functional} is a fixed-positive-noise result: its upper comparison constant diverges as $\varsigma\downarrow0$, so setting $\varsigma=0$ does not characterize the exact-quote risk.

The two moduli serve different purposes. Uniform quote control implies $L^2(Q)$ control, so
\(
 \varpi_g(\epsilon)\le\omega_g(\epsilon;\Pp(B_R))
\) and 
\( \varpi_g^{\alpha,L}(\epsilon)\le\omega_g(\epsilon;\mathcal C_{\alpha,L}).
\) 
Thus a lower bound for the uniform-quote modulus also supplies a lower bound for the statistical modulus on the same class. The converse comparison requires additional assumptions on the design. Likewise, a representing measure $m$ in Lemma~\ref{lem:replication-bound} need not have an $L^2(Q)$ density. Membership in $V_E$ alone therefore does not give an $L^2(Q)$ replication bound. 

Theorem~\ref{thm:minimax-functional} is a random-design version of the connection between optimal recovery and affine estimation in \cite{donoho1994}. Its proof also applies to a measurable family of continuous test functions $\phi_z$ on a compact metric space with $\sup_{z,x}|\phi_z(x)|\le B$. In that setting, $B^2$ replaces $4R^2$ in the risk and variance bounds, and the same modulus, dual formulation, and affine construction apply to any nonempty convex weakly closed class of probability laws.

\subsection{Confidence Intervals and Optimal Length}\label{sec:price-confidence}
The same bias and variance bounds give confidence intervals without a Gaussian or tail assumption on the errors. The weights below are chosen from the known class and design, independently of the observed sample.

\begin{corollary}[Confidence intervals for a price]\label{cor:price-ci}
Under the setting of Theorem~\ref{thm:minimax-functional}, let $g\in C(K)$ and $0<\delta<\tfrac14$. For a deterministic $w\in L^2(Q)$, let $\widehat p_w$ and $c_w$ be as above, and put
\(
 \widehat I_w=[\widehat p_w-\ell_w,\ \widehat p_w+\ell_w]
 \)
  and 
 \(
 \ell_w=\frac12\osc_{\mathcal C}(g-\cA^*(wQ))
 +\sqrt{\frac{4R^2+\varsigma^2}{n\delta}}\norm w_{L^2(Q)}.
\)
\begin{enumerate}[label=(\roman*)]
\item For every such $w$,
\(
 \inf_{\mu\in\mathcal C}\Pr_\mu\left(\int g\dd\mu\in\widehat I_w\right)\ge1-\delta.
\) 
For every $\eta>0$, one can choose $w$ so that
\(
 \ell_w\le\delta^{-1/2}\sqrt{1+\frac{4R^2}{\varsigma^2}}
 \left\{\omega_g\left(\frac{\varsigma}{\sqrt n};\mathcal C\right)+\eta\right\}.
\) 
The full interval length is $|\widehat I_w|=2\ell_w$.
\item If the error density has finite Fisher information $I_\xi$, every data-measurable interval $\widehat I$ satisfying
\(
 \inf_{\mu\in\mathcal C}\Pr_\mu\left(\int g\dd\mu\in\widehat I\right)\ge1-\delta
\)
\,obeys\,
\(
 \sup_{\mu\in\mathcal C}\E_\mu|\widehat I|
 \ge\frac{\tfrac12-2\delta}{\sqrt{1\vee I_\xi}}\,
 \omega_g\left(\frac{\varsigma}{\sqrt n};\mathcal C\right).
\)
\end{enumerate}
\end{corollary}

For fixed $\delta\in(0,1/4)$ and $\varsigma>0$, these bounds show that the smallest possible worst-case expected length among all intervals with uniform coverage at least $1-\delta$ is comparable to $\omega_g(\varsigma/\sqrt n;\mathcal C)$. The intervals constructed above attain this order. This is an optimality statement for worst-case expected length under uniform coverage. Unlike the exact-quote price intervals of Section~\ref{sec:price-intervals}, these are statistical confidence intervals.

\section{Directional Pricing and Exchange Options}\label{sec:directional-pricing}
Long-only basket quotes identify the joint law, yet the accuracy of an inferred price depends strongly on the payoff. Polynomial claims admit parametric estimation at fixed positive noise, while an exchange option in an unquoted spread direction has only logarithmic recovery. We first explain why missing directions make inversion unstable, then identify claims that remain well priced. We next establish sharp logarithmic moduli for unquoted kinks, compare uncertainty across classes of claims, and extend the exchange-price result to exact and noisy quotes.

\subsection{Missing Projection Directions}
Identifying a law does not imply a stable inverse. The relevant geometry concerns projection lines rather than oriented directions: $t\mapsto\widehat\mu(t\bu)$ covers both signs of $\bu$. Thus a hemisphere meets every projection line, whereas long-only directions miss lines with mixed signs, including those near the spread direction $(\be_1-\be_2)/\sqrt2$.

\begin{theorem}[Missing-cone instability]\label{thm:missing-cone}
Suppose there are $\be\in S^{d-1}$, $\kappa>0$ such that
\begin{equation}\label{eq:missing-cone}
 \abs{\bu\cdot\be}\le1-\kappa\qquad\text{for all }\bu\in E.
\end{equation}
For every $\alpha\ge0$ and every sufficiently large $L$, there are densities $f_N^+,f_N^-\in\F_\alpha(L)$ satisfying
\begin{equation}\label{eq:limited-w1}
 W_1(f_N^+\dd\bx,f_N^-\dd\bx)\ge cN^{-\alpha-1},
\end{equation}
whereas, for each $A>0$, there is $C_A<\infty$ such that
\begin{equation}\label{eq:limited-price}
 \sup_{\bu\in E,\,k\in[-R,R]}
       \abs{C_{f_N^+}(\bu,k)-C_{f_N^-}(\bu,k)}
                              \le C_A N^{-A}.
\end{equation}
Consequently, for every $\gamma>0$, no finite constant $C$ can make
\[
 W_1(f\dd\bx,g\dd\bx)
 \le C\Big(\sup_{\bu\in E,\,k\in[-R,R]}
                \abs{C_f(\bu,k)-C_g(\bu,k)}\Big)^\gamma
\]
hold uniformly over $f,g\in\F_\alpha(L)$. The same conclusion holds with an $L^p$ pricing norm for any probability design supported on $E\times[-R,R]$ and $1\le p\le\infty$.
\end{theorem}

\begin{corollary}[Long-only basket directions]\label{cor:long-only}
Let $E=\{\bu\in S^{d-1}:u_j\ge0\text{ for all }j\}$. Call prices on $E\times[-R,R]$ identify the  joint law, but admit no uniform H\"older inverse bound as in Theorem~\ref{thm:missing-cone} on a sufficiently rich bounded Sobolev class.
\end{corollary}

\subsection{Stable Claims and Polynomial Prices}\label{sec:long-only}
Instability of the joint law need not make every claim difficult to price. A payoff represented by cash and a portfolio of quoted calls inherits the stability bound of Lemma~\ref{lem:replication-bound}. The next proposition gives approximation bounds.

\begin{proposition}[Stable claims]\label{prop:positive}
Let $\bv\in E$, $G:[-R,R]\to\R$, and $g(\bx)=G(\bv\cdot\bx)$.
\begin{enumerate}[label=(\roman*)]
\item If $G'$ has bounded variation, then $\varpi_g(\epsilon)\le(\abs{G'(-R+)}+\abs{G''}((-R,R)))\,\epsilon$.
\item If $G$ is $1$-Lipschitz, then $\varpi_g(\epsilon)\le4\sqrt{2R\epsilon}+\epsilon$.
\item If $E$ has nonempty relative interior in $S^{d-1}$, every polynomial $g$ belongs to $V_E$ defined in \eqref{def:V_E} and is Lipschitz stable, i.e., $\varpi_g(\epsilon)\le C_g\epsilon$.
\end{enumerate}
\end{proposition}

\begin{example}\label{ex:squared-spread}
Let $E=\{\bu\in S^{d-1}:u_j\ge0\ \forall j\}$ and $\bv=(\be_1+\be_2)/\sqrt2\in E$. The squared spread $(x_1-x_2)^2=2x_1^2+2x_2^2-2(\bv\cdot\bx)^2$ is a static portfolio of long-only calls by Proposition~\ref{prop:positive}(i), hence it is Lipschitz stable. In particular, index and single-asset second moments determine the covariance, and hence the second moment of a spread.
\end{example}

\begin{example}[Polynomial claims from long-only quotes]\label{ex:polynomial-rate}
Let $E$ have nonempty relative interior, let $Q\ge q_*Q_E$ be known, and fix $\varsigma>0$. For every nonconstant polynomial $g$,
\[
 cn^{-\frac{1}{2}}\le\mathfrak R_n^Q\left(g;\Pp(B_r)\right)\le Cn^{-\frac{1}{2}}.
\]
The derivation is in the Appendix. 
Since $\omega_g(\epsilon;\Pp(B_r))\asymp\epsilon$ for small $\epsilon$, Corollary~\ref{cor:price-ci} gives confidence intervals of length of order $n^{-1/2}$ for these polynomial claims.
\end{example}

The coverage of the restricted quote domain determines how well these stable claims can be priced. For a finite design $z\subset\Z_E$, let $h_E(z)=\sup_{w\in\Z_E}\min_i\rho(w,z_i)$ be the fill distance. For the restricted-design measure $Q_E$ defined in Section~\ref{sec:modulus}, we use the regularity condition
\begin{equation}\label{eq:ahlfors}
 \mathcal H^{d-1}(E\cap B(\bu,t))\ge c_Et^{d-1},\qquad\bu\in E,\ 0<t\le1,
\end{equation}
where $\mathcal H^{d-1}$ is surface measure. The long-only set satisfies \eqref{eq:ahlfors}. 

\begin{corollary}[Finitely many quotes]\label{cor:finite} 
For every continuous $g$ and feasible price vector $v$, the price interval of $g$ over $\K_n(z,v)$ has width at most $\varpi_g(2Rh_E(z))$. If \eqref{eq:ahlfors} holds and the quote locations are i.i.d.\ $Q$ with $Q\ge q_*Q_E$, then with probability at least $1-\delta$, 
\[ 
h_E(Z_1,\ldots,Z_n) \le C_{d,q_*, c_E} \Big(\frac{\log(en/\delta)}{n}\Big)^{1/d}. 
\] 
On this event, simultaneously for every $g=c+\cA^*m\in V_E$ and $v$, the price interval has width at most $2RC\norm m\allowbreak(\log(en/\delta)/n)^{1/d}$. The exponent $1/d$ exceeds the exponent $2/(d(d+3))$ in the bound of Theorem~\ref{thm:exact-general} for all Lipschitz claims under full directional coverage. 
\end{corollary}

The missing-direction phenomenon also affects statistical recovery of the full law. In the noisy model of Section~\ref{sec:minimax-functional}, with $Q$ any probability measure on $\Z_E$ and $\varsigma>0$ fixed, define
\[
 \mathfrak R_n^Q(\alpha,L)
   =\inf_{\widehat\mu}\sup_{f\in\F_\alpha(L)}
                \E_{f,Q}W_1(\widehat\mu,f\dd\bx), \quad \text{where} \quad \widehat\mu = \widehat\mu((Z_i,Y_i)_{i=1}^n).
\]

\begin{corollary}[Absence of a polynomial noisy rate]\label{cor:missing-noisy}
Suppose $E$ satisfies \eqref{eq:missing-cone}. For every $\alpha\ge0$, $\gamma>0$ and every sufficiently large $L$, 
 $\lim_{n\to\infty}n^\gamma\mathfrak R_n^Q(\alpha,L)=+\infty$. 
In particular, this holds for noisy long-only basket quotes.
\end{corollary}

\subsection{Quoted and Unquoted Kinks}\label{sec:log-modulus}
We now turn from recovery of the law to the price of one kinked payoff. If its direction is quoted, directly or with the opposite orientation, the claim is statically replicable and its price is Lipschitz stable. Otherwise, the price can be identified but is much less stable. The next theorem gives the logarithmic order of this loss for a strike inside the support. Throughout, $E\subseteq S^{d-1}$ is closed with nonempty relative interior, $\be\in S^{d-1}$, and $g_k(\bx)=(\be\cdot\bx-k)_+$. For a finite signed measure $\eta$ on $B_R$, $\eta_{\be}$ denotes its image under $\bx\mapsto\be\cdot\bx$. Its Fourier transform satisfies $\widehat{\eta_{\be}}(t)=\widehat\eta(t\be)$.

\begin{theorem}[Quoted and unquoted kinks]\label{thm:log-modulus} 
If $\be\in E$ or $-\be\in E$, then $\varpi_{g_k}(\epsilon)\le2\epsilon$ for all $\abs k\le R$ and $\epsilon\ge0$. If $\pm\be\notin E$ and $\abs k<r$, then for every $\alpha\ge0$ and every sufficiently large $L$ there are $\epsilon_0\in(0,1)$ and $0<c<C<\infty$ such that, for $0<\epsilon\le\epsilon_0$, 
\begin{equation}\label{eq:log-modulus} 
c\Big(\log\frac1\epsilon\Big)^{-(\alpha+\frac32)} \le\varpi_{g_k}^{\alpha,L}(\epsilon) \le C\Big(\log\frac1\epsilon\Big)^{-(\alpha+\frac32)}. 
\end{equation} 
The lower bound holds for every closed $E$ satisfying \eqref{eq:missing-cone} with direction $\be$. 
\end{theorem}

\begin{proof}[Proof sketch] 
The quoted-direction bound follows from
Lemma~\ref{lem:replication-bound} and the put--call parity identity $(y-k)_+=(-y+k)_++y-k$.
For the upper bound in \eqref{eq:log-modulus}, express
$\int g_k\dd\eta$ through the cosine transform of
$\eta_{\be}$ and split at frequency $N$.
Lemma~\ref{lem:long-only-powers} bounds the low-frequency
part by $C\epsilon e^{cN}$, while
Lemma~\ref{lem:line-restriction} bounds the tail by
$CLN^{-(\alpha+3/2)}$.
Taking $N$ to be a sufficiently small multiple of
$\log(1/\epsilon)$ proves the upper bound.

For the lower bound, Lemmas~\ref{lem:window}
and~\ref{lem:profile} give a localized oscillatory
perturbation with price gap of order $N^{-(\alpha+3/2)}$
and quote discrepancy at most $Ce^{-cN}$ on $\Z_E$.
A nonnegative envelope preserves positivity and the uniform
$H^\alpha$ bound. Taking $N$ to be a sufficiently large
multiple of $\log(1/\epsilon)$ proves the lower bound.
Full details are in the Appendix.
\end{proof}

The upper bound guarantees price error at most $a$ when $\epsilon\le\exp(-Ca^{-1/(\alpha+3/2)})$; the lower bound requires $\exp(-ca^{-1/(\alpha+3/2)})$. Since adding a payoff in $V_E$ changes the modulus by at most a linear term in $\epsilon$, the bounds in Theorem~\ref{thm:log-modulus} extend to $g_E+c\,g_k$ with $g_E\in V_E$ and $c\ne0$.

\begin{example}\label{ex:long-only-exotics}
Let $E$ be the long-only directions and $\be=(\be_1-\be_2)/\sqrt2$, so that $\pm\be\notin E$. Since $x_j=\phi_{(\be_j,-R)}-R\in V_E$, Theorem~\ref{thm:log-modulus} gives moduli of order $(\log(1/\epsilon))^{-(\alpha+3/2)}$ over $\F_\alpha(L)$ for the exchange option $(x_1-x_2)_+=\sqrt2(\be\cdot\bx)_+$, the spread call $(x_1-x_2-K)_+$ with $\abs K<\sqrt2r$, the best-of payoff $\max(x_1,x_2)=x_2+(x_1-x_2)_+$, and the worst-of payoff $\min(x_1,x_2)=x_1-(x_1-x_2)_+$. All four prices are identified by the full long-only pricing surface, but become Lipschitz stable when the spread direction $\be$ is quoted. Together with Example~\ref{ex:squared-spread}, this shows that polynomial and kinked payoffs of the same unquoted spread can have very different pricing stability.
\end{example}

\begin{corollary}[Minimax pricing of a kink from noisy quotes]\label{cor:log-minimax}
Let $E$ be closed and satisfy the missing-cone condition with direction $\be$, let $Q$ be a probability measure on $\Z_E$, $\varsigma>0$ fixed, $\abs k<r$, $\alpha\ge0$, $L$ sufficiently large, $g_E\in V_E$, $c\ne0$, and $g=g_E+c\,g_k$.
\begin{enumerate}[label=(\roman*)]
\item There are $c'>0$ and $n_0$ such that, for $n\ge n_0$,
\(
 \mathfrak R_n^Q(g;\mathcal C_{\alpha,L})
 \ge c'(\log n)^{-(\alpha+\frac{3}{2})}.
\)
\item If $E$ also has nonempty relative interior and satisfies \eqref{eq:ahlfors}, and $Q\ge q_*Q_E$, there are $C<\infty$ and $n_0$ such that, for $n\ge n_0$,
\(
 \mathfrak R_n^Q(g;\mathcal C_{\alpha,L})
 \le C(\log n)^{-(\alpha+\frac{3}{2})}.
\)
\end{enumerate}
\end{corollary}

Corollary~\ref{cor:price-ci} then gives confidence intervals for the exchange price of length of order $(\log n)^{-(\alpha+3/2)}$ when $q_*Q_E\le Q$.

\subsection{Claim-Indexed Uncertainty}\label{sec:claim-indexed}
By \eqref{eq:price-diameter}, the $W_1$-diameter of $\K\subseteq\Pp(B_R)$ is the largest price uncertainty among all $1$-Lipschitz claims. The preceding results show why it is useful to measure uncertainty over a specified class of claims rather than over all Lipschitz payoffs. 

For $\mathcal G\subseteq C(B_R)$ closed under $g\mapsto-g$, let $d_{\mathcal G}(\mu,\nu)=\sup_{g\in\mathcal G}\int g\dd(\mu-\nu)$. Then every compact $\varnothing \neq \K\subseteq\Pp(B_R)$ satisfies $\diam_{d_{\mathcal G}}\K =\sup_{g\in\mathcal G} \big(\max_{\nu\in\K}\int g\dd\nu -\min_{\nu\in\K}\int g\dd\nu\big)$. 
The single claim $\mathcal G=\{\pm g\}$ gives the width of the price interval of $g$, whose worst case under quote errors is the modulus $\varpi_g$ of Section~\ref{sec:modulus}. The class $\mathcal G=\{\text{1-Lipschitz functions}\}$ gives $W_1$. Between these extremes, for a target set $D\subseteq S^{d-1}$ of basket directions let
\[
 \mathcal G_D=\left\{\bx\mapsto G(\bv\cdot\bx):\ \bv\in D,\ G:[-R,R]\to\R\ \text{1-Lipschitz}\right\},
 \qquad d_D=d_{\mathcal G_D}.
\]
Then $d_D(\mu,\nu)=\sup_{\bv\in D}W_1(\mu_{\bv},\nu_{\bv})$ is the max-sliced Wasserstein distance over $D$ \cite{deshpande2019}, where $\mu_{\bv}$ is the law of $\bv\cdot X$ under $\mu$. Let $\varpi_D(\epsilon)$ be the supremum of $d_D(\mu,\nu)$ over $\mu,\nu\in\Pp(B_R)$ satisfying $\sup_{\Z_E}\abs{C_\mu-C_\nu}\le\epsilon$, and $\varpi_D^{\alpha,L}(\epsilon)$ its restriction to densities in $\F_\alpha(L)$.

\begin{theorem}[Directional moduli]\label{thm:directional}
\begin{enumerate}[label=(\alph*)]

\item For every nonempty $D$ and closed $E$, $\varpi_D(\epsilon)\ge(\sqrt2\pi)^{-1}\sqrt{R\epsilon}$ for $0<\epsilon\le R/(2\pi^2)$. If $D\subseteq E$, then $\varpi_D(\epsilon)\le4\sqrt{2R\epsilon}+\epsilon$ for all $\epsilon>0$.

\item Suppose $\be\in D$ and $E$ satisfies \eqref{eq:missing-cone} with direction $\be$. Then for every $\alpha\ge0$ and sufficiently large $L$ there are $c,\epsilon_0>0$ such that $\varpi_D^{\alpha,L}(\epsilon) \ge c(\log(1/\epsilon))^{-(\alpha+1)}$ for $0<\epsilon\le\epsilon_0$.

\item If $E$ has nonempty relative interior, $\pm\be\notin E$, and $D\subseteq E\cup\{\be\}$, then $\varpi_D^{\alpha,L}(\epsilon)\le C(\log(1/\epsilon))^{-(\alpha+1)}$ for $0<\epsilon\le\epsilon_0$. If moreover $\be\in D$ and $L$ is sufficiently large, then $\varpi_D^{\alpha,L}(\epsilon)\asymp(\log(1/\epsilon))^{-(\alpha+1)}$.
\end{enumerate}
\end{theorem}

Even full directional coverage cannot improve the exponent $1/2$ in (a) for arbitrary one-dimensional Lipschitz claims. The exponent $\alpha+1$ concerns worst-case uncertainty over all Lipschitz payoffs along an unquoted direction, whereas $\alpha+\tfrac32$ in \eqref{eq:log-modulus} concerns a fixed kink. The Lipschitz class includes oscillatory payoffs that can vary with the quote tolerance.

\begin{corollary}\label{cor:indexed-finite}
Assume \eqref{eq:ahlfors} and i.i.d. quote locations with $Q\ge q_*Q_E$. With probability at least $1-\delta$, simultaneously for all feasible price vectors, the unrestricted admissible set satisfies $\diam_{d_D}\K_n\le C(\log(en/\delta)/n)^{1/(2d)}$ whenever $D\subseteq E$. If $E$ has nonempty relative interior, $\pm\be\notin E$, and $D\subseteq E\cup\{\be\}$, the Sobolev-restricted admissible set satisfies $\diam_{d_D}\K_{n,\alpha,L}\le C(\log n)^{-(\alpha+1)}$ for $n\ge n_0(\delta)$.
\end{corollary}

\begin{table}[ht]
\centering
\caption{Pricing moduli for different classes of claims. }
\small
\renewcommand{\arraystretch}{1.15}
\setlength{\tabcolsep}{4pt}
\begin{tabularx}{\textwidth}{@{}>{\raggedright\arraybackslash}p{0.36\textwidth}>{\raggedright\arraybackslash}p{0.24\textwidth}>{\raggedright\arraybackslash}X>{\raggedright\arraybackslash}p{0.14\textwidth}@{}}
\toprule
Claims & Information & Modulus & Result\\
\midrule
Static portfolios $c+\cA^*m$ & Any direction set $E$ & $\le\norm m \epsilon$ & Lem.~\ref{lem:replication-bound}\\
Lipschitz payoffs on a quoted basket & $D\subseteq E$, all laws & Order $\epsilon^{1/2}$ & Thm.~\ref{thm:directional}(a)\\
A kink on an unquoted spread & Long-only; $\F_\alpha(L)$ & $(\log(1/\epsilon))^{-(\alpha+3/2)}$ & Thm.~\ref{thm:log-modulus}\\
All Lipschitz payoffs on that spread & Long-only; $\F_\alpha(L)$ & $(\log(1/\epsilon))^{-(\alpha+1)}$ & Thm.~\ref{thm:directional}(b,c)\\
All Lipschitz payoffs ($W_1$) & Full directions; all laws & $\le C\epsilon^{2/(d+3)}$ & Thm.~\ref{thm:modulus}\\
\bottomrule
\end{tabularx}
\label{tab:summary}
\end{table}

\subsection{Exact Quotes and Uniformity in the Noise Level}\label{sec:exchange-uniform}
The preceding minimax bound assumes fixed positive noise. Does the exchange price become easier to estimate when the quotes are exact? Under two-sided bounds on the sampling density, the answer is no: the same logarithmic rate holds uniformly down to zero noise. The key difficulty is not measurement noise but the information missing from finitely many long-only quotes.

Throughout this subsection $E$ has nonempty relative interior and satisfies \eqref{eq:ahlfors}, $\pm\be\notin E$, and there is $\boldsymbol w\in S^{d-1}$ with $\boldsymbol w\cdot\be=0$ and $\min_{\bu\in E}\bu\cdot\boldsymbol w>0$. For the long-only directions and $\be=(\be_1-\be_2)/\sqrt2$ one can take $\boldsymbol w=d^{-1/2}(1,\ldots,1)$, with $\min_{\bu\in E}\bu\cdot\boldsymbol w\ge d^{-1/2}$. We write $\mathfrak R_n^Q(g;\mathcal C,\varsigma)$ for the risk defined in Section~\ref{sec:minimax-functional}, now displaying its dependence on the noise level and allowing $\varsigma=0$.

\begin{theorem}[Pricing an unquoted kink uniformly in the noise level]\label{thm:exchange-uniform}
Let $q_*Q_E\le Q\le q^*Q_E$, $\abs k<r$, $\alpha\ge0$, $L$ sufficiently large, and $\overline\varsigma<\infty$. There are $0<c<C<\infty$ and $n_0$, independent of $n$ and $\varsigma$, such that for all $n\ge n_0$ and $0\le\varsigma\le\overline\varsigma$,
\[
 c(\log n)^{-(\alpha+\frac{3}{2})}
 \le\mathfrak R_n^Q\bigl(g_k;\mathcal C_{\alpha,L},\varsigma\bigr)
 \le C(\log n)^{-(\alpha+\frac{3}{2})}.
\]
\end{theorem}
Overall, long-only quotes identify the joint law, but identification alone gives no useful polynomial rate for an unquoted spread payoff. Theorem~\ref{thm:exchange-uniform} shows that even exact quotes yield only $(\log n)^{-(\alpha+3/2)}$ minimax risk for the exchange price. Yet the same quotes estimate polynomial prices at the parametric rate. Adding the spread direction $(\be_1-\be_2)/\sqrt2$ makes all four payoffs in Example~\ref{ex:long-only-exotics} Lipschitz stable. Thus pricing stability depends on a claim's direction relative to the quoted directions, not merely on quote precision.

\section{Numerical Experiments}\label{sec:numerical}
We illustrate three consequences of the theory in dimension $d=2$: the benefit of residual correction, recovery under nonuniform designs, and the effect of adding a spread direction when pricing a claim. Call prices are evaluated analytically with Gaussian quote noise $N(0,\varsigma^2)$. The reference design is $Q_0=\operatorname{Unif}(S^1)\otimes\operatorname{Unif}([-1,1])$. We compute $W_1$ error by exact Euclidean transport between cell masses on an $81\times81$ grid. 
The recovery and calibration densities have the form
\[
 f_{\beta,\boldsymbol c,A}(\boldsymbol x)
 =\frac{\beta+1}{\pi|\det A|}
   \bigl(1-\|A^{-1}(\boldsymbol x-\boldsymbol c)\|^2\bigr)_+^\beta,
 \qquad A=\mathsf U_\varphi\operatorname{diag}(a,b),
\]
where $\mathsf U_\varphi$ is the rotation through angle $\varphi$. All three laws use $\beta=0.6$ and have $H^1$ densities supported in $B_r$ with $r=0.85$.

\begin{itemize}
\setlength{\itemsep}{1pt}\setlength{\parskip}{0pt}\raggedright
\item \textit{Ellipse.} A single density with $\boldsymbol c=(0.08,-0.04)$ and $(a,b,\varphi)=(0.64,0.38,0.61)$.
\item \textit{Mixture.} Two components with weights $(0.58,0.42)$; centers $(-0.25,-0.13)$ and $(0.27,0.18)$; $(a,b,\varphi)$ equal to $(0.30,0.21,-0.35)$ and $(0.31,0.40,0.80)$, respectively.

\item \textit{Calibration.} A mixture with weights $(0.65,0.35)$; centers $(-0.12,0.13)$ and $(0.24,-0.14)$; $(a,b,\varphi)$ equal to $(0.46,0.36,-0.45)$ and $(0.26,0.32,0.20)$, respectively.
\end{itemize}

\subsection{Residual correction}
\begin{figure}[!htbp]
\centering
\includegraphics[width=\linewidth]{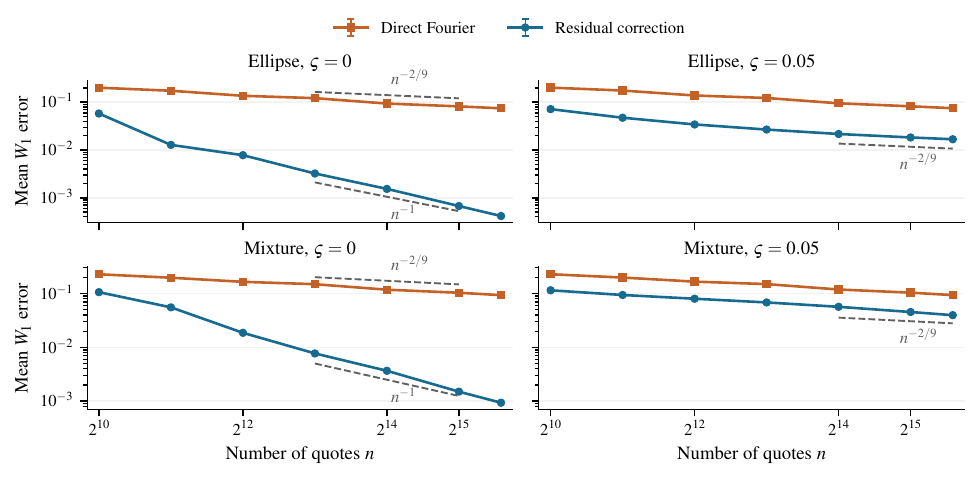}
\caption{Mean $W_1$ error of the probability estimates. Dashed lines give theoretical rates.}\label{fig:numerical-recovery}
\end{figure}

For $\alpha=1$, the direct Fourier bound in Proposition~\ref{prop:noisy-upper} gives $n^{-2/9}$, whereas the minimax benchmark in Theorem~\ref{thm:nu-smooth} is $n^{-1}+(\varsigma^2/n)^{2/9}$. Residual correction is designed to recover the faster exact-quote behavior. Figure~\ref{fig:numerical-recovery} compares the two estimators across the exact-quote and noisy regimes. 
We use torus side $5$, noise levels $\varsigma\in\{0,0.05\}$, and cutoffs $M_{\mathrm{dir}}=1.4\,n^{1/9}$ and $M_{\mathrm{cv}}=\min\{0.1\sqrt n,\, 0.8(n/\varsigma^2)^{1/9}\}$, with the second term interpreted as infinity when $\varsigma=0$.
The pilot is a local quartic fit on the $60$ nearest quotes. We average the two estimates obtained by exchanging the sample halves. After a smooth support cutoff and mass correction, we project onto probability histograms supported in $B_{0.85}$ in a discrete $H^{-1}$ norm, using $\|\cdot\|_{\mathrm{Lip}^*}\le C\|\cdot\|_{H^{-1}}$. The spatial grid has $J=\max(64,8\lceil M\rceil)$ cells across $[-1,1]$. Results are averaged over $100$ replications.

With exact quotes, residual correction reduces mean $W_1$ error by about two orders of magnitude at the largest sample size. The corrected curves decrease faster than the $n^{-1}$ reference line. Notably, the simulated densities belong to $H^\alpha$ for every $\alpha<1.1$. At $\varsigma=0.05$, Theorem~\ref{thm:nu-smooth} predicts the slower rate $n^{-2/9}$, which closely matches direct Fourier recovery. For $n\ge2048$, the corrected curves have empirical exponents $0.32$ and $0.27$ for the ellipse and mixture. At the largest sample size, residual correction reduces the ellipse and mixture errors by factors of about 4.5 and 2.4, respectively.

\subsection{Recovery under general designs}
Corollary~\ref{cor:nu-general-design} extends the rate $n^{-1}+(\varsigma^2/n)^{2/9}$ to designs with densities bounded above and away from zero relative to $Q_0$. We compare recovery under $Q_0$ and $Q_1=\tfrac12Q_0+\tfrac12Q_{\mathrm{c}}$, constructing $\widehat\mu_h$ by smoothing and projection onto a finite mixture dictionary. The estimator uses the observed quote locations, not the design density. Its bandwidth is selected by comparing candidate reconstructions, with a threshold motivated by the stochastic term in \eqref{eq:nu-direct-upper}. Under $Q_{\mathrm{c}}$, the angle has a von Mises distribution with mean $\pi/4$ and concentration $2$, and the strike is $2U-1$ for an independent $U\sim\operatorname{Beta}(2,2)$.

For bandwidth $h$, we use a local quartic pilot with $\max(60,\lceil4nh^2\rceil)$ neighbors. Its values are truncated to $[0,2R]$ as in Lemma~\ref{lem:nu-pilot}. For the three sample sizes in Table~\ref{tab:numerical-design}, the pricing grids have $48\times81$, $64\times113$ and $80\times129$ equally spaced direction--strike nodes, respectively. We use an angular Fourier and strike cosine basis for smoothing. Projection minimizes the discrete $H^{3/2}$ distance to the smoothed surface over convex mixtures of $f_{0.6,\boldsymbol c,aI}$. 

We select the largest $h$ in $\mathcal H_n=\{2^{j/4}:j\in\mathbb Z\}\cap[n^{-1/2},1/4]$ satisfying $W_1(\widehat\mu_h,\widehat\mu_{h'})\le0.35\,\widehat\varsigma n^{-1/2}(h')^{-5/2}$ for every $h'<h$ in $\mathcal H_n$. The common multiplier was chosen from a fixed grid by mean $W_1$ on $16$ independent calibration datasets, using the calibration law specified above and both designs at $n=4096$ and $8192$. We estimate $\varsigma$ by the root mean square of residuals from leave-one-out quartic fits on the $60$ nearest neighbors, after dividing each residual by $(1+\sum_j a_{ij}^2)^{1/2}$, where $a_{ij}$ are the local regression weights. Table~\ref{tab:numerical-design} reports results from $50$ replications. At each noise level $\varsigma\in\{0.005,0.05\}$, the oracle error is the mean of the smallest candidate $W_1$ error on each replication.

Increasing $n$ from $2,048$ to $16,384$ reduces the exact-quote errors by factors of $6.5$ to $9.1$, close to the factor $8$ predicted by $n^{-1}$. At $\varsigma=0.005$ and $0.05$, the factors are $1.7$ to $1.8$ and $1.8$ to $1.9$, respectively, compared with $8^{2/9}\approx1.6$ from the noise term. The mean errors under $Q_1$ differ from those under $Q_0$ by at most $26\%$. At the largest sample size, the selected errors are at most 12\% above the oracle errors.

\begin{table}[!htbp]
\centering\small
\caption{Recovery under the two designs. Entries are $10^3$ times mean $W_1$ error. Parentheses give Monte Carlo standard errors.}\label{tab:numerical-design}
\setlength{\tabcolsep}{5pt}
\begin{tabular}{lcrrrrrr}
\toprule
Law & Design & $n$ & \multicolumn{3}{c}{Selected, $\varsigma$} & \multicolumn{2}{c}{Oracle, $\varsigma$}\\
\cmidrule(lr){4-6}\cmidrule(lr){7-8}
 & & & $0$ & $0.005$ & $0.05$ & $0.005$ & $0.05$\\
\midrule
Ellipse & $Q_0$ & 2,048 & 2.24 (0.02) & 11.79 (0.19) & 34.27 (1.20) & 10.95 & 31.82\\
 &  & 8,192 & 0.56 ($<0.01$) & 8.22 (0.14) & 23.04 (0.77) & 7.75 & 20.82\\
 &  & 16,384 & 0.35 ($<0.01$) & 6.91 (0.11) & 19.48 (0.52) & 6.52 & 17.50\\
\addlinespace[3pt]
 & $Q_1$ & 2,048 & 2.60 (0.03) & 11.58 (0.25) & 37.97 (1.26) & 10.62 & 33.06\\
 &  & 8,192 & 0.62 ($<0.01$) & 8.60 (0.18) & 25.26 (0.92) & 7.93 & 22.83\\
 &  & 16,384 & 0.37 ($<0.01$) & 6.81 (0.10) & 19.95 (0.50) & 6.30 & 18.58\\
\midrule
Mixture & $Q_0$ & 2,048 & 5.62 (0.04) & 16.96 (0.20) & 65.02 (1.23) & 16.19 & 58.99\\
 &  & 8,192 & 1.22 ($<0.01$) & 11.32 (0.23) & 45.68 (0.80) & 10.85 & 39.85\\
 &  & 16,384 & 0.72 ($<0.01$) & 9.52 (0.11) & 35.28 (0.86) & 9.27 & 32.73\\
\addlinespace[3pt]
 & $Q_1$ & 2,048 & 7.06 (0.10) & 18.04 (0.25) & 66.64 (1.39) & 16.36 & 57.29\\
 &  & 8,192 & 1.48 (0.01) & 11.62 (0.21) & 47.54 (0.87) & 11.30 & 43.05\\
 &  & 16,384 & 0.77 ($<0.01$) & 10.12 (0.08) & 37.34 (0.80) & 9.44 & 33.93\\
\bottomrule
\end{tabular}
\end{table}

\subsection{Quoted directions and claim prices}
We compare $E=\{(\cos\theta,\sin\theta):0\le\theta\le\pi/2\}$ with $E\cup\{\boldsymbol e\}$, where $\boldsymbol e=(1,-1)/\sqrt2$. The claims are the squared spread $(\boldsymbol e\cdot\boldsymbol x)^2$ and the exchange option $g_+(\boldsymbol x)=(\boldsymbol e\cdot\boldsymbol x)_+$. We take $\mathcal C$ to be the convex hull of $299$ component laws with radial densities $f_{3,\boldsymbol c,aI}$,  $a\in\{0.12,0.24,0.40\}$, and centers on a grid of spacing $0.12$, retaining supports inside $B_{0.92}$.

Table~\ref{tab:numerical-moduli} bounds the uniform-quote modulus $\varpi_g(\epsilon)$ on $\mathcal C$. The dual of a linear program on a direction--strike grid certifies upper bounds through the finite-class form of the replication bound in Lemma~\ref{lem:replication-bound}. For lower bounds, we mix both measures in the LP solution with their midpoint using the same weight. Derivative bounds between grid points determine the weight needed to keep the absolute quote difference at most $\epsilon$ over every quoted direction and strike.

\begin{table}[!htbp]
\centering\small
\caption{Lower and upper bounds for pricing moduli, multiplied by $10^3$.}\label{tab:numerical-moduli}
\begin{tabular}{crrrr}
\toprule
$\epsilon$ & \multicolumn{2}{c}{Squared spread} & \multicolumn{2}{c}{Exchange option}\\
\cmidrule(lr){2-3}\cmidrule(lr){4-5}
& Long-only & With spread & Long-only & With spread\\
\midrule
$10^{-2}$ & $[84.79,86.53]$ & $[41.39,43.72]$ & $[132.4,135.2]$ & $[9.638,10]$\\
$10^{-3}$ & $[8.829,9.008]$ & $[4.39,4.48]$ & $[67.59,68.98]$ & $[0.9488,1]$\\
$10^{-4}$ & $[0.9086,0.9328]$ & $[0.4603,0.4709]$ & $[28.17,28.79]$ & $[0.09417,0.1]$\\
\bottomrule
\end{tabular}
\end{table}

The squared-spread modulus decreases nearly linearly with $\epsilon$ under both designs, as expected from polynomial replication. By contrast, the exchange modulus decreases much more slowly with long-only quotes. Once the spread direction is included, the exchange payoff itself is quoted, and Lemma~\ref{lem:replication-bound} gives the upper bound $\epsilon$. The numerical intervals then decrease linearly with $\epsilon$. 
These computations illustrate the directional effect of Theorem~\ref{thm:log-modulus}, whose exact logarithmic order concerns the full Sobolev class, not the finite mixture class used here.

\section{Conclusion and Open Problems}\label{sec:discussion}
Random basket quotes link distributional recovery and the estimation of individual claim prices through one linear inverse operator. Its negative-Sobolev geometry gives sharp Wasserstein recovery rates on Sobolev density classes, uniformly from exact quotes to bounded noise levels. A bounded systematic pricing bias creates an additional error that does not vanish with more observations. For a specified claim, a modulus of continuity and its static-replication dual characterize minimax estimation at fixed positive noise. Under long-only directions, the joint law is identified, yet the exchange price is recovered only at a logarithmic rate, even from exact quotes.

\par\medskip\noindent\textbf{Arbitrary laws.} For exact quotes without a density assumption, Theorem~\ref{thm:nu-all-laws} leaves a gap between the exponents $1/d$ and $(d+3)/(d(2d+3))$. A possible route to closing it is a pricing-surface pilot with expected squared $L^2(Q_0)$ error bounded by $C(m^{-1-3/d}+\varsigma^2)$ uniformly over $\Pp(B_r)$. If such a pilot is available, Proposition~\ref{prop:nu-transfer}, with $M\asymp\min\{n^{1/d},(n/\varsigma^2)^{1/(2d+3)}\}$, would give the matching rate $n^{-1/d}+(\varsigma^2/n)^{1/(2d+3)}$.

\par\medskip\noindent\textbf{Post-design uniformity.} Theorem~\ref{thm:exact-smooth} controls all Sobolev-restricted quote fibers on a common random-design event, with a logarithmic factor. The minimax expected-risk rate of Theorem~\ref{thm:nu-smooth} has no such factor. Whether it is necessary for the post-design supremum remains open. The two criteria are distinguished in Section~\ref{sec:exact-lower-bounds}.

\par\medskip\noindent\textbf{Adaptation.} The bandwidth of the estimator in Corollary~\ref{cor:nu-general-design} depends on $\alpha,L,\varsigma$, and $\epsilon$. Section~\ref{sec:numerical} uses a Lepski-type selection rule without establishing an adaptive minimax theorem. Proving the theoretical guarantees and allowing heteroskedastic or price-dependent errors are natural next steps.

\putbib
\end{bibunit}

\clearpage
\begin{appendix}
\begin{bibunit}

\section{Proofs for Section~\ref{sec:model}}
We use $H^2(P,P')=\int(\sqrt{\mathrm dP}-\sqrt{\mathrm dP'})^2$ for the squared Hellinger distance.

\begin{lemma}[Hellinger bound for quote experiments]\label{lem:hellinger}
Let $Q$ be a design, $\varsigma>0$, and let $P^{(n)}_\mu$ denote the law of $n$ i.i.d. triples $(Z_i,C_\mu(Z_i)+\varsigma\xi_i)$ with $Z_i\sim Q$. If $\xi_1$ has a density with finite Fisher information $I_\xi$, then
\[
 \operatorname{TV}\left(P^{(n)}_\mu,P^{(n)}_\nu\right)^2\le H^2\left(P^{(n)}_\mu,P^{(n)}_\nu\right)\le\frac{nI_\xi}{4\varsigma^2}\norm{C_\mu-C_\nu}_{L^2(Q)}^2.
\]
\end{lemma}

\begin{proof}
Since $\sqrt{p_\xi}$ is absolutely continuous with $\int((\sqrt{p_\xi})')^2=I_\xi/4$, the Cauchy--Schwarz inequality gives $\int(\sqrt{p_\xi(y-a)}-\sqrt{p_\xi(y)})^2\dd y\le a^2I_\xi/4$ for every shift $a$. Conditioning on $Z_1$ and rescaling by $\varsigma$ bounds the squared Hellinger distance of one triple by $I_\xi\norm{C_\mu-C_\nu}_{L^2(Q)}^2/(4\varsigma^2)$. Since Hellinger affinities multiply over independent coordinates, $H^2$ is subadditive, and $\operatorname{TV}\le H$.
\end{proof}

\begin{proof}[Proof of Theorem~\ref{thm:direction-identification}]
The compact-support argument follows the homogeneous-polynomial criterion of \cite{cuesta2007}. Suppose (ii) holds and $C_\mu=C_\nu$ on $E\times[-R,R]$. The projected measures agree in each direction in $E$, by the strike-extension argument preceding the theorem. For every integer $m\ge0$, the homogeneous polynomial
\[
 p_m(\bu)=\int(\bu\cdot\bx)^m\,(\mu-\nu)(\mathrm d\bx)
\]
vanishes on $E$ and therefore vanishes identically. Expansion in monomials shows that all mixed moments of $\mu-\nu$ vanish. Polynomials are uniformly dense in the continuous functions on $B_R$, by the real Stone--Weierstrass theorem. Thus $\mu=\nu$.

Conversely, suppose a nonzero homogeneous polynomial $p$ vanishes on $E$. Multiplication by a nonzero homogeneous factor, if necessary, ensures that its degree is at least two. Let $g\in C_c^\infty$ be a nonzero real function supported in a small cube inside $B_R$ and define $\eta=p(\partial)g$. The function $\eta$ is not identically zero: its Fourier transform is $p(\mathrm i\bxi)\widehat g(\bxi)$, and a nonzero polynomial cannot vanish on an open set on which $\widehat g$ is nonzero. Integration by parts shows that $\eta$ has zero mass and first moments. For each $\bu\in E$ and $t\in\R$,
\[
 \widehat\eta(t\bu)=\mathrm i^{\deg p}t^{\deg p}
                                      p(\bu)\widehat g(t\bu)=0.
\]
Therefore every corresponding projected signed measure is zero. Choose a smooth probability density $f_0$ positive on $\supp \eta$ and a sufficiently small $\epsilon>0$. The distinct densities $f_0\pm\epsilon \eta$ have identical projected laws, hence identical call prices, in all directions of $E$.

A homogeneous polynomial vanishing on a relatively open sphere patch vanishes on the associated open cone and hence is identically zero. For finite $E=\{\bu_1,\ldots,\bu_m\}$, choose nonzero $\bv_j$ perpendicular to $\bu_j$ and use $p(\bxi)=\prod_j(\bv_j\cdot\bxi)$. This polynomial is nonzero and vanishes on $E$. In dimension two, a nonzero homogeneous polynomial has finitely many zeros on the real projective line, and an infinite set of distinct unoriented directions therefore satisfies (ii).
\end{proof}

\begin{proof}[Proof of Theorem~\ref{thm:set-consistency}]
Choose a countable base for the compact metric space $\supp Q$. Every nonempty relative open set in this base has positive $Q$ probability and is visited almost surely by the i.i.d. sample. Taking the intersection of these probability-one events shows that $\{Z_i:i\ge1\}$ is dense in $\supp Q$. Continuity in \eqref{eq:lipschitz} gives
\[
 \bigcap_{n\ge1}\K_n(\mu)=\K_\infty(\mu;Q).
\]
Suppose the Hausdorff convergence fails on this event. Since $\K_\infty\subseteq\K_n$, there are $\epsilon>0$, integers $n_j\uparrow\infty$, and $\nu_j\in\K_{n_j}(\mu)$ with
$\dist_{W_1}(\nu_j,\K_\infty)\ge\epsilon$.
A subsequence converges in $W_1$ to some $\nu$. For each $m$, all sufficiently late terms belong to the closed set $\K_m(\mu)$, hence $\nu\in\K_m(\mu)$. Thus $\nu\in\K_\infty$, which contradicts the displayed distance bound. Under identification, $\K_\infty=\{\mu\}$.
\end{proof}

\section{Proofs for Section~\ref{sec:geometry}}
\begin{lemma}\label{lem:mean}
If $\eta(B_R)=0$ and $\ba=\int\bx\,\eta(\mathrm d\bx)$, then
\begin{equation}\label{eq:parity}
 \cA\eta(\bu,k)-\cA\eta(-\bu,-k)=\bu\cdot\ba,
 \qquad
 \ba=2d\int_\Z\bu\,\cA\eta(\bu,k)\dd Q_0.
\end{equation}
Moreover, $\abs{\ba}\le2\sqrt d\norm{\cA\eta}_{L^2(Q_0)}$.
\end{lemma}

\begin{proof}
The identity $t_+-(-t)_+=t$ gives the first equality. The transformation $(\bu,k)\mapsto(-\bu,-k)$ preserves $Q_0$, and $\int\bu\bu^{\mathsf T}\dd\sigma=I_d/d$, giving the second. The odd part of $\cA\eta$ under this transformation is $\bu\cdot\ba/2$. Its $L^2(Q_0)$ norm is $\abs{\ba}/(2\sqrt d)$ and cannot exceed the norm of $\cA\eta$.
\end{proof}

Fix a smooth probability density $\varrho$ supported strictly inside $B_R$ and define $\tau_j=-\partial_j\varrho$, $1\le j\le d$. Integration by parts gives
\begin{equation}\label{eq:tau-moments}
 \int\tau_j\dd\bx=0,
 \qquad \int x_l\tau_j(\bx)\dd\bx=\mathbf 1_{\{j=l\}}.
\end{equation}
Let $\eta$ be a finite signed measure on $B_R$ with $\eta(B_R)=0$, and let $\ba=\int\bx\,\eta(\mathrm d\bx)$. Then
\begin{equation}\label{eq:eta0}
 \eta_0=\eta-\sum_{j=1}^d a_j\tau_j\dd\bx
\end{equation}
has zero mass and zero first moments. If $\eta=\mu-\nu$ with $\mu,\nu\in\Pp(B_R)$ and $D=\norm{\cA\eta}_{L^2(Q_0)}$, then
\begin{equation}\label{eq:correction-bounds}
 \norm{\eta_0}_{\TV}\le C_{d,R},\qquad
 \norm{\cA\eta_0}_{L^2(Q_0)}\le C_{d,R}D.
\end{equation}
These follow from the bound $\abs{\ba}\le2\sqrt d\,D$ of Lemma~\ref{lem:mean}, from $D\le2R$, and from the smooth functions $\tau_j$.

\begin{proof}[Proof of Theorem~\ref{thm:operator-identity}]
For each $\bu$, let $\eta_{\bu}$ be the pushforward of $\eta$ under $\bx\mapsto\bu\cdot\bx$. The zero mass and first moment imply $\cA\eta(\bu,k)=0$ for $\abs{k}\ge R$. Its second distributional derivative in $k$ is $\eta_{\bu}$. Fourier transformation in $k$ gives
\begin{equation}\label{eq:slice}
 \widehat{\cA\eta}(\bu,t)=-\frac{\widehat\eta(t\bu)}{t^2},\qquad t\ne0.
\end{equation}
Since $\widehat\eta(\bxi)=O(\abs{\bxi}^2)$ near zero, the expressions below are integrable there. At infinity, boundedness of $\widehat\eta$ gives integrability. Plancherel's identity and polar coordinates yield
\begin{align}
 \norm{\cA\eta}_{L^2(Q_0)}^2
 &=\frac{1}{4\pi R}\int_{S^{d-1}}\int_\R
       \frac{\abs{\widehat\eta(t\bu)}^2}{t^4}\dd t\dd\sigma(\bu)\,=\frac{1}{2\pi R\omega_{d-1}}
       \int_{\R^d}\abs{\widehat\eta(\bxi)}^2\abs{\bxi}^{-d-3}\dd\bxi.
 \label{eq:polar-calculation}
\end{align}

For $a,b\in[-R,R]$, direct integration gives
\begin{align*}
 \int_{-R}^R(a-k)_+(b-k)_+\dd k
 ={}&\frac{\abs{a-b}^3}{12}+\frac{R^3}{3}
       +\frac{R^2}{2}(a+b)+Rab-\frac{a^3+b^3}{12}+\frac{a^2b+ab^2}{4}.
\end{align*}
Substitute $a=\bu\cdot\bx$ and $b=\bu\cdot\by$ and integrate against $\eta\otimes\eta$. Every polynomial term vanishes because each contains either a factor of total mass or a factor of first moment. Thus
\[
 \int_{-R}^R\abs{\cA\eta(\bu,k)}^2\dd k
 =\frac1{12}\iint\abs{\bu\cdot(\bx-\by)}^3
                    \,\eta(\mathrm d\bx)\eta(\mathrm d\by).
\]
Rotation invariance gives $\int\abs{\bu\cdot\bx}^3\dd\sigma=m_{3,d}\abs{\bx}^3$. The value of $m_{3,d}$ follows by integrating the density proportional to $(1-t^2)^{(d-3)/2}$ of the first coordinate of a uniform direction. Division by $2R$ proves \eqref{eq:cubic-identity}.
\end{proof}\begin{proof}[Proof of Theorem~\ref{thm:modulus}]
Let $\eta=\mu-\nu$ and form $\eta_0$ as in \eqref{eq:eta0}. A 1-Lipschitz function on $B_R$, normalized to vanish at the origin, has a compactly supported Lipschitz extension $F$ with
\begin{equation}\label{eq:lip-extension}
 \norm{F}_{H^1}+\norm{F}_\infty+\Lip(F)\le C_{d,R}.
\end{equation}
Let $F_\epsilon=F*\phi_\epsilon$, where $\phi$ is a smooth compactly supported probability density. Then
\[
 \norm{F-F_\epsilon}_\infty\le C\epsilon,
 \qquad
 \norm{F_\epsilon}_{\dot H^s}\le C\epsilon^{-(s-1)}\norm{F}_{H^1}.
\]
For the second inequality, use
$\sup_{\bxi}\abs{\bxi}^{s-1}\abs{\widehat\phi(\epsilon\bxi)}
\le C\epsilon^{-(s-1)}$.
Theorem~\ref{thm:operator-identity} and \eqref{eq:correction-bounds} give
\begin{align*}
 \abs{\int F\dd\eta}
 &\le\norm{\eta_0}_{\TV}\norm{F-F_\epsilon}_\infty
    +\norm{\eta_0}_{\dot H^{-s}}\norm{F_\epsilon}_{\dot H^s}
    +\sum_j\abs{a_j}\abs{\int F\tau_j\dd\bx}\\
 &\le C\left(\epsilon+D\epsilon^{-(s-1)}+D\right).
\end{align*}
Choose $\epsilon=D^{1/s}$ for $0<D\le1$. The cases $D=0$ and $D>1$ follow by a limit and the bound $W_1\le2R$, respectively. Taking the supremum over Lipschitz tests proves \eqref{eq:universal-modulus}.

For the density assertion, $\norm{\eta_0}_{H^\alpha}\le C(1+L)$ and
$\norm{\eta_0}_{H^{-s}}\le\norm{\eta_0}_{\dot H^{-s}}\le CD$.
Weighted H\"older's inequality in the Fourier integral gives
\begin{equation}\label{eq:interpolation}
 \norm{\eta_0}_{H^{-1}}
 \le\norm{\eta_0}_{H^{-s}}^{(\alpha+1)/(\alpha+s)}
    \norm{\eta_0}_{H^\alpha}^{(s-1)/(\alpha+s)}.
\end{equation}
Pairing with \eqref{eq:lip-extension} and adding the mean correction proves \eqref{eq:smooth-modulus}.

For sharpness, take a nonnegative $\zeta_0\in C_c^\infty(\R^d)$ supported in a small cube inside $B_R$ with positive integral. Choose a smooth probability density $f_0$ that is bounded below by a positive constant on that support. For integers $N\ge1$, define
\begin{equation}\label{eq:oscillation}
 H_N(\bx)=bN^{-\alpha-2}\partial_1^2
                    \left(\zeta_0(\bx)\cos(Nx_1)\right),
 \qquad f_N^\pm=f_0\pm H_N.
\end{equation}
For a sufficiently small $b>0$, these are nonnegative probability densities. Integration by parts gives zero mass and first moments for $H_N$. The product rule gives $\norm{H_N}_\infty\le CbN^{-\alpha}$, and Fourier translation gives $\norm{H_N}_{H^\alpha}\le Cb$.

The Fourier transform is
\[
 \widehat H_N(\bxi)
 =-\frac b2 N^{-\alpha-2}\xi_1^2
 \left(\widehat\zeta_0(\bxi-N\be_1)
                 +\widehat\zeta_0(\bxi+N\be_1)\right).
\]
On $\abs{\bxi}\ge N/2$, the weighted integral in \eqref{eq:fourier-identity} is bounded by
$Cb^2N^{-2\alpha-2s}$, since $4-2s=1-d\le0$ and the translated transforms have uniformly bounded $L^2$ norms. On $\abs{\bxi}<N/2$, Schwartz decay bounds each translated transform by $C_PN^{-P}$ for every $P$; the radial integral of $\abs{\bxi}^{4-2s}$ over this ball is $O(N)$. Choosing $P\ge d/2$ gives the same upper bound. On a small ball about $N\be_1$, continuity and $\widehat\zeta_0(0)>0$ give the matching lower bound. Therefore
\begin{equation}\label{eq:oscillation-rates}
 \norm{\cA H_N}_{L^2(Q_0)}\asymp N^{-\alpha-s}.
\end{equation}
The function $F_N(\bx)=N^{-1}\cos(Nx_1)$ is 1-Lipschitz, and
\[
 \abs{\int F_NH_N\dd\bx}
 =bN^{-\alpha-1}\int\zeta_0(\bx)\cos^2(Nx_1)\dd\bx
 \ge cN^{-\alpha-1}
\]
for all sufficiently large $N$. A single integration by parts against any Lipschitz test also gives the upper bound $C N^{-\alpha-1}$. Thus
\[
 W_1(f_N^+\dd\bx,f_N^-\dd\bx)\asymp N^{-\alpha-1}.
\]
Together with \eqref{eq:oscillation-rates}, this prevents an increase of the exponent in \eqref{eq:smooth-modulus}. The case $\alpha=0$ proves sharpness of \eqref{eq:universal-modulus}.
\end{proof}

\begin{proof}[Proof of Lemma~\ref{lem:sobolev-map}]
Choose the functions in \eqref{eq:tau-moments} supported strictly inside $B_r$. For $\eta_0=\eta-\sum_j(\int x_j\eta\dd\bx)\tau_j$, both mass and first moments vanish. The moment functionals are bounded on functions supported on $B_r$ in every Sobolev norm, by pairing with smooth cutoffs times $1,x_1,\ldots,x_d$. The finite sum removed from $\eta$ has a smooth pricing surface on $\Z$. It suffices to estimate $\cA \eta_0$, extended by zero outside $[-r,r]$ in the strike coordinate.

Begin with a smooth compactly supported function $v$ with zero mass and first moments. Equation~\eqref{eq:slice} applies, and the estimates below extend the operator by density. Let $m$ be a nonnegative integer. The product Sobolev $H^m$ norm on $S^{d-1}\times\R$ is bounded by the $L^2$ norms of derivatives obtained from at most $m$ strike derivatives and smooth angular vector fields. Applying $q$ strike derivatives and $j$ angular derivatives to \eqref{eq:slice}, with $q+j\le m$, gives a finite sum of terms
\begin{equation}\label{eq:angular-terms}
 b(\bu)t^{q+l-2}\widehat{p_l v}(t\bu),
 \qquad 0\le l\le j,
\end{equation}
where $b$ is a smooth function on the sphere and $p_l$ is a polynomial of degree at most $l$. Terms arising from differentiation of angular coefficients have the same form with a smaller $l$.
On $\abs{t}\ge1$, polar coordinates bound the squared $L^2$ norm of \eqref{eq:angular-terms} by
\begin{align}
 C\int_{\abs{\bxi}\ge1}
    \abs{\widehat{p_l v}(\bxi)}^2
               \abs{\bxi}^{2(q+l)-d-3}\dd\bxi
 &\le C\norm{p_l v}_{H^{m-s}}^2\,\le C\norm{v}_{H^{m-s}}^2.
 \label{eq:high-frequency-map}
\end{align}
For the last step, multiply $p_l$ by a smooth cutoff equal to one on the common compact support; multiplication by this function is bounded on each Sobolev space in question.

On $\abs{t}\le1$, use moment cancellation before differentiating:
\begin{equation}\label{eq:low-frequency-map}
 \frac{\widehat v(t\bu)}{t^2}
 =-\int_0^1(1-a)
      \ip{v}{(\bu\cdot\bx)^2 e^{-\mathrm iat\bu\cdot\bx}}\dd a.
\end{equation}
A cutoff on the compact support is understood in this pairing. Every required angular derivative of the test function, for $\abs{t}\le1$, is bounded in $H^{s-m}$ uniformly in $a,t,\bu$. Strike derivatives in the original variable multiply by $t^q$ after Fourier transformation. Hence the low-frequency contribution to the $H^m$ norm is bounded by $C\norm{v}_{H^{m-s}}$ as well. Combining the two regions yields
\begin{equation}\label{eq:integer-map}
 \norm{\cA v}_{H^m(S^{d-1}\times\R)}
          \le C_m\norm{v}_{H^{m-s}(\R^d)}.
\end{equation}

Interpolation is carried out on unconstrained spaces. Choose a smooth cutoff $\chi$ with compact support and equal to one on $B_r$, and a smooth function $\tau_0$ with mass one and mean zero. Define the bounded finite-rank correction
\[
 \Pi w=w-\left(\int w\dd\bx\right)\tau_0
                 -\sum_j\left(\int x_jw\dd\bx\right)\tau_j
\]
for compactly supported distributions $w$. Applying \eqref{eq:integer-map} to the operator $v\mapsto\cA\Pi(\chi v)$ gives bounded maps from $H^{m-s}(\R^d)$ to $H^m(S^{d-1}\times\R)$ at every integer $m\ge0$. Interpolation between neighboring Hilbert Sobolev spaces gives the map from $H^\alpha$ to $H^{\alpha+s}$. For an input in $H^\alpha$ with $\alpha\ge0$, approximation on a compact support in $L^2$ also gives convergence in $L^1$, and hence uniform convergence of the call integrals. Thus the extended operator agrees with the operator in \eqref{eq:call}. Adding the smooth mean correction proves \eqref{eq:forward-sobolev}.

For \eqref{eq:inverse-sobolev}, let $\eta=f-g$, $\ba=\int\bx \eta\dd\bx$, and $\eta_0=\eta-\sum_j a_j\tau_j$. From \eqref{eq:slice} and polar coordinates,
\begin{equation}\label{eq:hminusone}
 \int_{S^{d-1}}\int_\R
       \abs{t}^{2(s-1)}\abs{\widehat{\cA \eta_0}(\bu,t)}^2
                          \dd t\dd\sigma(\bu)
 =\frac{2}{\omega_{d-1}}
     \int_{\R^d}\abs{\widehat \eta_0(\bxi)}^2\abs{\bxi}^{-2}\dd\bxi.
\end{equation}
The left side is bounded by the squared $H^{s-1}$ norm on the full cylinder. Since $\cA \eta_0$ vanishes for $\abs{k}\ge r$, multiplying a bounded extension from $\Z$ by a smooth strike cutoff supported in $(-R,R)$ and equal to one on $[-r,r]$ shows that this full-cylinder norm is bounded by $C\norm{\cA \eta_0}_{H^{s-1}(\Z)}$. Lemma~\ref{lem:mean} gives
\[
 \abs{\ba}\le C\norm{\cA \eta}_{L^2(Q_0)},\qquad
 \norm{\cA \eta_0}_{H^{s-1}(\Z)}
       \le C\norm{\cA \eta}_{H^{s-1}(\Z)}.
\]
Pairing $\eta_0$ with the extension in \eqref{eq:lip-extension}, using \eqref{eq:hminusone} and $\norm{\eta_0}_{H^{-1}}\le\norm{\eta_0}_{\dot H^{-1}}$, bounds every Lipschitz test by $C\norm{\cA \eta}_{H^{s-1}(\Z)}$. The mean correction has the same bound. Taking the supremum proves \eqref{eq:inverse-sobolev}.
\end{proof}

Set $a_{d,R}=R\omega_{d-1}(2\pi)^{1-d}$.
\begin{proof}[Proof of Lemma~\ref{lem:dual-features}]
Define $G_\psi$ on the full cylinder by its Fourier transform in $k$:
\begin{equation}\label{eq:dual-feature-fourier}
 \widehat G_\psi(\bu,t)
       =-a_{d,R}\abs{t}^{d+1}\widehat\psi(t\bu).
\end{equation}
Plancherel and polar coordinates give
\begin{equation}\label{eq:dual-feature-norm}
 \int_{S^{d-1}}\int_\R\abs{G_\psi(\bu,k)}^2
                   \frac{\dd k}{2R}\dd\sigma(\bu)
       =a_{d,R}\norm{\psi}_{\dot H^s(\R^d)}^2.
\end{equation}
For every finite signed measure $\eta$ with zero mass and first moments supported on $B_R$, \eqref{eq:slice}, \eqref{eq:dual-feature-fourier}, and the same polar calculation show
\begin{equation}\label{eq:dual-zero-moments}
 \int\psi\dd\eta=\int_\Z G_\psi\cA\eta\dd Q_0.
\end{equation}
The full-cylinder integral reduces to $\Z$, since $\cA\eta$ vanishes outside $[-R,R]$. The full-cylinder Parseval integral equals
\[
 \frac{a_{d,R}}{4\pi R}
 \int_{S^{d-1}}\int_\R\abs{t}^{d-1}
          \widehat\psi(t\bu)\widehat\eta(-t\bu)\dd t\dd\sigma(\bu)
 =(2\pi)^{-d}\int\widehat\psi(\bxi)\widehat\eta(-\bxi)\dd\bxi,
\]
which is the required pairing. The Fourier integral converges absolutely since $\widehat\psi$ is rapidly decreasing and $\widehat\eta$ is bounded. The pairing of $G_\psi$ and $\cA\eta$ is a Plancherel pairing in $L^2$ of the full cylinder, since $G_\psi\in L^2$ by \eqref{eq:dual-feature-norm} and $\cA\eta\in L^2$ by \eqref{eq:fourier-identity}.

Let $\bm=\int\bx\,\mu(\mathrm d\bx)$ and use $\eta=\mu-\delta_0-\sum_jm_j\tau_j\dd\bx$ in \eqref{eq:dual-zero-moments}. Define
\begin{equation}\label{eq:dual-corrections}
 \begin{split}
 b_{\psi,j}&=\int\psi\tau_j\dd\bx
                     -\int_\Z G_\psi\cA\tau_j\dd Q_0,\\
 a_\psi&=\psi(0)-\int_\Z G_\psi C_{\delta_0}\dd Q_0.
 \end{split}
\end{equation}
Then $\int\psi\dd\mu=a_\psi+\int G_\psi C_\mu\dd Q_0+\bm\cdot\boldsymbol b_\psi$.
For a probability measure, the parity identity reads
$C_\mu(\bu,k)-C_\mu(-\bu,-k)=\bu\cdot\bm-k$.
After multiplication by $\bu$ and integration, the strike term vanishes and
\[
 \bm=2d\int_\Z\bu C_\mu(\bu,k)\dd Q_0.
\]
Thus \eqref{eq:dual-functional} holds with
\begin{equation}\label{eq:tilde-feature}
 \widetilde G_\psi(\bu,k)=G_\psi(\bu,k)
                                  +2d\boldsymbol b_\psi\cdot\bu.
\end{equation}
Equation~\eqref{eq:dual-feature-norm}, Cauchy--Schwarz in \eqref{eq:dual-corrections}, and the smoothness of $\tau_j$ bound $\norm{\widetilde G_\psi}_2$ by $C\norm{\psi}_{H^s}$. Finally, $s>d/2$ gives $\abs{\psi(0)}\le C\norm{\psi}_{H^s}$, proving the bound for $a_\psi$.
\end{proof}

\section{Proofs for Section~\ref{sec:exact}}\label{AppendixC}
\begin{lemma}[Random coverage]\label{lem:coverage}
Assume $Q\ge q_*Q_0$ as measures, where $q_*>0$. There are constants $c,C>0$ depending on $d$ such that
\begin{equation}\label{eq:fill-tail}
 \Pr(h_n>2t)\le Ct^{-d}\exp\left(-cq_*nt^d\right),\qquad 0<t<1.
\end{equation}
Consequently, for every $0<\delta<1$,
\begin{equation}\label{eq:fill-highprob}
 \Pr\left(h_n>C_{d,q_*}
         \left(\frac{\log(en/\delta)}{n}\right)^{\frac{1}{d}}\right)\le\delta,
\end{equation}
with a constant enlarged to include the cases where the displayed threshold is not small. Also $h_n=O((\log n/n)^{1/d})$ almost surely.
\end{lemma}\begin{proof}
Let $w_1,\ldots,w_M$ be a $t$-net with $M\le Ct^{-d}$. If every ball $B_\rho(w_j,t)$ contains a sampled point, then $h_n\le2t$. For each $j$, its probability of being empty is at most
$\left(1-cq_*t^d\right)^n\le\exp(-cq_*nt^d)$.
The union bound gives \eqref{eq:fill-tail}. Taking $t^d$ to be a sufficiently large multiple of $\log(en/\delta)/n$ gives \eqref{eq:fill-highprob}. Taking $\delta=n^{-2}$ and applying the Borel--Cantelli lemma gives the almost-sure assertion.
\end{proof}\begin{proof}[Proof of Proposition~\ref{prop:no-uniform-rate}]
With probability $(1-\epsilon)^n$ all $n$ locations equal $z_0=(\bu_0,k_0)$, and the data are $n$ copies of $(z_0,C_\mu(z_0))$. Let $\boldsymbol w$ be a unit vector orthogonal to $\bu_0$, and put $\bx=r\boldsymbol w$ and $\by=-r\boldsymbol w$. The point masses $\delta_{\bx}$ and $\delta_{\by}$ have the same call price $(-k_0)_+$ at $z_0$ and $W_1(\delta_{\bx},\delta_{\by})=2r$. The event depends only on the design, and on it the two laws produce identical data. The triangle inequality gives $\E_{\bx}W_1(\widehat\mu,\delta_{\bx})+\E_{\by}W_1(\widehat\mu,\delta_{\by})\ge2r(1-\epsilon)^n$, and the larger of the two risks is at least $r(1-\epsilon)^n$.
\end{proof}

\begin{proof}[Proof of Theorem~\ref{thm:exact-general}]
If $\mu,\nu\in\K_n(z,v)$, the function $C_\mu-C_\nu$ is zero at all $z_i$ and has $\rho$-Lipschitz constant at most $2R$. Therefore
\[
 \norm{C_\mu-C_\nu}_{L^2(Q_0)}
 \le\norm{C_\mu-C_\nu}_\infty\le2Rh(z).
\]
Apply \eqref{eq:universal-modulus} and take the supremum over the pair. The probability event in Lemma~\ref{lem:coverage} depends only on the design, giving simultaneous validity over the price vectors.
\end{proof}\begin{lemma}[Sampling inequalities]\label{lem:zeros}
Let $t>d/2$ and $0\le q\le t$. There are $h_*>0$ and $C>0$ such that every finite set $z\subset\Z$ with fill distance $h\le h_*$ and every $G\in H^t(\Z)$ satisfy
\begin{equation}\label{eq:sampling-nonzero}
 \norm{G}_{L^2(\Z)}\le C\left(h^t\norm{G}_{H^t(\Z)}+\max_{z_i\in z}\abs{G(z_i)}\right).
\end{equation}
If $G$ vanishes on $z$, then
\begin{equation}\label{eq:zeros}
 \norm{G}_{H^q(\Z)}\le C h^{t-q}\norm{G}_{H^t(\Z)}.
\end{equation}
\end{lemma}

\begin{proof}
The Euclidean $L^2$ sampling inequality on a bounded Lipschitz domain, in the case $p=q=2$, $l=0$ with discrete exponent $\infty$ of Theorem~1.1 in Corrigan, Wallin, and Wanner \cite{corrigan2009}, which restates the sampling theorem of Arcang\'eli and coauthors, gives
$\norm{G}_{L^2}\le C(h^t\norm{G}_{H^t}+\max_i\abs{G(z_i)})$ for every real $t>d/2$.

To transfer this estimate to $\Z$, use a finite collection of sphere charts with smaller subcharts covering the sphere and smooth cutoffs supported in the larger charts. In each chart, apply the Euclidean result to the cutoff times $G$. At original sample points within the chart, this product is bounded by $\max_i\abs{G(z_i)}$. In the portion outside the cutoff support, add a finite $h$-net of points, where the product is zero. The positive distance between the smaller and larger chart boundaries ensures, for sufficiently small $h$, a chart fill distance at most $Ch$. The strike endpoints are covered by Lipschitz half-boxes, to which the same Euclidean result applies. Bounded coordinate maps and the finite partition of unity give \eqref{eq:sampling-nonzero}. If $G$ vanishes on $z$, its right side is $Ch^t\norm{G}_{H^t(\Z)}$, and Hilbert interpolation gives
\[
 \norm{G}_{H^q(\Z)}
 \le C\norm{G}_{L^2(\Z)}^{1-q/t}\norm{G}_{H^t(\Z)}^{q/t}
 \le Ch^{t-q}\norm{G}_{H^t(\Z)}.
\]
\end{proof}\begin{proof}[Proof of Theorem~\ref{thm:exact-smooth}]
For two feasible densities, $G=\cA(f-g)$ vanishes at every sampled quote and, by \eqref{eq:forward-sobolev}, satisfies
$\norm{G}_{H^{\alpha+s}(\Z)}\le2CL$.
Since $\alpha+s>d/2$, Lemma~\ref{lem:zeros} applies with $t=\alpha+s$ and $q=s-1$. Lemma~\ref{lem:sobolev-map} then gives
\[
 W_1(f\dd\bx,g\dd\bx)
 \le C\norm{G}_{H^{s-1}(\Z)}
 \le Ch(z)^{\alpha+1}\norm{G}_{H^{\alpha+s}(\Z)}
 \le CLh(z)^{\alpha+1}.
\]
Lemma~\ref{lem:coverage} gives \eqref{eq:smooth-exact-random}. When its fill-distance threshold exceeds $h_*$, use $W_1\le2r$ and the positive lower bound on the $L^2$ norm of a probability density supported on $B_r$ to enlarge the constant.
\end{proof}

\begin{proof}[Proof of Proposition~\ref{prop:bounded-errors}]
For two feasible densities $f,g$, write $G=\cA(f-g)$ and $t=\alpha+s$. Lemma~\ref{lem:sobolev-map} gives $\norm G_{H^t(\Z)}\le CL$, and $\max_i\abs{G(z_i)}\le2\epsilon$. Inequality \eqref{eq:sampling-nonzero} gives $\norm G_{L^2(\Z)}\le C(Lh^t+\epsilon)$. Interpolating between orders $0$ and $t$ at order $s-1$, and using $1-(s-1)/t=\vartheta$ and $t\vartheta=\alpha+1$,
\[
 \norm G_{H^{s-1}(\Z)}\le C\norm G_{L^2(\Z)}^\vartheta\norm G_{H^t(\Z)}^{1-\vartheta}\le C(Lh^t+\epsilon)^\vartheta L^{1-\vartheta}\le C\left[Lh^{\alpha+1}+L^{1-\vartheta}\epsilon^\vartheta\right].
\]
Inequality \eqref{eq:inverse-sobolev} proves \eqref{eq:bounded-errors}, and Lemma~\ref{lem:coverage} gives the random version.
\end{proof}

\begin{corollary}[Simultaneous control and stopping]\label{cor:stopping}
Under the assumptions of Theorem~\ref{thm:exact-smooth}, with probability at least $1-\delta$, its diameter bound holds for every $n\ge1$ simultaneously with right side
\begin{equation}\label{eq:time-uniform}
 CL\left(\frac{\log(e n^3/\delta)}{n}\right)^{\frac{\alpha+1}{d}}.
\end{equation}
It therefore remains valid at any almost surely finite stopping time based on the observed triples. For $0<t<h_*$, the coverage stopping time
$\tau_t=\inf\{n\ge1:h_n\le t\}$ satisfies
\begin{equation}\label{eq:coverage-stopping}
 \E\tau_t\le C_{d,q_*}t^{-d}\log(e/t).
\end{equation}
At $\tau_t$, every Sobolev-restricted quote fiber has $W_1$ diameter at most $CLt^{\alpha+1}$.
\end{corollary}

\begin{proof}[Proof of Corollary~\ref{cor:stopping}]
Apply Theorem~\ref{thm:exact-smooth} at sample size $n$ with failure probability $6\delta/(\pi^2n^2)$ and take a union bound. Since the resulting event holds for all sample sizes and price vectors, it also holds at the stopping time. Monotonicity of $h_n$ and \eqref{eq:fill-tail} give
\[
 \Pr(\tau_t>n)\le\min\left(1,Ct^{-d}e^{-cq_*nt^d}\right).
\]
Summing this bound over $n\ge0$, splitting the sum at a constant multiple of $t^{-d}\log(e/t)$, proves \eqref{eq:coverage-stopping}. The last assertion is the deterministic bound \eqref{eq:smooth-exact-diameter}.
\end{proof}

For the minimum-norm interpolant \eqref{eq:minnorm-estimator}, if the true density is in $\F_\alpha(L)$, the feasible set is nonempty and the minimizer has norm at most $L$. A minimizing sequence is bounded in the Hilbert space $H^\alpha$, and its weak limit preserves support, nonnegativity, mass, and the finite integral constraints. Lower semicontinuity gives existence and strict convexity gives uniqueness. The minimizer is a measurable function of the observations, by the measurable minimum theorem applied to the weakly compact metrizable norm sublevels and their countable union.
\begin{lemma}[Localized perturbations]\label{lem:bumps}
Fix $\alpha\ge0$. There are a smooth probability density $f_0$ supported on $B_r$, a nonzero smooth function $\psi$ supported in the unit cube with zero mass and first moments, and constants $b_0,C,c>0$ with the following property. For each sufficiently small $\Delta>0$, one can place $M\asymp \Delta^{-d}$ disjoint translated cubes of side $\Delta$ in a region where $f_0$ is bounded below. For any $\btheta\in[-1,1]^M$ and $0<b\le b_0$, define
\begin{equation}\label{eq:bump-family}
 H_{\btheta}(\bx)=b \Delta^\alpha\sum_{j=1}^M \theta_j
                  \psi\left(\frac{\bx-\bx_j}{\Delta}\right).
\end{equation}
Then $f_{\btheta}=f_0+H_{\btheta}$ is a probability density with the one-dimensional coordinate marginals of $f_0$, $\norm{H_{\btheta}}_{H^\alpha}\le Cb$, and
\begin{equation}\label{eq:bump-w1}
 W_1\left(f_{\btheta}\dd\bx,f_{\btheta'}\dd\bx\right)
       \ge cb \Delta^{\alpha+d+1}\norm{\btheta-\btheta'}_1,
 \qquad\btheta,\btheta'\in[-1,1]^M.
\end{equation}
\end{lemma}

\begin{proof}
Take $\psi=\partial_1\partial_2\zeta$ for a nonzero smooth bump $\zeta$ compactly supported inside the unit cube. Integration by parts gives its zero mass and first moments, and integrating out $x_1$ or $x_2$ shows that every one-dimensional coordinate marginal of $\psi$ vanishes. Choose a smooth function $\varphi$ compactly supported in the same cube with $\Lip(\varphi)\le1$ and $\int\varphi\psi\dd\bx>0$. Such a function is obtained by taking a sufficiently small positive multiple of $\psi$ itself. Place the cubes with a fixed proportion of their side length between the supports. Choose $f_0$ smooth and positive on the containing cube. Small $b_0$ guarantees positivity, and the zero mass of every bump preserves total mass.

If $k$ is an integer with $k\ge\alpha$ and $k\ge1$, disjoint supports give
\[
 \norm{H_{\btheta}}_{L^2}\le Cb \Delta^\alpha,
 \qquad \norm{H_{\btheta}}_{H^k}\le Cb \Delta^{\alpha-k}.
\]
Indeed, the squared norm of a derivative of order $l\le k$ is bounded by $Cb^2 \Delta^{2\alpha-2l+d}\sum_j \theta_j^2\le Cb^2\Delta^{2\alpha-2l}$. Fourier interpolation yields $\norm{H_{\btheta}}_{H^\alpha}\le Cb$; for $\alpha=0$ the $L^2$ bound suffices.

The smooth function
\[
 F(\bx)=\Delta\sum_j\operatorname{sgn}(\theta_j-\theta_j')
                  \varphi\left(\frac{\bx-\bx_j}{\Delta}\right)
\]
has gradient norm at most one everywhere, since the summands have disjoint supports. Its integral against $H_{\btheta}-H_{\btheta'}$ is
$b \Delta^{\alpha+d+1}(\int\varphi\psi)\norm{\btheta-\btheta'}_1$.
\end{proof}\begin{proof}[Proof of Theorem~\ref{thm:exact-lower}]
Use Lemma~\ref{lem:bumps} with a grid containing $M$ cells, where $2n\le M\le C_dn$ and $\Delta\asymp n^{-1/d}$. Set
\[
 T_{ij}=\int\phi_i(\bx)
           \psi\left(\frac{\bx-\bx_j}{\Delta}\right)\dd\bx,
 \qquad S=\ker(T)\subseteq\R^M.
\]
Then $\dim S\ge M-n$. Choose a vertex $\btheta$ of the compact polytope $S\cap[-1,1]^M$. At least $\dim S$ coordinates of $\btheta$ have absolute value one. Otherwise, there is a nonzero vector in $S$ whose coordinates vanish at all the saturated coordinates of $\btheta$; sufficiently small perturbations in both signs of this vector stay in the polytope, contradicting the vertex property. Hence
\[
 T\btheta=0,\qquad \sum_j\abs{\theta_j}\ge M-n\ge M/2.
\]
Let $L_0=\norm{f_0}_{H^\alpha}$, and choose $b>0$ sufficiently small that the positivity condition holds and $L_0+Cb\le L$. The densities $f^\pm=f_0\pm H_{\btheta}$ are in $\F_\alpha(L)$, and $T\btheta=0$ gives all the quote equalities. Equation~\eqref{eq:bump-w1} yields
\[
 W_1(f^+\dd\bx,f^-\dd\bx)
 \ge cb\Delta^{\alpha+d+1}(M-n)
 \ge c'\Delta^{\alpha+1}\ge c''n^{-(\alpha+1)/d}.
\]
The constants can be reduced to include any finitely many small $n$ outside the initial grid construction.
\end{proof}

\begin{lemma}[Bayesian bound for exact linear observations]\label{lem:prior-nullspace}
Let $T\in\R^{n\times M}$ be deterministic. Let $\boldsymbol\Theta$ have independent coordinates with density
\begin{equation}\label{eq:compact-prior}
 p_0(t)=\cos^2\left(\frac{\pi t}{2}\right)\mathbf1_{\{-1<t<1\}}.
\end{equation}
For every measurable map $\widehat{\btheta}:\R^n\to[-1,1]^M$, based on the exact observation $T\boldsymbol\Theta$,
\begin{equation}\label{eq:prior-lower}
 \E\norm{\widehat{\btheta}(T\boldsymbol\Theta)-\boldsymbol\Theta}_2^2
       \ge\frac{M-\rank(T)}{\pi^2},\qquad
 \E\norm{\widehat{\btheta}(T\boldsymbol\Theta)-\boldsymbol\Theta}_1
       \ge\frac{M-\rank(T)}{2\pi^2}.
\end{equation}
\end{lemma}

\begin{proof}
The function $p_0$ integrates to one, vanishes at its support endpoints, and satisfies
\[
 \int_{-1}^1p_0'(t)\dd t=0,\qquad
 \int_{-1}^1\frac{p_0'(t)^2}{p_0(t)}\dd t
   =\pi^2\int_{-1}^1\sin^2\left(\frac{\pi t}{2}\right)\dd t
   =\pi^2.
\]
Write $p(\btheta)=\prod_{j=1}^M p_0(\theta_j)$ and
$\boldsymbol S(\btheta)=\nabla\log p(\btheta)$ on the open cube. Independence and the preceding integrals give
\begin{equation}\label{eq:prior-score}
 \E\boldsymbol S(\boldsymbol\Theta)=0,\qquad
 \E\left[\boldsymbol S(\boldsymbol\Theta)
              \boldsymbol S(\boldsymbol\Theta)^{\mathsf T}\right]=\pi^2 I_M.
\end{equation}
Let $\boldsymbol e(\btheta)=\widehat{\btheta}(T\btheta)-\btheta$. For any unit vector $\bv\in\ker(T)$ and any coordinate $j$, Fubini's theorem along lines parallel to $\bv$ yields
\begin{equation}\label{eq:prior-invariance}
 \int_{\R^M}\widehat\theta_j(T\btheta)\,\partial_{\bv}p(\btheta)\dd\btheta=0.
\end{equation}
Indeed, $\widehat\theta_j(T(\btheta+t\bv))$ is constant in $t$, whereas the integral of $\partial_t p(\btheta+t\bv)$ on each such line is zero. The product density is continuously differentiable with compact support, and the decoder is bounded, which justifies the integration. Integration by parts also gives
\[
 \int_{\R^M}\theta_j\partial_{\bv}p(\btheta)\dd\btheta=-v_j.
\]
Subtracting and using $\partial_{\bv}p=p\,\bv\cdot\boldsymbol S$ proves
\begin{equation}\label{eq:prior-score-identity}
 \E\left[e_j(\boldsymbol\Theta)\,
                       \bv\cdot\boldsymbol S(\boldsymbol\Theta)\right]=v_j.
\end{equation}
Let $P$ be the orthogonal projection onto $\ker(T)$ and let $k=\operatorname{tr}(P)=M-\rank(T)$. Sum \eqref{eq:prior-score-identity} over an orthonormal basis of the kernel, after taking its inner product with the corresponding basis vector. Equations~\eqref{eq:prior-score} and \eqref{eq:prior-score-identity} give
\[
 \E\left[(P\boldsymbol e)\cdot(P\boldsymbol S)\right]=k,
 \qquad \E\norm{P\boldsymbol S}_2^2=\pi^2 k.
\]
For $k>0$, Cauchy--Schwarz therefore implies
\[
 k^2\le \E\norm{P\boldsymbol e}_2^2\,\pi^2 k
       \le\E\norm{\boldsymbol e}_2^2\,\pi^2 k.
\]
The case $k=0$ is immediate. Since both the coefficient vector and the decoder lie in $[-1,1]^M$, each coordinate error has absolute value at most two. Thus
$\norm{\boldsymbol e}_2^2\le2\norm{\boldsymbol e}_1$, proving both inequalities in \eqref{eq:prior-lower}.
\end{proof}

\begin{proof}[Proof of Theorem~\ref{thm:exact-information}]
Draw $\boldsymbol\Theta$ from the product prior \eqref{eq:compact-prior}, independently of the design and of the randomization. Fix $\varepsilon>0$ and choose a measurable $\widehat{\btheta}_\varepsilon\in[-1,1]^M$ with $\mathsf d(\widehat f,f_{\widehat{\btheta}_\varepsilon})\le\inf_{\btheta}\mathsf d(\widehat f,f_{\btheta})+\varepsilon$; such a choice follows by enumerating a countable dense subset of the cube and selecting the first coefficient vector within $\varepsilon$ of the infimum over that subset, the objective being continuous in $\btheta$ and measurable in the data. The triangle inequality and the separation hypothesis give
\[
 \delta\norm{\widehat{\btheta}_\varepsilon-\boldsymbol\Theta}_1
   \le\mathsf d(f_{\widehat{\btheta}_\varepsilon},f_{\boldsymbol\Theta})
   \le2\,\mathsf d(\widehat f,f_{\boldsymbol\Theta})+\varepsilon.
\]
Conditional on $(\mathsf L_1,\ldots,\mathsf L_n,a)$ and the randomization, the data are $a'+T\boldsymbol\Theta$ with a fixed vector $a'\in\R^n$ and a fixed matrix $T\in\R^{n\times M}$ of rank at most $n$, and $\widehat{\btheta}_\varepsilon$ is a measurable function of $T\boldsymbol\Theta$. Lemma~\ref{lem:prior-nullspace} gives a conditional expectation of $\norm{\widehat{\btheta}_\varepsilon-\boldsymbol\Theta}_1$ of at least $(M-n)/(2\pi^2)$, since the conditional prior is still \eqref{eq:compact-prior}. Taking expectations, letting $\varepsilon\downarrow0$, and bounding the prior average by the supremum proves \eqref{eq:exact-information}.
\end{proof}

\begin{theorem}[Exact random-design minimax lower bound]\label{thm:random-exact-lower}
Fix $\alpha\ge0$ and $L>L_0$ as in Theorem~\ref{thm:exact-lower}. There is $c>0$, independent of $\mathsf Q_n$, such that
\begin{equation}\label{eq:random-exact-lower-app}
 \mathfrak R_n^{\mathrm{ex}}(\mathsf Q_n,\alpha,L)
           \ge c n^{-(\alpha+1)/d}.
\end{equation}
For an i.i.d. design with $Q\ge q_*Q_0$, there is also $C<\infty$ such that
\begin{equation}\label{eq:random-exact-sandwich}
 c n^{-\frac{\alpha+1}{d}}
 \le\mathfrak R_n^{\mathrm{ex}}(Q^{\otimes n},\alpha,L)
 \le C\left(\frac{\log(en)}{n}\right)^{\frac{\alpha+1}{d}}.
\end{equation}
The lower assertion remains valid for any $n$ continuous linear payoffs whose joint selection is exogenous.
\end{theorem}

\begin{proof}[Proof of Theorem~\ref{thm:random-exact-lower} and Corollary~\ref{cor:random-exact-lower}]
Use the localized family of Lemma~\ref{lem:bumps} with $2n\le M\le C_dn$, $\Delta\asymp n^{-1/d}$, and a sufficiently small $b>0$. All densities
$f_{\btheta}=f_0+H_{\btheta}$, $\btheta\in[-1,1]^M$, lie in $\F_\alpha(L)$.

Conditional on a realized design $z$, the observations are $V_i=C_{f_0}(z_i)+(T_z\btheta)_i$ with
\[
 (T_z)_{ij}=b \Delta^\alpha\int_{B_r}
       (\bu_i\cdot\bx-k_i)_+
              \psi\left(\frac{\bx-\bx_j}{\Delta}\right)\dd\bx,
\]
which is affine in $\btheta$, and $\btheta\mapsto f_{\btheta}\dd\bx$ is continuous for $W_1$. Theorem~\ref{thm:exact-information} with $\mathsf d=W_1$ and the separation \eqref{eq:bump-w1}, $\delta=cb\Delta^{\alpha+d+1}$, gives
\begin{equation}\label{eq:decoder-loss}
 \sup_{f\in\F_\alpha(L)}\E_{f,\mathsf Q_n}W_1(\widehat\mu,f\dd\bx)
   \ge\frac{cb}{4\pi^2}\Delta^{\alpha+d+1}(M-n)
   \ge c_3 n^{-\frac{\alpha+1}{d}}
\end{equation}
for every estimator, since $M-n\ge M/2\asymp\Delta^{-d}$. The same calculation replaces the basket payoff by any specified continuous linear payoff, proving the final assertion.

For the upper bound, use the minimum-norm interpolant in \eqref{eq:minnorm-estimator}. Write $\beta=(\alpha+1)/d$. For $n\ge2$, apply \eqref{eq:smooth-exact-random} with $\delta=n^{-\beta-2}$ and use $W_1\le2r$ on its complement. Uniformly over the density class,
\[
 \E_fW_1(\widehat\mu,f\dd\bx)
   \le CL\left(\frac{\log(e n^{\beta+3})}{n}\right)^\beta
                  +2r n^{-\beta-2}
   \le C'\left(\frac{\log(en)}n\right)^\beta.
\]
\end{proof}

\par\medskip\noindent\textbf{Adaptive and dependent quote arrivals.}
The coverage argument of Lemma~\ref{lem:coverage} uses independence only through the probability that a fixed ball is missed. Suppose that the conditional law of each new quote location, given the past quotes and any selection randomization, dominates $q_iQ_0$ as a measure for deterministic $q_i\in[0,1]$. Iterating the conditional miss probability for each ball of a $t$-net gives $\Pr(h_n>2t)\le Ct^{-d}\exp(-ct^d\Lambda_n)$ with $\Lambda_n=\sum_{i\le n}q_i$, and the deterministic bound \eqref{eq:smooth-exact-diameter} then yields the analogue of \eqref{eq:smooth-exact-random} with $n$ replaced by $1+\Lambda_n$ and $\log(en/\delta)$ replaced by $\log(en^3/\delta)$, simultaneously over $n$.

\section{Proofs for Section~\ref{sec:noise-uniform}}

The torus $\mathbb T_T^d$, the cutoff $\chi$, the tests $\psi_{\bk}$, the features $\widetilde G_{\bk}$ and the constants $t$, $D_\alpha$, $\vartheta$ are as in Section~\ref{sec:noise-uniform}.

\begin{proof}[Proof of Proposition~\ref{prop:noisy-upper}]
For the projection in \eqref{eq:kr-projection}, fix a countable $W_1$-dense subset of $\Pp(B_r)$ and choose its first element whose objective is within $1/n$ of the infimum over that subset. The objectives are measurable functions of the finitely many estimated moments, and continuity in $\nu$ makes this infimum equal to the infimum in \eqref{eq:kr-projection}. The construction works for any positive tolerance, including $e^{-n}$ in \eqref{eq:nu-projection}.

Let $f_M=\E_\mu\widetilde f_M$, a signed filtered version of $\mu$. For a mass-zero trigonometric density $g$ on the torus, normalization of a Lipschitz test at a point bounds its $H^1$ norm by a constant depending only on the torus. Sobolev duality therefore gives
\begin{equation}\label{eq:torus-dual}
 \norm{g\dd\bx}_{\Lip^*}\le C\norm{g}_{H^{-1}(\mathbb T_T^d)}.
\end{equation}
Orthogonality of the Fourier basis and \eqref{eq:moment-variance} give
\begin{align}
 \E_\mu\norm{\widetilde f_M-f_M}_{H^{-1}(\mathbb T_T^d)}^2
 &\le\frac Cn\sum_{\abs{\bk}\le2M}
                 (1+\abs{\bk}^2)^{-1}(1+\abs{\bk})^{2s}\,\le\frac Cn M^{2s+d-2}
  =\frac Cn M^{2d+1}.
 \label{eq:stochastic-bound}
\end{align}

If $\mu=f\dd\bx$ with $f\in\F_\alpha(L)$, its zero extension, regarded as a periodic function, has $H^\alpha(\mathbb T_T^d)$ norm at most $CL$. The multiplier $1-w(\bk/M)$ vanishes for $\abs{\bk}\le M$, and hence
\begin{equation}\label{eq:smooth-bias}
 \norm{f-f_M}_{H^{-1}(\mathbb T_T^d)}^2
 \le C M^{-2\alpha-2}\norm{f}_{H^\alpha(\mathbb T_T^d)}^2
 \le C L^2 M^{-2\alpha-2}.
\end{equation}

For an arbitrary probability measure, let $\Phi_M$ be the periodic convolution kernel with Fourier multiplier $w(\bk/M)$. Its integral is one. Poisson summation expresses it as the periodization of a rescaled Schwartz function, and therefore
\begin{equation}\label{eq:kernel-firstmoment}
 \int_{\mathbb T_T^d}\dist(\bx,0)\abs{\Phi_M(\bx)}\dd\bx\le C M^{-1}.
\end{equation}
Indeed, the corresponding integral on $\R^d$ is $M^{-1}$ times the finite first absolute moment of that Schwartz function; folding to the torus only decreases distance. Since $f_M\dd\bx=\Phi_M*\mu$, each 1-Lipschitz $F$ satisfies
\begin{align*}
 \abs{\int F\dd\mu-\int F f_M\dd\bx}
 &=\abs{\int\!\int \Phi_M(\bz)
          \left(F(\bx)-F(\bx+\bz)\right)\dd\bz\,\mu(\mathrm d\bx)}\\
 &\le CM^{-1}.
\end{align*}

The true $\mu$ is feasible for the projection in \eqref{eq:kr-projection}. The triangle inequality yields
\[
 W_1(\widehat\mu_M,\mu)
 \le2\norm{\widetilde f_M\dd\bx-\mu}_{\Lip^*}+\frac1n.
\]
Combine this inequality with \eqref{eq:torus-dual}--\eqref{eq:kernel-firstmoment} and Cauchy--Schwarz to obtain both conclusions.
\end{proof}

\begin{proof}[Proof of Proposition~\ref{prop:noisy-lower}]
Use the sign family $f_{\btheta}=f_0+H_{\btheta}$ from Lemma~\ref{lem:bumps}, with a sufficiently small $b>0$ such that every density lies in $\F_\alpha(L)$. The full observation law under $f_{\btheta}$ is denoted by $P_{\btheta}^{(n)}$. If $\btheta$ and $\btheta'$ differ only in coordinate $j$, the density difference is
\[
 2b \Delta^\alpha\theta_j\psi\left(\frac{\bx-\bx_j}{\Delta}\right).
\]
It has zero mass and first moments. Scaling its Fourier transform in \eqref{eq:fourier-identity} gives
\begin{equation}\label{eq:adjacent-pricing}
 \norm{C_{f_{\btheta}}-C_{f_{\btheta'}}}_{L^2(Q_0)}^2
     =c_\psi b^2 \Delta^{2\alpha+2d+3},
\end{equation}
where $0<c_\psi<\infty$. To verify the exponent, the transform contributes $\Delta^{2\alpha+2d}$, the frequency substitution contributes $\Delta^{-d}$, and the weight $\abs{\bxi}^{-d-3}$ contributes $\Delta^{d+3}$. The integral is finite at zero by the two vanishing moments of $\psi$ and at infinity by smoothness.

Lemma~\ref{lem:hellinger} gives
\begin{equation}\label{eq:adjacent-kl}
 H^2\left(P_{\btheta}^{(n)},P_{\btheta'}^{(n)}\right)
 \le\frac{nI_\xi}{4\varsigma^2}
       \norm{C_{f_{\btheta}}-C_{f_{\btheta'}}}_{L^2(Q_0)}^2
 =\frac{nI_\xi c_\psi b^2}{4\varsigma^2}
                                      \Delta^{2\alpha+2d+3}.
\end{equation}
Choose $\Delta\asymp n^{-1/(2\alpha+2d+3)}$ with a sufficiently small proportionality constant, so that \eqref{eq:adjacent-kl} is at most $1/16$. A grid with $M\asymp \Delta^{-d}$ disjoint cells is available, and the total variation distance between adjacent observation laws is at most $1/4$.

For an arbitrary estimator $\widehat\mu$, choose a closest member $f_{\widehat\btheta}$ of the finite family in $W_1$, resolving ties deterministically. Since
$W_1(\widehat\mu,f_{\widehat\btheta}\dd\bx)
\le W_1(\widehat\mu,f_{\btheta}\dd\bx)$,
\eqref{eq:bump-w1}, with $\norm{\btheta-\widehat\btheta}_1=2\sum_j\mathbf1_{\{\widehat\theta_j\ne\theta_j\}}$, and the triangle inequality imply
\[
 W_1(\widehat\mu,f_{\btheta}\dd\bx)
 \ge\frac c2 b \Delta^{\alpha+d+1}
             \sum_{j=1}^M\mathbf1_{\{\widehat\theta_j\ne\theta_j\}}.
\]
Average over the $2^M$ sign vectors. For each coordinate $j$, pair every vector with the vector obtained by flipping its $j$th sign. For any binary decision, the sum of its two error probabilities is at least one minus the total variation distance of the corresponding laws. The adjacent bound therefore gives
\[
 2^{-M}\sum_{\btheta}
       P_{\btheta}^{(n)}(\widehat\theta_j\ne\theta_j)\ge\frac38.
\]
It follows that
\[
 \sup_{\btheta}\E_{\btheta} W_1(\widehat\mu,f_{\btheta}\dd\bx)
 \ge cb\Delta^{\alpha+d+1}M
 \ge c'\Delta^{\alpha+1}
 \ge c''n^{-(\alpha+1)/(2\alpha+2d+3)}.
\]
Adjusting the constants covers the finitely many small $n$.
\end{proof}

The Fourier estimator of Section~\ref{sec:noise} uses the known uniform design. For a known density $q=\mathrm dQ/\mathrm dQ_0$ bounded below by a positive constant, replace $\widetilde G_\psi(Z_i)Y_i$ by $\widetilde G_\psi(Z_i)Y_i/q(Z_i)$. The upper proof is unchanged up to constants. Matching lower bounds under nonuniform designs follow from the same testing construction if $q$ is also bounded from above.

\subsection{Pointwise Control of Inverse Features}
We control the weighted sum of squared inverse features uniformly over quote locations. The bound transfers integrated pricing-surface error directly to the negative-Sobolev variance of the reconstructed law.

Use the torus, real cutoff $\chi$, and Fourier tests $\psi_{\bk}$ defined above.

\begin{lemma}[Spectral envelope]\label{lem:nu-envelope}
For every $M\ge1$, there is a jointly measurable version of the inverse features such that
\begin{equation}\label{eq:nu-envelope}
 \sup_{z\in\Z}
 \sum_{0<\abs{\bk}\le2M}(1+\abs{\bk}^2)^{-1}
                   \abs{\widetilde G_{\bk}(z)}^2
             \le C M^{2d+1}.
\end{equation}
The constant depends only on the fixed support, torus, and smooth cutoffs.
\end{lemma}

\begin{proof}
Let $\boldsymbol\omega_{\bk}=2\pi\bk/T$. The uncorrected feature from \eqref{eq:dual-feature-fourier} satisfies
\[
 \widehat G_{\bk}(\bu,t)
       =-a_{d,R}\abs{t}^{d+1}\widehat\chi(t\bu+\boldsymbol\omega_{\bk}).
\]
Fourier inversion defines a continuous function of $(\bu,k)$, since the integrand is integrable. Write $b=\bu\cdot\boldsymbol\omega_{\bk}$ and $a=\dist(\boldsymbol\omega_{\bk},\R\bu)$. For every sufficiently large $P$, Schwartz decay gives
\begin{align*}
 \sup_k\abs{G_{\bk}(\bu,k)}
 &\le C_P\int_\R\abs{t}^{d+1}
              \left(1+\sqrt{(t+b)^2+a^2}\right)^{-P}\dd t\\
 &\le C_P'(1+\abs{\bk})^{d+1}
              \int_\R(1+\abs{v})^{d+1}(1+\abs{v}+a)^{-P}\dd v.
\end{align*}
In the second line $v=t+b$ and $\abs{b}\le C\abs{\bk}$. Choosing $P$ sufficiently large for each $A$ and integrating after $v=(1+a)w$ gives
\begin{equation}\label{eq:nu-feature-tube}
 \sup_k\abs{G_{\bk}(\bu,k)}
 \le C_A(1+\abs{\bk})^{d+1}
              \left(1+\dist(\bk,\R\bu)\right)^{-A}.
\end{equation}

Uniformly in $\bu$, the number of integer lattice points with $\abs{\bk}\le2M$ and distance at most $j+1$ from $\R\bu$ is at most
\[
 C(M+j+2)(j+2)^{d-1}\le C'M(j+2)^d.
\]
To verify this, place disjoint unit cubes about these lattice points. Their union lies in a tube of length at most $CM+C(j+1)$ and radius at most $C(j+1)$; its volume gives the bound. Summing over integer distance shells, for $2A>d+1$,
\[
 \sum_{\abs{\bk}\le2M}
       \left(1+\dist(\bk,\R\bu)\right)^{-2A}\le CM.
\]
Equation~\eqref{eq:nu-feature-tube} consequently gives
\begin{equation}\label{eq:nu-main-envelope}
 \sup_{z\in\Z}\sum_{0<\abs{\bk}\le2M}
       (1+\abs{\bk}^2)^{-1}\abs{G_{\bk}(z)}^2
       \le CM^{2d}\,M=CM^{2d+1}.
\end{equation}
For the correction in \eqref{eq:tilde-feature}, the proof of Lemma~\ref{lem:dual-features} gives
$\abs{\boldsymbol b_{\psi_{\bk}}}\le C\norm{\psi_{\bk}}_{H^s}
\le C(1+\abs{\bk})^s$. Hence
\[
 \sup_{\bu\in S^{d-1}}\sum_{0<\abs{\bk}\le2M}
      (1+\abs{\bk}^2)^{-1}
       \abs{2d\boldsymbol b_{\psi_{\bk}}\cdot\bu}^2
       \le CM^{2s+d-2}=CM^{2d+1}.
\]
Combine this with \eqref{eq:nu-main-envelope} and $\abs{x+y}^2\le2\abs{x}^2+2\abs{y}^2$.
\end{proof}

\begin{proof}[Proof of Lemma~\ref{lem:nu-pilot}]
We describe one coordinate cube $D=[0,1]^d$. The sampling law dominates a positive multiple of Lebesgue measure on this cube; samples outside the cube are ignored. The same construction will be applied in finitely many sphere--strike charts.

Fix a coordinatewise polynomial degree $p$. Choose $J=(p+1)^d$ tensor-product interpolation nodes strictly inside $D$, with mutually disjoint small neighborhoods $A_1,\ldots,A_J$. Write $b_1,\ldots,b_J$ for a basis of the coordinatewise polynomial space. The matrix $B=(b_j(a_l))_{l,j}$ at the reference nodes is invertible. By continuity, the neighborhoods can be chosen so small that every matrix whose $l$th row is an average of $(b_j(x))_j$ over points in $A_l$ has inverse norm at most $2\norm{B^{-1}}$.
For a dyadic cube $D'$ of side $a$, transport these neighborhoods and basis functions by its affine coordinate map. Call $D'$ good if every transported neighborhood contains at least one of the $m$ samples. Each neighborhood has sampling probability at least $c a^d$. Thus
\begin{equation}\label{eq:nu-bad-cube}
 \Pr(D'\text{ is not good})\le J\exp(-cma^d).
\end{equation}
On a good cube, fit a polynomial by equating its average over the actual sample points in each neighborhood to the corresponding average observed price. This fit reproduces every polynomial in the chosen space exactly. Moreover, on the rescaled unit cube its absolute value is bounded by a constant times the largest absolute data value. Its conditional variance at each point is at most
\begin{equation}\label{eq:nu-pilot-variance}
 C\varsigma^2\sum_{l=1}^J\frac{1}{\mathsf N_l}\le CJ\varsigma^2,
\end{equation}
where $\mathsf N_l\ge1$ is the count in neighborhood $l$. Distinct neighborhoods in this fit use disjoint samples.

Choose a dyadic base side $h\ge m^{-1/d}$. A prescribed bandwidth in $[m^{-1/d},1]$ is rounded upward to a dyadic side, by a factor smaller than two. We retain the notation $h$ for this side. For each base cube, use the fit on its smallest good dyadic ancestor, including itself. If no ancestor is good, use zero on that base cube. The resulting piecewise polynomial operator, denoted by $P_m$, is a measurable function of the design and is linear in the data.

Let $G\in H^t(D)$, $t>d/2$, and choose $p\ge\lceil t\rceil$. On a cube $D'$ of side $a$, polynomial approximation and Sobolev embedding give a polynomial $P$ in the fitted space with
\begin{equation}\label{eq:nu-local-approx}
 \norm{G-P}_{L^\infty(D')}
       \le C a^{t-d/2}|G|_{H^t(D')}.
\end{equation}
Here $|\cdot|_{H^t}$ is the top-order derivative seminorm for integer $t$, and the fractional seminorm of the derivatives of order $\lfloor t\rfloor$ otherwise. To obtain \eqref{eq:nu-local-approx}, on the unit cube subtract a bounded polynomial projection of degree $t-1$ for integer $t$, or degree $\lfloor t\rfloor$ for noninteger $t$. The generalized Poincar\'e inequality bounds the full $H^t$ norm of this remainder by its top-order seminorm. The embedding $H^t(D)\subset C(D)$ bounds the supremum norm. Affine rescaling then supplies $a^{t-d/2}$. 

Stability and polynomial reproduction give the same bound, up to a fixed constant, for the error of the fit to exact values of $G$ on a good $D'$. If a base cube uses an ancestor of side $a=2^\ell h$ with $\ell\ge1$, the preceding ancestor of side $a/2$ was not good. Its selection probability is therefore at most
\[
 \pi_\ell=\min\left(1,J\exp\left[-cm(2^{\ell-1}h)^d\right]\right),
 \qquad \pi_0=1.
\]
At each level, each cube $D'$ contains at most $a^d/h^d$ base cubes. Applying \eqref{eq:nu-local-approx} on their possible parent, then summing the squared integrated errors, gives
\begin{align}
 \E\norm{P_mG-G}_{L^2(D)}^2
 &\le C\sum_{\ell\ge0}\pi_\ell(2^\ell h)^{2t}
             \sum_{D'\text{ at level }\ell}|G|_{H^t(D')}^2
       +C e^{-cm}\norm{G}_{L^2(D)}^2\notag\\
 &\le C h^{2t}\norm{G}_{H^t(D)}^2.
 \label{eq:nu-pilot-sobolev}
\end{align}
The sums at each level are over disjoint cubes and are bounded by the global seminorm, also at fractional orders. The numerical series is finite uniformly in $m$, since $mh^d\ge1$ and exponential decay in $2^{\ell d}$ dominates $2^{2t\ell}$. If the root cube is not good, \eqref{eq:nu-bad-cube} bounds its contribution by the last term. This term is also at most $C h^{2t}\norm{G}_{H^t}^2$. For the sharper noise bound, let $\mathsf N$ be the count in one neighborhood of a cube of side $a$, with sampling probability $p_a\ge ca^d$. Since $\mathsf N\sim\operatorname{Bin}(m,p_a)$,
\begin{equation}\label{eq:nu-inverse-count}
 \E\left[\frac{\mathbf1_{\{\mathsf N\ge1\}}}{\mathsf N}\right]
 \le2\E\left[\frac1{\mathsf N+1}\right]
 =\frac{2[1-(1-p_a)^{m+1}]}{(m+1)p_a}
 \le\frac{C}{ma^d}.
\end{equation}
For a base cube, drop its ancestor-selection indicator when bounding each nonnegative noise-variance contribution in \eqref{eq:nu-pilot-variance}. Applying \eqref{eq:nu-inverse-count} at every possible ancestor and then integrating over the base cubes gives
\begin{equation}\label{eq:nu-averaged-noise}
 \E\norm{P_m(\varsigma\boldsymbol\xi)}_{L^2(D)}^2
 \le C\varsigma^2\sum_{\ell\ge0}\frac1{m(2^\ell h)^d}
 \le\frac{C'\varsigma^2}{mh^d}.
\end{equation}
Conditional on the design, $P_m(\varsigma\boldsymbol\xi)$ has mean zero. Thus its variance adds to the expected squared bias without a cross term.

For the sphere, take a finite collection of larger smooth coordinate cubes and smaller cubes whose images cover $S^{d-1}$, with a smooth partition of unity supported in the larger charts. Include the full strike interval by an affine coordinate. Chart Jacobians and their reciprocals are bounded, and $Q_0$ has a coordinate density bounded below. Construct $P_m$ in each larger chart and combine the outputs using the partition of unity. Coordinate Sobolev norms are bounded by the corresponding global norm. Squared errors and conditional variances of the finite sum are bounded by a constant multiple of the sum of the individual bounds. The strike endpoints are faces of the coordinate cubes and require no samples outside the observed interval. Hence the analogue of \eqref{eq:nu-pilot-sobolev} holds on $\Z$.

For $f\in\F_\alpha(L)$, subtract a smooth probability density $f_*$ supported inside $B_r$. Lemma~\ref{lem:sobolev-map} applies to $f-f_*$, and $C_{f_*}$ is smooth on $\Z$. Consequently
\[
 \norm{C_f}_{H^t(\Z)}\le C(1+L)\le C'L.
\]
The last inequality uses the positive lower bound on the $L^2$ norm of a probability density on $B_r$. Since $C_f$ takes values in $[0,2R]$, truncating the output pointwise to this interval decreases its squared error. Taking $h\asymp m^{-1/d}$ in the first inequality of \eqref{eq:nu-bandwidth-pilot} gives \eqref{eq:nu-smooth-pilot}.

For \eqref{eq:nu-universal-pilot}, fix $\bx\in B_r$ and write $G_{\bx}(\bu,k)=(\bu\cdot\bx-k)_+$. In each coordinate chart, its nonsmooth set is the graph $k=\bu\cdot\bx$. The volume of its $a$-neighborhood is at most $Ca$, uniformly in $\bx$: the defining function $\bu\cdot\bx-k$ has a nonzero strike derivative and uniformly bounded angular derivatives. Cubes of side $a$ meeting this graph have total volume at most $Ca$. On such cubes, the Lipschitz property and reproduction of constants give a fit error at most $Ca$. On the remaining cubes, $G_{\bx}$ is either zero or a smooth function with uniformly bounded second derivatives; reproduction of affine polynomials gives an error at most $Ca^2$. At each level, the squared integrated error is therefore bounded by
\[
 C\left(a\,a^2+a^4\right)\le Ca^3.
\]
Multiplication by the same selection probabilities $\pi_\ell$ and summation over levels yield
\begin{equation}\label{eq:nu-atom-pilot}
 \sup_{\bx\in B_r}\E\norm{P_mG_{\bx}-G_{\bx}}_{L^2(Q_0)}^2
                  \le C h^3.
\end{equation}
Before truncation, linearity gives $P_mC_\mu=\int P_mG_{\bx}\,\mu(\mathrm d\bx)$ for exact inputs. Jensen's inequality in $L^2(Q_0)$ and Tonelli's theorem imply
\[
 \E\norm{P_mC_\mu-C_\mu}_{L^2(Q_0)}^2
 \le\int\E\norm{P_mG_{\bx}-G_{\bx}}_{L^2(Q_0)}^2\,\mu(\mathrm d\bx)
 \le Ch^3.
\]
The noise contribution is bounded by \eqref{eq:nu-averaged-noise}, and truncation again decreases the error. Taking $h\asymp m^{-1/d}$ in the second inequality of \eqref{eq:nu-bandwidth-pilot} gives \eqref{eq:nu-universal-pilot}. For $Q\ge q_*Q_0$, every neighborhood in \eqref{eq:nu-bad-cube} still has probability at least $cq_*a^d$.
\end{proof}\begin{proof}[Proof of Proposition~\ref{prop:nu-transfer}]
Let $\mathcal D_0$ denote the first-group observations. Conditional on $\mathcal D_0$, Lemma~\ref{lem:dual-features} and independence give
\begin{equation}\label{eq:nu-cv-unbiased}
 \E\left[\widehat\theta_{\bk}^{\mathrm{cv}}\mid\mathcal D_0\right]
   =a_{\bk}+\int\widetilde G_{\bk}(C_\mu+e)\dd Q_0
   =\theta_{\bk}(\mu)+b_{\bk},
 \qquad b_{\bk}=\int\widetilde G_{\bk}\,e\dd Q_0.
\end{equation}
Write $R_0=C_\mu+e-\widehat C$, so that $\E\norm{R_0}_{L^2(Q_0)}^2\le2v_n(\mu)+2\epsilon^2$. The conditional variance obeys
\begin{equation}\label{eq:nu-cv-var}
 \E\left[\abs{\widehat\theta_{\bk}^{\mathrm{cv}}-\theta_{\bk}(\mu)-b_{\bk}}^2
                 \mid\mathcal D_0\right]
 \le\frac1{n_1}\int\abs{\widetilde G_{\bk}(z)}^2
                     \left[R_0(z)^2+\varsigma^2\right]\dd Q_0(z).
\end{equation}
Apply \eqref{eq:nu-envelope} inside the integral, and then take expectation over $\mathcal D_0$. Let $f_M$ be the deterministic filtered version of $\mu$ and
\[
 B(\bx)=\frac1{T^d}\sum_{\bk\ne0}w(\bk/M)\,b_{\bk}\,e^{2\pi\mathrm i\bk\cdot\bx/T}.
\]
This gives
\begin{equation}\label{eq:nu-cv-hminusone}
 \E\norm{\widetilde f_M^{\mathrm{cv}}-f_M-B}_{H^{-1}(\mathbb T_T^d)}^2
       \le\frac{CM^{2d+1}}{n}
                      \left[\varsigma^2+\epsilon^2+v_n(\mu)\right].
\end{equation}

Since $B$ has zero mass, we may restrict to mean-zero $1$-Lipschitz functions $F$ on the torus. Then $\widehat F(\boldsymbol0)=0$ and $\norm F_{H^1(\mathbb T_T^d)}\le C$. Since $\psi\mapsto\widetilde G_\psi$ is linear,
\[
 \int FB\dd\bx=\int_\Z\widetilde G_\Psi\,e\dd Q_0,\qquad\Psi=\chi\,(\Phi_M*F),
\]
with $\Phi_M$ the convolution kernel in \eqref{eq:kernel-firstmoment}. Lemma~\ref{lem:dual-features}, $\abs e\le\epsilon$, and $\norm{\Phi_M*F}_{H^s(\mathbb T_T^d)}\le CM^{s-1}\norm F_{H^1(\mathbb T_T^d)}$ give
\begin{equation}\label{eq:nu-systematic-bias}
 \norm{B\dd\bx}_{\Lip^*}\le C\epsilon M^{s-1}.
\end{equation}
The smooth and arbitrary-law bias bounds are exactly \eqref{eq:smooth-bias} and \eqref{eq:kernel-firstmoment}. The true law is feasible in \eqref{eq:nu-projection}, so
\[
 W_1(\widehat\mu_M^{\mathrm{cv}},\mu)
 \le2\norm{\widetilde f_M^{\mathrm{cv}}\dd\bx-\mu}_{\Lip^*}+e^{-n}.
\]
Use \eqref{eq:torus-dual}, \eqref{eq:nu-cv-hminusone}, \eqref{eq:nu-systematic-bias}, and Cauchy--Schwarz to obtain both inequalities.
\end{proof}

\begin{lemma}[Curved wave packets]\label{lem:curved-packets}
Let $\alpha\ge0$, let $f_0$ be the base density of Lemma~\ref{lem:bumps}, and let $L>L_0$ as in Theorem~\ref{thm:exact-lower}. There are $c,C>0$ and $N_0$ such that for every integer $N\ge N_0$ there are $f_N^\pm\in\F_\alpha(L)$, whose one-dimensional coordinate marginals equal those of $f_0$, with
\[
 W_1(f_N^+\dd\bx,f_N^-\dd\bx)\ge cN^{-\alpha-1},\qquad
 \sup_{\bu\in S^{d-1},\,k\in\R}\abs{C_{f_N^+}(\bu,k)-C_{f_N^-}(\bu,k)}\le CN^{-\alpha-s}.
\]
\end{lemma}

\begin{proof}
Fix $\bx_1$ and $\rho_1>0$ such that $f_0$ is bounded below on $B(\bx_1,\rho_1)$, let $\boldsymbol a=(\be_1+\be_2)/\sqrt2$, and let $\bx_c=\bx_1-8\rho_1\boldsymbol a$. Let $\chi\in C_c^\infty(B(\bx_1,\rho_1))$ be nonnegative and nonzero, $\varphi(\bx)=\abs{\bx-\bx_c}$, and
\[
 \eta_N=bN^{-\alpha-2}\partial_1\partial_2\left[\chi\cos(N\varphi)\right],\qquad f_N^\pm=f_0\pm\eta_N.
\]
Since the function $\chi\cos(N\varphi)$ has fixed support and derivatives of order $m$ bounded by $C_mN^m$, $\norm{\eta_N}_{H^m}\le C_mbN^{m-\alpha}$ for integers $m$, interpolation gives $\norm{\eta_N}_{H^\alpha}\le Cb$, and $\norm{\eta_N}_\infty\le CbN^{-\alpha}$. Since $\eta_N$ is a derivative $\partial_1\partial_2$ of a compactly supported function, it has zero mass and first moments, and integrating out $x_1$ or $x_2$ shows that its one-dimensional coordinate marginals vanish. Choosing $b$ so small that $f_0\pm\eta_N\ge0$ and $L_0+Cb\le L$, where $L_0=\norm{f_0}_{H^\alpha}$, gives $f_N^\pm\in\F_\alpha(L)$.

\emph{Separation.} Let $\tilde\chi\in C_c^\infty(B(\bx_1,2\rho_1))$ equal one on $B(\bx_1,\rho_1)$ with $0\le\tilde\chi\le1$, and $F=(2N)^{-1}\tilde\chi\cos(N\varphi)$. Since $\abs{\nabla\varphi}=1$, $F$ is $1$-Lipschitz for $N\ge\norm{\nabla\tilde\chi}_\infty$. On $\supp\chi$, the unit vector $\nabla\varphi$ makes an angle at most $\arcsin\tfrac18$ with $\boldsymbol a$, hence $\partial_1\varphi\,\partial_2\varphi\ge\tfrac12\cos(2\arcsin\tfrac18)\ge\tfrac14$, and $\partial_1\partial_2F=-\tfrac N2\cos(N\varphi)\,\partial_1\varphi\,\partial_2\varphi-\tfrac12\sin(N\varphi)\,\partial_1\partial_2\varphi$. Integrating by parts twice,
\[
 \begin{aligned}
 \int F\eta_N\dd\bx&=bN^{-\alpha-2}\int\chi\cos(N\varphi)\,\partial_1\partial_2F\dd\bx\\
 &=-\frac b4N^{-\alpha-1}\int\chi\,\partial_1\varphi\,\partial_2\varphi\,\big(1+\cos(2N\varphi)\big)\dd\bx\\
 &\qquad-\frac b4N^{-\alpha-2}\int\chi\,\partial_1\partial_2\varphi\,\sin(2N\varphi)\dd\bx.
 \end{aligned}
\]
One integration by parts along $\nabla\varphi$ bounds both oscillatory integrals by $C/N$. Hence $\abs{\int F\eta_N\dd\bx}\ge\tfrac b{32}N^{-\alpha-1}\int\chi$ for large $N$, and since $\pm F$ are $1$-Lipschitz, $W_1(f_N^+\dd\bx,f_N^-\dd\bx)\ge2\abs{\int F\eta_N\dd\bx}\ge\tfrac b{16}N^{-\alpha-1}\int\chi$.

\emph{Quotes.} For $\bu\in S^{d-1}$ and $k\in\R$, the distributional derivative $\partial_1\partial_2$ of $\bx\mapsto(\bu\cdot\bx-k)_+$ is $u_1u_2$ times surface measure on the hyperplane $H_{\bu,k}=\{\bu\cdot\bx=k\}$. Hence
\[
 C_{f_N^+}(\bu,k)-C_{f_N^-}(\bu,k)=2bN^{-\alpha-2}u_1u_2\int_{H_{\bu,k}}\chi\cos(N\varphi)\dd S,
\]
and $\abs{u_1u_2}\le1$.
Let $\delta=\dist(\bx_c,H_{\bu,k})$, and use coordinates $\by\in\R^{d-1}$ on $H_{\bu,k}$ centered at the foot of the perpendicular from $\bx_c$, so that $\varphi=(\delta^2+\abs\by^2)^{1/2}$. On $\supp\chi$, $\abs{\bx-\bx_c}\in[7\rho_1,9\rho_1]$. If $\delta\le4\rho_1$, the tangential gradient $\by/\varphi$ of $\varphi$ has norm at least $(1-16/49)^{1/2}$ on $H_{\bu,k}\cap\supp\chi$, and repeated integration by parts gives $O(N^{-m})$ for every $m$. If $\delta>4\rho_1$, the phase has the single critical point $\by=0$, with Hessian $\delta^{-1}I_{d-1}$ and $\abs{\nabla\varphi}\ge\abs\by/(9\rho_1)$ on the support, and the phase and amplitude have derivatives bounded uniformly in $(\bu,k)$. The stationary phase estimate \cite[Theorem~7.7.5]{hormander1983} gives $\abs{\int_{H_{\bu,k}}\chi\cos(N\varphi)\dd S}\le CN^{-(d-1)/2}$ uniformly in $(\bu,k)$. Thus the quote difference is at most $CbN^{-\alpha-2-(d-1)/2}=CbN^{-\alpha-s}$.
\end{proof}

The pair in Lemma~\ref{lem:curved-packets} lies in the fiber $\K^\epsilon_{n,\alpha,L}(z,v)$ of every design, with $v$ the average of the two quote vectors and $\epsilon=CN^{-\alpha-s}$. Hence, for each $L>L_0$, the exponent $\vartheta$ of $\epsilon$ in \eqref{eq:bounded-errors} is optimal.

\begin{proof}[Proof of Theorem~\ref{thm:nu-smooth}]
Before clipping to $[0,2R]$, the pilot of Lemma~\ref{lem:nu-pilot} is linear in the observed prices, and on each cube its fit is bounded by a constant times the largest absolute data value in that cube. The error surface therefore changes the unclipped pilot by at most $C\epsilon$ uniformly on $\Z$, and clipping does not increase the distance to $C_f$. Hence \eqref{eq:nu-smooth-pilot} holds with the additional term $C\epsilon^2$. For $n\ge2$, Proposition~\ref{prop:nu-transfer} gives
\begin{equation}\label{eq:nu-prebalance}
 \begin{aligned}
 \sup_{f\in\F_\alpha(L)}\ \sup_{\abs e\le\epsilon}\E W_1(\widehat\mu_M^{\mathrm{cv}},f\dd\bx)
 &\le C\Big[L M^{-\alpha-1}+\epsilon M^{s-1}\\
 &\qquad+\big(\varsigma n^{-\frac{1}{2}}+\epsilon n^{-\frac{1}{2}}+L n^{-\frac{1}{2}-\frac{t}{d}}\big)M^{d+\frac{1}{2}}\Big]+e^{-n}.
 \end{aligned}
\end{equation}
Take
\begin{equation}\label{eq:nu-cutoff}
 M=\max\left(1,\left\lfloor
       \min\left(n^{\frac{1}{d}},\left(\frac n{\varsigma^2}\right)^{\frac{1}{D_\alpha}},\left(\frac L\epsilon\right)^{\frac{1}{\alpha+s}}\right)
                           \right\rfloor\right),
\end{equation}
where the second and third entries are $+\infty$ at $\varsigma=0$ and $\epsilon=0$. If $M\le n^{1/d}$, then
\begin{align}
 n^{-\frac{1}{2}-\frac{t}{d}}M^{d+\frac{1}{2}}
 &=M^{-\alpha-1}
      n^{-D_\alpha/(2d)}M^{D_\alpha/2}\le M^{-\alpha-1},
 \label{eq:nu-design-balance}
\end{align}
and $\epsilon n^{-1/2}M^{d+1/2}\le\epsilon M^{s-1}$, since $M^{d+1/2-(s-1)}=M^{d/2}\le n^{1/2}$.
Likewise, $M\le(n/\varsigma^2)^{1/D_\alpha}$ implies
$\varsigma n^{-1/2}M^{d+1/2}\le M^{-\alpha-1}$, and $M\le(L/\epsilon)^{1/(\alpha+s)}$ implies $\epsilon M^{s-1}\le LM^{-\alpha-1}$.
Thus \eqref{eq:nu-prebalance} is at most a constant times
\[
 n^{-\frac{\alpha+1}{d}}
   +\left(\frac{\varsigma^2}{n}\right)^{\frac{\alpha+1}{D_\alpha}}+\epsilon^\vartheta.
\]
The exponentially small projection tolerance is absorbed for every $\alpha$. For the finitely many remaining sample sizes, and when the minimum in \eqref{eq:nu-cutoff} is below one, the estimator $\delta_0$ and the bound $W_1\le2r$ suffice after increasing the constant.

For the lower bound, the error surface $e=0$ is admissible, hence $\mathfrak R_n(\F_\alpha(L),\varsigma,\epsilon)\ge\mathfrak R_n(\F_\alpha(L),\varsigma)$. For the first lower term, use the localized family of Theorem~\ref{thm:random-exact-lower}, with $M_0\asymp n$ coefficients. The data are $C_{f_0}(Z_i)+(T_Z\btheta)_i+\varsigma\xi_i$, hence Theorem~\ref{thm:exact-information} applies with the exogenous shift $a_i=\varsigma\xi_i$ and gives
\begin{equation}\label{eq:nu-information-floor}
 \mathfrak R_n(\F_\alpha(L),\varsigma)
                    \ge c n^{-(\alpha+1)/d}
\end{equation}
uniformly in $\varsigma$, including zero.

For $\varsigma>0$, use the localized sign family of Proposition~\ref{prop:noisy-lower}. Adjacent alternatives at spatial scale $\Delta$ satisfy, by Lemma~\ref{lem:hellinger},
\[
 H^2(P_{\btheta}^{(n)},P_{\btheta'}^{(n)})
       \le\frac{nI_\xi c_\psi b^2}{4\varsigma^2}\Delta^{D_\alpha}.
\]
Choose $\Delta=c_0(\varsigma^2/n)^{1/D_\alpha}$. A fixed $c_0>0$, depending on $\overline\varsigma$ and the support geometry, ensures both that the bump construction applies for every $n$ and that this distance is at most $1/16$. There are $M_0\asymp \Delta^{-d}$ cells. The paired sign-testing argument in the proof of Proposition~\ref{prop:noisy-lower}, with \eqref{eq:bump-w1} and \eqref{eq:adjacent-kl}, gives
\begin{equation}\label{eq:nu-noise-lower}
 \mathfrak R_n(\F_\alpha(L),\varsigma)
        \ge c \Delta^{\alpha+1}
        =c'\left(\frac{\varsigma^2}{n}\right)^{\frac{\alpha+1}{D_\alpha}}.
\end{equation}

For the systematic term, let $\epsilon_1>0$ be small and $0<\epsilon\le\epsilon_1$, and take the pair of Lemma~\ref{lem:curved-packets} with $N=\lceil(C/\epsilon)^{1/(\alpha+s)}\rceil$, so that $\abs{C_{f_N^+}-C_{f_N^-}}\le\epsilon$ on $\Z$. The error surfaces $e_\pm=\mp\tfrac12(C_{f_N^+}-C_{f_N^-})$ satisfy $\abs{e_\pm}\le\epsilon$, and $C_{f_N^+}+e_+=C_{f_N^-}+e_-$. The data therefore have the same law under $(f_N^+,e_+)$ and $(f_N^-,e_-)$, and every estimator has risk at least $\tfrac12W_1(f_N^+\dd\bx,f_N^-\dd\bx)\ge cN^{-\alpha-1}\ge c'\epsilon^\vartheta$ under one of them. For $\epsilon_1<\epsilon\le\overline\epsilon$, monotonicity in $\epsilon$ gives the bound with a smaller constant. The largest of the three lower bounds is at least a third of their sum.
\end{proof}

\begin{proof}[Proof of Corollary~\ref{cor:nu-transition}]
Equate the powers in \eqref{eq:nu-smooth-rate} and use
$2\alpha+2d+3=2[\alpha+(d+3)/2]+d$ and $(\alpha+1)/\vartheta=\alpha+(d+3)/2$. Substitution of $\varsigma_n=n^{-\gamma}$ gives the second assertion.
\end{proof}

\begin{proof}[Proof of Corollary~\ref{cor:nu-general-design}]
Set $q=s-1$. There is a linear smoothing operator $J_h$ on $\Z$ such that
\begin{equation}\label{eq:nu-smoothing}
 \norm{J_hG}_{H^q}\le C h^{-q}\norm{G}_{L^2},\qquad
 \norm{J_hG-G}_{H^q}\le C h^{t-q}\norm{G}_{H^t}.
\end{equation}
To construct it, use a finite chart partition $\eta_j$, larger chart cutoffs $\chi_j$ equal to one on $\supp\eta_j$, and extension operators $E_j$ from the coordinate cubes to Euclidean space. Each $E_j$ is bounded simultaneously in $L^2$ and $H^t$. On a cube, finite linear combinations of reflections at each face, matching normal derivatives through an integer order above $t$, followed by a smooth cutoff, give such an extension; interpolation supplies the intermediate orders. In coordinates, the $j$th summand of $J_hG$ is
\[
 \chi_j\,\mathcal F^{-1}\left[w(h\bxi)\,
                      \mathcal F\left(E_j(\eta_jG)\right)(\bxi)\right],
\]
restricted to its chart and extended by zero off that chart. The sum of the unsmoothed summands is $G$. The multiplier bounds
\[
 (1+\abs{\bxi}^2)^{q/2}\abs{w(h\bxi)}\le Ch^{-q},\qquad
 (1+\abs{\bxi}^2)^{(q-t)/2}\abs{1-w(h\bxi)}\le Ch^{t-q}
\]
and boundedness of smooth multiplication and coordinate changes prove \eqref{eq:nu-smoothing}. Strike endpoints are handled by the same cube extensions; no boundary condition on $G$ is imposed.

For a dyadic bandwidth $h\in[n^{-1/d},1]$, let $\widehat C_{n,h}$ be the pilot of Lemma~\ref{lem:nu-pilot} built from all $n$ observations, and set $\widetilde C_n=J_h\widehat C_{n,h}$. Let $\mathcal S=\{C_f:f\in\F_\alpha(L)\}\subset H^q(\Z)$. The set $\mathcal S$ is convex since $\F_\alpha(L)$ is convex and $\cA$ is linear. The class $\F_\alpha(L)$ is weakly compact in $H^\alpha$. By Lemma~\ref{lem:sobolev-map} and the compact embedding of $H^t(\Z)$ into $H^q(\Z)$, its image $\mathcal S$ is compact. Let $\widehat C_n$ be the metric projection of $\widetilde C_n$ onto $\mathcal S$ in $H^q(\Z)$. This projection is unique and nonexpansive. The inverse bound \eqref{eq:inverse-sobolev} makes $C_f\mapsto f\dd\bx$ continuous from $\mathcal S$ to $W_1$. Hence the law $\widehat\mu_n=f\dd\bx$ with $C_f=\widehat C_n$ is a measurable probability-valued estimator.

As in the proof of Theorem~\ref{thm:nu-smooth}, the error surface changes the unclipped pilot by at most $C\epsilon$ uniformly on $\Z$. Hence \eqref{eq:nu-bandwidth-pilot} holds with the additional term $C\epsilon^2$, and \eqref{eq:nu-smoothing} gives
\begin{equation}\label{eq:nu-direct-risk}
 \begin{aligned}
 \E\norm{\widetilde C_n-C_f}_{H^q}^2
 &\le2\E\norm{J_h(\widehat C_{n,h}-C_f)}_{H^q}^2+2\norm{J_hC_f-C_f}_{H^q}^2\\
 &\le C\left[L^2h^{2(t-q)}+\frac{\varsigma^2}{nh^{d+2q}}+\frac{\epsilon^2}{h^{2q}}\right].
 \end{aligned}
\end{equation}
Since $C_f\in\mathcal S$ and the projection is nonexpansive, $\norm{\widehat C_n-C_f}_{H^q}\le\norm{\widetilde C_n-C_f}_{H^q}$, and \eqref{eq:inverse-sobolev} gives \eqref{eq:nu-direct-upper}.
Take $h$ of order
$\max(n^{-1/d},(\varsigma^2/n)^{1/D_\alpha},(\epsilon/L)^{1/(\alpha+s)})$,
rounded to a dyadic side and capped at one. At $h\ge(\varsigma^2/n)^{1/D_\alpha}$ the second term is at most $(\varsigma^2/n)^{(\alpha+1)/D_\alpha}$, since $\tfrac12D_\alpha-d-\tfrac12=\alpha+1$. At $h\ge(\epsilon/L)^{1/(\alpha+s)}$ the third term is at most $L^{1-\vartheta}\epsilon^\vartheta$, since $1-(s-1)/(\alpha+s)=\vartheta$. The first term at the chosen $h$ is a constant multiple of the largest of the three rates. At zero noise and $\epsilon=0$, the choice $h\asymp n^{-1/d}$ gives the upper bound in \eqref{eq:nu-general-design}.

The exact lower bound is Theorem~\ref{thm:random-exact-lower}, and the conditioning argument for \eqref{eq:nu-information-floor} preserves it under noise. If $Q\le q^*Q_0$, adjacent sign alternatives satisfy
\[
 H^2(P_{\btheta,Q}^{(n)},P_{\btheta',Q}^{(n)})
 \le\frac{nI_\xi}{4\varsigma^2}\norm{C_{f_{\btheta}}-C_{f_{\btheta'}}}_{L^2(Q)}^2
 \le\frac{nI_\xi q^*c_\psi b^2}{4\varsigma^2}\Delta^{2\alpha+2d+3}.
\]
The scale choice and testing proof of \eqref{eq:nu-noise-lower} therefore apply with constants depending on $q^*$. The argument for the systematic term in the proof of Theorem~\ref{thm:nu-smooth} applies to every design, since the two data laws coincide.
\end{proof}\begin{proof}[Proof of Corollary~\ref{cor:nu-marginals}]
By Lemma~\ref{lem:bumps}, every $f_{\btheta}=f_0+H_{\btheta}$ has the coordinate marginals of $f_0$, so the proofs of \eqref{eq:nu-information-floor} and \eqref{eq:nu-noise-lower} apply verbatim on $\mathcal P_0$. For the systematic term, the pair of Lemma~\ref{lem:curved-packets} lies in $\mathcal P_0$, and the last argument in the proof of Theorem~\ref{thm:nu-smooth} applies. The upper bounds hold since $\mathcal P_0\subseteq\F_\alpha(L)$.
\end{proof}

With exact quotes, $\varsigma=0$ and $\epsilon=0$, an estimate can also be made exactly price feasible without changing its risk order. For any estimator $\widehat\mu$, choose a nearest law $\overline\mu$ to $\widehat\mu$ in the unrestricted exact quote fiber on $B_r$. A nearest law exists by compactness. The constraint graph is closed in the design, price vector, and candidate law. Its compact sections form a measurable correspondence on the set of feasible data. The measurable minimum theorem therefore supplies a measurable choice. Since the true law belongs to this fiber,
\[
 W_1(\overline\mu,\mu)
 \le W_1(\overline\mu,\widehat\mu)+W_1(\widehat\mu,\mu)
 \le2W_1(\widehat\mu,\mu).
\]

Theorem~\ref{thm:nu-smooth} and Corollary~\ref{cor:nu-general-design} carry no logarithmic factor since the local estimate \eqref{eq:nu-pilot-sobolev} averages the bad-cell probabilities for each pricing surface, whereas \eqref{eq:smooth-exact-random} controls all feasible price vectors on a common design event, as discussed in Section~\ref{sec:exact-lower-bounds}.

\begin{proof}[Proof of Theorem~\ref{thm:nu-all-laws}]
The upper bound follows from \eqref{eq:nu-universal-pilot} and \eqref{eq:nu-transfer-all}:
\[
 \sup_\mu\E W_1(\widehat\mu_M^{\mathrm{cv}},\mu)
       \le C\left[M^{-1}
          +n^{-\frac{1}{2}}\left(\varsigma+n^{-3/(2d)}\right)M^{d+\frac{1}{2}}\right]+e^{-n}.
\]
Choose $M$, up to integer rounding and a lower truncation at one, of order
\[
 \left(  \frac{n}{\varsigma^2+n^{-3/d}} \right)^{\frac{1}{2d+3}}.
\]
The two displayed terms are then bounded by
$C[(\varsigma^2+n^{-3/d})/n]^{1/(2d+3)}$.
Small $n$ are covered as in Theorem~\ref{thm:nu-smooth}. Apply that theorem with $\alpha=0$ and one sufficiently large $L$ to obtain the lower side, since this density class is contained in $\Pp(B_r)$. If $\varsigma^2\ge c_0n^{-3/d}$, the noise term bounds the exact upper term up to a fixed factor, while it is itself a lower bound.
\end{proof}

\section{Proofs for Section~\ref{sec:stable-pricing}}
\begin{proof}[Proof of Proposition~\ref{prop:price-diameter}]
The left side of \eqref{eq:price-diameter} equals
\[
 \sup_{\Lip(g)\le1}\sup_{\nu,\nu'\in\K}\int g\dd(\nu-\nu')
 =\sup_{\nu,\nu'\in\K}\sup_{\Lip(g)\le1}\int g\dd(\nu-\nu')
 =\sup_{\nu,\nu'\in\K}W_1(\nu,\nu').
\]
For the convergence of endpoints, the functional $\nu\mapsto\int g\dd\nu$ is continuous on the compact space $\Pp(B_R)$ and hence uniformly continuous. The Hausdorff convergence in Theorem~\ref{thm:set-consistency} therefore gives convergence of its maximum and minimum over the admissible sets.
\end{proof}

\begin{proof}[Proof of Lemma~\ref{lem:replication-bound}]
Let $\mu,\nu$ be admissible in the definition of $\varpi_g(\epsilon)$ and write $\eta=\mu-\nu$. Then $\eta(B_R)=0$ and $\norm\eta\le2$. For $h=g-\cA^*m$ and $c$ the midpoint of the range of $h$, $\int h\dd\eta=\int(h-c)\dd\eta\le\tfrac12\osc(h)\norm\eta\le\osc(h)$. Fubini's theorem gives $\int\cA^*m\dd\eta=\int C_\eta\dd m\le\epsilon\norm m$, since $\abs{C_\eta}\le\epsilon$ on $\Z_E$. Adding the two bounds and taking the supremum over admissible pairs proves \eqref{eq:replication-bound}.
\end{proof}

\begin{lemma}[Continuity under identification]\label{lem:modulus-continuity}
For every $g\in C(B_R)$, the modulus $\varpi_g$ is nondecreasing. If the pricing surfaces on $\Z_E$ identify $\Pp(B_R)$, then $\varpi_g(\epsilon)\to0$ as $\epsilon\downarrow0$.
\end{lemma}

\begin{proof}
Monotonicity is immediate from the definition. Suppose that identification holds and that there are $\delta>0$, $\epsilon_j\downarrow0$, and $\mu_j,\nu_j\in\Pp(B_R)$ with $\sup_{\Z_E}\abs{C_{\mu_j}-C_{\nu_j}}\le\epsilon_j$ and $\int g\dd(\mu_j-\nu_j)\ge\delta$. By compactness of $\Pp(B_R)$ in $W_1$, a subsequence satisfies $\mu_j\to\mu$ and $\nu_j\to\nu$. By \eqref{eq:lipschitz},
\[
 \sup_{\Z_E}\abs{C_\mu-C_\nu}\le\epsilon_j+W_1(\mu_j,\mu)+W_1(\nu_j,\nu)\longrightarrow0,
\]
hence $C_\mu=C_\nu$ on $\Z_E$ and $\mu=\nu$ by identification. Continuity of $\nu\mapsto\int g\dd\nu$ in $W_1$ gives $\int g\dd(\mu-\nu)\ge\delta$, a contradiction.
\end{proof}

\begin{proof}[Proof of Theorem~\ref{thm:minimax-functional}]
(i) Let $D_{\mathcal C}=\{\mu-\nu:\mu,\nu\in\mathcal C\}$. Since $\mathcal C$ is a closed subset of the compact space $\Pp(K)$, $D_{\mathcal C}$ is convex, symmetric, and weak-$*$ compact in the finite signed measures on $K$, and $\osc_{\mathcal C}(h)=\sup_{\eta\in D_{\mathcal C}}\int h\dd\eta$. Apply the Fenchel--Rockafellar theorem \cite[Ch.~III, Theorem 4.1]{ekeland1976} to $V=L^2(Q)$, $Y=C(K)$, the bounded operator $\Lambda w=\cA^*(wQ)$, $F(w)=\epsilon\norm w_{L^2(Q)}$, and $G(h)=\osc_{\mathcal C}(g-h)$. Both functions are convex and finite, and $G$ is $2$-Lipschitz in the uniform norm. Hence the infimum in \eqref{eq:l2-modulus-duality} equals $\max_\eta\{-F^*(\Lambda^*\eta)-G^*(-\eta)\}$ over the finite signed measures $\eta$ on $K$. The functional $\Lambda^*\eta$ acts by $w\mapsto\int C_\eta w\dd Q$ and has norm $\norm{C_\eta}_{L^2(Q)}$, hence $F^*(\Lambda^*\eta)$ is $0$ if $\norm{C_\eta}_{L^2(Q)}\le\epsilon$ and $+\infty$ otherwise. Substituting $q=g-h$,
\[
 G^*(-\eta)=-\int g\dd\eta+\sup_{q\in C(K)}\left\{\int q\dd\eta-\sup_{\eta'\in D_{\mathcal C}}\int q\dd\eta'\right\}.
\]
The last supremum is $0$ if $\eta\in D_{\mathcal C}$. If $\eta\notin D_{\mathcal C}$, the Hahn--Banach theorem separates $\eta$ from the weak-$*$ closed convex set $D_{\mathcal C}$ by some $q\in C(K)$, and scaling $q$ makes the supremum $+\infty$. The dual problem is therefore $\sup\{\int g\dd\eta:\eta\in D_{\mathcal C},\ \norm{C_\eta}_{L^2(Q)}\le\epsilon\}=\omega_g(\epsilon;\mathcal C)$.

(ii) Let $a_\xi=\sqrt{1\vee I_\xi}$ and let $\mu,\nu\in\mathcal C$ satisfy $\norm{C_\mu-C_\nu}_{L^2(Q)}\le\varsigma/(a_\xi\sqrt n)$. By Lemma~\ref{lem:hellinger}, the joint laws of the $n$ triples have total variation distance at most $\tfrac12$. Let $a=\int g\dd\mu$, $b=\int g\dd\nu$, and $\phi=\mathbf1\{\abs{\widehat p-a}\ge\abs{a-b}/2\}$. Then $\abs{\widehat p-a}\ge\tfrac12\abs{a-b}\phi$ and $\abs{\widehat p-b}\ge\tfrac12\abs{a-b}(1-\phi)$, so
\[
 \E_\mu\abs{\widehat p-a}+\E_\nu\abs{\widehat p-b}\ge\tfrac12\abs{a-b}\left(1-\operatorname{TV}\right)\ge\tfrac14\abs{a-b},
\]
and the larger of the two risks is at least $\abs{a-b}/8$. The supremum over such pairs gives $\mathfrak R_n^Q(g;\mathcal C)\ge\omega_g(\varsigma/(a_\xi\sqrt n);\mathcal C)/8$. Since $\mathcal C$ is convex, convex combinations of feasible pairs are feasible. Hence $\omega_g$ is concave with $\omega_g(0)\ge0$, and $\omega_g(\epsilon/a_\xi)\ge\omega_g(\epsilon)/a_\xi$.

For the upper bound, fix $w\in L^2(Q)$ and $c\in\R$. Fubini's theorem gives $\E \widehat p_w=c+\int wC_\mu\dd Q=c+\int\cA^*(wQ)\dd\mu$, hence the bias of $\widehat p_w$ at $\mu$ is $\int(c+\cA^*(wQ)-g)\dd\mu$, which has absolute value at most $\tfrac12\osc_{\mathcal C}(g-\cA^*(wQ))$ when $c$ is the midpoint of the price range of $g-\cA^*(wQ)$ over $\mathcal C$. Since $0\le C_\mu\le2R$,
\[
 \operatorname{Var}_\mu(\widehat p_w)\le\frac1n\E\left[w(Z_1)^2Y_1^2\right]\le\frac{4R^2+\varsigma^2}{n}\norm w_{L^2(Q)}^2.
\]
Hence, with $\varkappa_\varsigma=\sqrt{1+4R^2/\varsigma^2}\ge1$,
\[
 \begin{aligned}
 \sup_{\mu\in\mathcal C}\E_\mu\abs{\widehat p_w-\int g\dd\mu}
 &\le\tfrac12\osc_{\mathcal C}\left(g-\cA^*(wQ)\right)
                  +\varkappa_\varsigma\frac{\varsigma}{\sqrt n}\norm w_{L^2(Q)}\\
 &\le\varkappa_\varsigma\left[\osc_{\mathcal C}\left(g-\cA^*(wQ)\right)
                  +\frac{\varsigma}{\sqrt n}\norm w_{L^2(Q)}\right].
 \end{aligned}
\]
Take the infimum over $w$ and use \eqref{eq:l2-modulus-duality}.
\end{proof}

\begin{proof}[Proof of Corollary~\ref{cor:price-ci}]
(i) The proof of Theorem~\ref{thm:minimax-functional} bounds the bias of $\widehat p_w$ by $\tfrac12\osc_{\mathcal C}(g-\cA^*(wQ))$ and its variance by $(4R^2+\varsigma^2)\norm w_{L^2(Q)}^2/n$, uniformly over $\mathcal C$. Chebyshev's inequality gives the coverage. The length $2\ell_w$ is at most $2\delta^{-1/2}\sqrt{1+4R^2/\varsigma^2}\,[\osc_{\mathcal C}(g-\cA^*(wQ))+(\varsigma/\sqrt n)\norm w_{L^2(Q)}]$, and \eqref{eq:l2-modulus-duality} gives the bound.

(ii) Let $a_\xi=\sqrt{1\vee I_\xi}$, and let $\mu,\nu\in\mathcal C$ satisfy $\norm{C_\mu-C_\nu}_{L^2(Q)}\le\varsigma/(a_\xi\sqrt n)$, so that the observation laws have total variation distance at most $\tfrac12$ by Lemma~\ref{lem:hellinger}. With $a=\int g\dd\mu$ and $b=\int g\dd\nu$,
\[
 \Pr_\mu(a\in\widehat I,\ b\in\widehat I)\ge\Pr_\mu(a\in\widehat I)-\Pr_\nu(b\notin\widehat I)-\tfrac12\ge\tfrac12-2\delta,
\]
hence $\E_\mu\abs{\widehat I}\ge(\tfrac12-2\delta)\abs{a-b}$. The supremum over such pairs and the concavity of $\omega_g$ established in the proof of Theorem~\ref{thm:minimax-functional} give the bound.
\end{proof}

\section{Proofs for Section~\ref{sec:directional-pricing}}

Throughout this section, $E$, $\be$, $g_k$ and $\eta_{\be}$ are as in Section~\ref{sec:log-modulus}.

\begin{proof}[Proof of Theorem~\ref{thm:missing-cone}]
Rotate the construction in \eqref{eq:oscillation} to the direction $\be$:
\[
 H_N(\bx)=bN^{-\alpha-2}(\be\cdot\nabla)^2
                    \left(\zeta_0(\bx)\cos(N\be\cdot\bx)\right),
 \qquad f_N^\pm=f_0\pm H_N.
\]
The positivity, moment cancellation, and $H^\alpha$ bounds in the proof of Theorem~\ref{thm:modulus} are unchanged under rotation. Testing against $N^{-1}\cos(N\be\cdot\bx)$ proves \eqref{eq:limited-w1}.

Condition~\eqref{eq:missing-cone} implies that there is $c_\kappa>0$ with
\begin{equation}\label{eq:cone-frequency}
 \abs{t\bu\pm N\be}\ge c_\kappa(\abs{t}+N),
 \qquad t\in\R,\ \bu\in E.
\end{equation}
To see this, the squared left side is at least
$t^2+N^2-2(1-\kappa)\abs{t}N$, which is bounded below by a positive multiple of $(\abs{t}+N)^2$.
The two derivatives in the perturbation cancel the denominator in \eqref{eq:slice}, giving
\[
 \widehat{\cA H_N}(\bu,t)
 =\frac b2N^{-\alpha-2}(\bu\cdot\be)^2
  \left(\widehat\zeta_0(t\bu-N\be)
                      +\widehat\zeta_0(t\bu+N\be)\right).
\]
For every $P>1$, Schwartz decay and \eqref{eq:cone-frequency} imply
\begin{align*}
 \sup_{\bu\in E}\int_\R
                  \abs{\widehat{\cA H_N}(\bu,t)}\dd t
 &\le C_PN^{-\alpha-2}
                  \int_\R(1+N+\abs{t})^{-P}\dd t\,\le C_P'N^{-\alpha-1-P}.
\end{align*}
Fourier inversion gives \eqref{eq:limited-price}, with $P$ chosen sufficiently large for each $A$. If a H\"older inequality with exponent $\gamma$ held, choose $A$ with $A\gamma>\alpha+1$ and combine it with \eqref{eq:limited-w1}; the resulting inequality fails as $N\to\infty$. Every $L^p$ norm under a probability design is bounded by the supremum norm, giving the final assertion.
\end{proof}

\begin{proof}[Proof of Corollary~\ref{cor:long-only}]
The set $E$ contains a relatively open patch, hence Theorem~\ref{thm:direction-identification} gives identification. For $\be=(\be_1-\be_2)/\sqrt2$ and $\bu\in E$,
$\abs{\bu\cdot\be}=\abs{u_1-u_2}/\sqrt2\le1/\sqrt2$.
Thus \eqref{eq:missing-cone} holds with $\kappa=1-1/\sqrt2$.
\end{proof}

\begin{proof}[Proof of Proposition~\ref{prop:positive}]
(i) For $y\in[-R,R]$, the spanning formula \cite{breeden1978,carr2001} gives $G(y)=G(-R)+G'(-R+)(y+R)+\int_{(-R,R)}(y-k)_+G''(\mathrm dk)$, and $\bv\cdot\bx+R=\phi_{(\bv,-R)}(\bx)$ on $B_R$. Hence $g=G(-R)+\cA^*m$ with $\norm m$ equal to the stated constant, and Lemma~\ref{lem:replication-bound} applies.

(ii) Let $N=\lceil\sqrt{2R/\epsilon}\,\rceil$ and $h=2R/N\le\sqrt{2R\epsilon}$. The piecewise linear interpolant $I$ of $G$ on the uniform grid of mesh $h$ satisfies $\norm{G-I}_\infty\le h$, and its slopes $a_1,\ldots,a_N$ lie in $[-1,1]$. By (i), $I(\bv\cdot\bx)=c+\cA^*m$ with $\norm m\le\abs{a_1}+\sum_{j\ge2}\abs{a_j-a_{j-1}}\le2N-1$. Lemma~\ref{lem:replication-bound} gives $\varpi_g(\epsilon)\le2h+\epsilon(2N-1)\le4\sqrt{2R\epsilon}+\epsilon$.

(iii) By (i), $(\bv\cdot\bx)^p\in V_E$ for every $\bv\in E$ and $p\ge0$. Let $U$ be a nonempty relatively open subset of $E$. Under the Fischer inner product $\langle\bx^a,\bx^b\rangle=a!\,\mathbf1\{a=b\}$ on homogeneous polynomials of degree $p$, $\langle q,(\bv\cdot\bx)^p\rangle=p!\,q(\bv)$. A homogeneous $q$ orthogonal to every $(\bv\cdot\bx)^p$ with $\bv\in U$ vanishes on $U$, hence on the open cone it generates, hence identically. Thus these powers span all homogeneous polynomials of degree $p$.
\end{proof}

\begin{proof}[Derivation for Example~\ref{ex:polynomial-rate}]
For the upper bound, we construct $w\in L^2(Q)$ with $g=c+\cA^*(wQ)$ on $B_r$. For $\bv\in E$ and integers $p\ge2$, the univariate spanning formula gives
\[
 \begin{split}
 (\bv\cdot\bx)^p={}&(-r)^p+p(-r)^{p-1}(\bv\cdot\bx+r)\,+\int_{-r}^{r}(\bv\cdot\bx-k)_+\,p(p-1)k^{p-2}\dd k.
 \end{split}
\]
The linear case is represented using the support buffer:
\[
 \bv\cdot\bx
 =\frac1{R-r}\int_{-R}^{-r}\phi_{(\bv,k)}(\bx)\dd k
                  -\frac{R+r}{2},\qquad \bx\in B_r.
\]
Consequently, for every $p\ge1$ there are a constant $c_p$ and $w_p\in L^\infty[-R,R]$, independent of $\bv$, such that $(\bv\cdot\bx)^p=c_p+\int\phi_{(\bv,k)}(\bx)w_p(k)\dd k$. For a nonempty relatively open $U\subseteq E$, the map
\[
 \beta\longmapsto\int_U(\bv\cdot\bx)^p\beta(\bv)\,
                           \mathcal H^{d-1}(\mathrm d\bv)
\]
from $L^2(U)$ to the homogeneous polynomials of degree $p$ is surjective. Indeed, a polynomial Fischer-orthogonal to its range vanishes on $U$, hence identically, as in Proposition~\ref{prop:positive}(iii). Every polynomial therefore equals $c+\int_{\Z_E}\phi_z w_E(z)\,Q_E(\mathrm dz)$ with $w_E\in L^2(Q_E)$. Set $w=w_E\,\mathrm dQ_E/\mathrm dQ$. Then
\[
 \norm w_{L^2(Q)}^2\le q_*^{-1}\norm{w_E}_{L^2(Q_E)}^2,
 \qquad \omega_g(\epsilon;\Pp(B_r))\le\epsilon\norm w_{L^2(Q)},
\]
and Theorem~\ref{thm:minimax-functional} gives the upper rate. For the lower rate, choose $\mu,\nu$ with $\int g\dd(\mu-\nu)\ne0$, and write $\eta=\mu-\nu$. The representation implies $\norm{C_\eta}_{L^2(Q)}>0$. The pair
\[
 \mu_t=\tfrac12(\mu+\nu)+\tfrac t2\eta,\qquad
 \nu_t=\tfrac12(\mu+\nu)-\tfrac t2\eta,\qquad 0\le t\le1,
\]
gives $\omega_g(\epsilon;\Pp(B_r))\ge c\epsilon$ for all sufficiently small $\epsilon$. The lower bound in \eqref{eq:minimax-functional} then gives $\mathfrak R_n^Q(g;\Pp(B_r))\ge c'n^{-1/2}$.

To estimate correlation, fix $v_*>0$ and restrict attention to the nonempty class
\[
 \mathcal C(v_*)=\left\{\mu\in\Pp(B_r):
       \operatorname{Var}_\mu(X_j)\ge v_*,\quad j=1,2\right\}.
\]
Let $m_j=\E_\mu X_j$, $s_j=\E_\mu X_j^2$, and $q_{12}=\E_\mu X_1X_2$. The polynomial estimators above, clipped respectively to $[-r,r]$, $[0,r^2]$, and $[-r^2,r^2]$, give estimates $\widehat m_j,\widehat s_j,\widehat q_{12}$ whose expected absolute errors are uniformly $O(n^{-1/2})$. Define
\[
 \widehat v_j=\max\left(\frac{v_*}{2},\widehat s_j-\widehat m_j^2\right),
 \qquad
 \widehat\rho=\operatorname{clip}_{[-1,1]}
    \left(\frac{\widehat q_{12}-\widehat m_1\widehat m_2}
                 {\sqrt{\widehat v_1\widehat v_2}}\right).
\]
For $v_j=s_j-m_j^2\ge v_*$,
\[
 \abs{\widehat v_j-v_j}\le\abs{\widehat s_j-s_j}
                                    +2r\abs{\widehat m_j-m_j},
\]
and
\[
 \abs{\widehat q_{12}-\widehat m_1\widehat m_2-(q_{12}-m_1m_2)}
 \le\abs{\widehat q_{12}-q_{12}}
       +r\sum_{j=1}^2\abs{\widehat m_j-m_j}.
\]
The map $(c,v_1,v_2)\mapsto c/(v_1v_2)^{1/2}$ is Lipschitz on the bounded range in question with $v_1,v_2\ge v_*/2$. Clipping cannot increase its error relative to a correlation in $[-1,1]$. Thus
\[
 \sup_{\mu\in\mathcal C(v_*)}
 \E_\mu\abs{\widehat\rho-\operatorname{Corr}_\mu(X_1,X_2)}
 \le C_{r,v_*,Q,\varsigma}n^{-\frac{1}{2}}.
\]
\end{proof}

\begin{proof}[Proof of Corollary~\ref{cor:finite}]
For $\mu,\nu\in\K_n(z,v)$, $C_\mu-C_\nu$ vanishes on $z$ and is $2R$-Lipschitz in $\rho$ by \eqref{eq:lipschitz}, hence $\sup_{\Z_E}\abs{C_\mu-C_\nu}\le2Rh_E(z)$. For the random bound, \eqref{eq:ahlfors} gives $Q_E(B_\rho(w,t))\ge ct^d$ for $w\in\Z_E$ and $t\le1$, and $\Z_E$ has a $t$-net of cardinality at most $C_dt^{-d}$. The proof of Lemma~\ref{lem:coverage} applies verbatim.
\end{proof}

\begin{proof}[Proof of Corollary~\ref{cor:missing-noisy}]
Let $g_0(\bx)=(\be\cdot\bx)_+$, which is $1$-Lipschitz, so every estimator $\widehat\mu$ satisfies $\abs{\int g_0\dd\widehat\mu-\int g_0\dd\mu}\le W_1(\widehat\mu,\mu)$ and $\mathfrak R_n^Q(\alpha,L)\ge\mathfrak R_n^Q(g_0;\mathcal C_{\alpha,L})$. Corollary~\ref{cor:log-minimax}(i) gives $\mathfrak R_n^Q(g_0;\mathcal C_{\alpha,L})\ge c(\log n)^{-(\alpha+3/2)}$ for large $n$,which proves the claimed divergence. The long-only assertion follows from Corollary~\ref{cor:long-only}.
\end{proof}

\begin{lemma}\label{lem:line-restriction}
There is $C=C(d,r,\alpha)$ such that every $\eta\in H^\alpha(\R^d)$ supported in $B_r$ satisfies
\[
 \int_\R\abs{\widehat\eta(t\be)}^2(1+t^2)^\alpha\dd t\le C\norm\eta_{H^\alpha}^2.
\]
\end{lemma}

\begin{proof}
Rotate so that $\be=\be_1$, write $\bx=(y,\bx')$, and let $\mathcal F_y\eta(t,\bx')$ denote the Fourier transform of $\eta$ in $y$ alone. Since $\eta$ is supported in $B_r$, $\widehat\eta(t\be)=\int_{\abs{\bx'}\le r}\mathcal F_y\eta(t,\bx')\dd\bx'$, and Cauchy--Schwarz gives $\abs{\widehat\eta(t\be)}^2\le C_{d,r}\int\abs{\mathcal F_y\eta(t,\bx')}^2\dd\bx'$. Plancherel's identity in $\bx'$ and $(1+t^2)^\alpha\le(1+t^2+\abs{\bxi'}^2)^\alpha$ give
\[
 \int_\R(1+t^2)^\alpha\abs{\widehat\eta(t\be)}^2\dd t
 \le C_{d,r}\int_{\R^d}(1+\abs\bxi^2)^\alpha\abs{\widehat\eta(\bxi)}^2\dd\bxi
 \le C\norm\eta_{H^\alpha}^2.
\]
\end{proof}

\begin{lemma}[Extrapolation of powers]\label{lem:long-only-powers}
Let $U\subseteq S^{d-1}$ be closed with nonempty relative interior and let $\be\in S^{d-1}$. There is $\Lambda=\Lambda(U,\be)\ge1$ such that for every integer $i\ge1$ there is a finite signed measure $m_i$ on $U$ with $\norm{m_i}\le\Lambda^i$ and
\[
 (\be\cdot\bx)^i=\int_U(\bv\cdot\bx)^i\,m_i(\mathrm d\bv),\qquad\bx\in\R^d.
\]
\end{lemma}

\begin{proof}
Choose $\bu_0$ in the relative interior of $U$ with $\bu_0\cdot\be\ne0$, which is possible since $\be^\perp\cap S^{d-1}$ has empty relative interior, and $\rho\in(0,1]$ such that the closed cap $U_0=\{\bv\in S^{d-1}:\abs{\bv-\bu_0}\le\rho\}$ lies in $U$. Let $\mathsf H_i$ be the space of homogeneous polynomials of degree $i$ on $\R^d$. The powers $(\bv\cdot\bx)^i$, $\bv\in U_0$, form a compact set spanning $\mathsf H_i$, by the argument of Proposition~\ref{prop:positive}(iii). By the bipolar theorem in $\mathsf H_i$, the least total variation of a representing measure on $U_0$ equals $\sup\abs{\lambda((\be\cdot\bx)^i)}$ over linear functionals $\lambda$ on $\mathsf H_i$ with $\abs{\lambda((\bv\cdot\bx)^i)}\le1$ for $\bv\in U_0$, and the infimum is attained. Every such $\lambda$ has the form $\lambda(P)=\langle q,P\rangle/i!$ for some $q\in\mathsf H_i$ under the Fischer inner product, and then $\lambda((\bv\cdot\bx)^i)=q(\bv)$. It therefore suffices to show $\abs{q(\be)}\le\Lambda^i\sup_{U_0}\abs q$ for $q\in\mathsf H_i$.

Let $\boldsymbol p_{\be}=\be/(\be\cdot\bu_0)$ and $D=\abs{\boldsymbol p_{\be}-\bu_0}$. Since $\boldsymbol p_{\be}\cdot\bu_0=1$, the vector $\boldsymbol p_{\be}-\bu_0$ is orthogonal to $\bu_0$. If $D=0$, then $\be=\pm\bu_0$ and $\abs{q(\be)}\le\sup_{U_0}\abs q$. Otherwise let $P(s)=q(\bu_0+sD^{-1}(\boldsymbol p_{\be}-\bu_0))$, a polynomial of degree at most $i$. For $\abs s\le\rho$, the point $\by=\bu_0+sD^{-1}(\boldsymbol p_{\be}-\bu_0)$ satisfies $\abs\by^2=1+s^2$ and $\abs{\by/\abs\by-\bu_0}^2=2-2(1+s^2)^{-1/2}\le s^2$, hence $\by/\abs\by\in U_0$ and, by homogeneity,
\[
 \abs{P(s)}=\abs\by^i\abs{q(\by/\abs\by)}\le2^{i/2}\sup_{U_0}\abs q,\qquad\abs s\le\rho.
\]
If $D\le\rho$, this bounds $\abs{P(D)}$. If $D>\rho$, Chebyshev's extremal property \cite[Chapter 2]{borwein1995} gives $\abs{P(D)}\le\abs{T_i(D/\rho)}\sup_{\abs s\le\rho}\abs{P(s)}$, where $T_i$ is the Chebyshev polynomial, and $\abs{T_i(x)}\le(2x)^i$ for $x\ge1$. Since $q(\be)=(\be\cdot\bu_0)^iq(\boldsymbol p_{\be})=(\be\cdot\bu_0)^iP(D)$ and $\abs{\be\cdot\bu_0}\le1$,
\[
 \abs{q(\be)}\le\Lambda^i\sup_{U_0}\abs q,
 \qquad\Lambda=\sqrt2\max(1,2D/\rho).
\]
\end{proof}

\begin{lemma}\label{lem:moments}
Let $\eta$ be a finite signed measure on $B_R$ with $\eta(B_R)=0$ and $\sup_{\Z_E}\abs{C_\eta}\le\epsilon$. For $\abs k\le R$ and every integer $j\ge1$,
\[
 \abs{\int(\be\cdot\bx-k)^j\,\eta(\mathrm d\bx)}\le\frac{3j}{R}\,\varrho^j\epsilon,\qquad \varrho=R(\Lambda+1),
\]
where $\Lambda=\Lambda(E,\be)$ is the constant of Lemma~\ref{lem:long-only-powers}.
\end{lemma}

\begin{proof}
For $\bv\in E$ and $i\ge1$, Proposition~\ref{prop:positive}(i) with $G(y)=y^i$ writes $(\bv\cdot\bx)^i$ as a constant plus $\cA^*m$ with $\norm m\le3iR^{i-1}$, hence $\abs{\int(\bv\cdot\bx)^i\dd\eta}\le3iR^{i-1}\epsilon$. Lemma~\ref{lem:long-only-powers} with $U=E$ gives $\abs{\int(\be\cdot\bx)^i\dd\eta}\le3i\Lambda^iR^{i-1}\epsilon$. Expanding $(\be\cdot\bx-k)^j$ binomially, the constant term integrates to zero, and the remaining terms sum to at most $(3j/R)(\abs k+\Lambda R)^j\epsilon$.
\end{proof}

\begin{lemma}[Compactly supported profiles]\label{lem:profile}
Let $d'\ge1$, let $\zeta_0\in C_c^\infty(\R^{d'})$ be nonnegative with $\int\zeta_0>0$, and let $a_\perp>0$. For $N\ge1$ let $u_N=(N/2a_\perp)\mathbf1_{[-a_\perp/N,a_\perp/N]}$, $\kappa_N(\bx)=\prod_{j=1}^{d'}u_N^{*N}(x_j)$, and $\zeta_N=\zeta_0*\kappa_N$. Then $\zeta_N\ge0$ is smooth, supported in $\supp\zeta_0+[-a_\perp,a_\perp]^{d'}$, with $\int\zeta_N=\int\zeta_0$ and $\norm{\zeta_N}_{H^{\alpha'}}\le\norm{\zeta_0}_{H^{\alpha'}}$ for every $\alpha'\ge0$. For every $c_0>0$ and $J\ge0$ there are $C,c>0$, independent of $N$, with
\[
 \abs{\widehat\zeta_N(\bxi)}\le Ce^{-cN}(1+\abs\bxi)^{-J},\qquad\abs\bxi\ge c_0N.
\]
\end{lemma}

\begin{proof}
$\kappa_N$ is a probability density supported in $[-a_\perp,a_\perp]^{d'}$, which gives the support, sign, mass, and Sobolev assertions by Young's inequality. Its transform is $\widehat\kappa_N(\bxi)=\prod_j\operatorname{sinc}(a_\perp\xi_j/N)^N$ with $\operatorname{sinc}v=\sin v/v$. If $\abs\bxi\ge c_0N$, some coordinate satisfies $\abs{\xi_j}\ge c_0N/\sqrt{d'}$, hence $\abs{a_\perp\xi_j/N}\ge\delta:=a_\perp c_0/\sqrt{d'}$ and $\abs{\widehat\kappa_N(\bxi)}\le q_\delta^N$ with $q_\delta=\sup_{\abs v\ge\delta}\abs{\operatorname{sinc}v}<1$. Since $\abs{\widehat\zeta_0(\bxi)}\le C_J(1+\abs\bxi)^{-J}$, the claim holds with $c=\log(1/q_\delta)$.
\end{proof}

\begin{lemma}\label{lem:window}
Fix $r'>0$. For each $N\ge1$ there is $\mathsf K_N\in C_c^\infty(-r',r')$ with $\mathsf K_N(0)=1$ such that, with constants depending only on $r'$ and $m$,
\begin{enumerate}[label=(\alph*)]
\item $\norm{\mathsf K_N}_{L^2}^2\le C/N$ and $\norm{\mathsf K_N^{(m)}}_\infty\le C_mN^m$;
\item $\norm{\widehat{\mathsf K}_N}_\infty\le C/N$;
\item $\abs{\widehat{\mathsf K}_N(\tau)}\le C_mN^me^{-r'N/256}(1+\abs\tau)^{-m}$ for $\abs\tau\ge N/8$.
\end{enumerate}
\end{lemma}

\begin{proof}
Let $\chi_0\in C_c^\infty(-1,1)$ with $0\le\chi_0\le1$ and $\chi_0=1$ on $[-\tfrac12,\tfrac12]$. Set $\mathsf v_N^2=N/(32r')$, $\mathcal W_N(y)=e^{-\mathsf v_N^2y^2}\chi_0(y/r')$, $B_N(y)=\sin(Ny/16)/(Ny/16)$, and $\mathsf K_N=\mathcal W_NB_N$. Then $\mathsf K_N(0)=1$, $\abs{\mathsf K_N}\le\abs{B_N}$, and $\norm{B_N}_{L^2}^2=16\pi/N$; the derivative bounds follow from Leibniz's rule. Next, $\widehat{\mathsf K}_N=(2\pi)^{-1}\widehat{\mathcal W}_N*\widehat B_N$ with $\widehat B_N=(16\pi/N)\mathbf1_{[-N/16,N/16]}$, and $\widehat{\mathcal W}_N=\mathsf g_N+\mathcal R_N$. Here $\mathsf g_N(\tau)=(\sqrt\pi/\mathsf v_N)e^{-\tau^2/(4\mathsf v_N^2)}$ has $\norm{\mathsf g_N}_{L^1}=2\pi$, and $\mathcal R_N$ is the transform of the smooth function $e^{-\mathsf v_N^2y^2}(\chi_0(y/r')-1)$, supported in $\abs y\ge r'/2$. Its derivatives give $\abs{\mathcal R_N(\tau)}\le C_m\mathsf v_N^me^{-\mathsf v_N^2r'^2/8}(1+\abs\tau)^{-m}$, and $\mathsf v_N^2r'^2/8=r'N/256$. Hence $\norm{\widehat{\mathcal W}_N}_{L^1}\le C$, which gives (b). For $\abs\tau\ge N/8$ and $\abs\zeta\le N/16$, $\abs{\tau-\zeta}\ge\max(N/16,\abs\tau/2)$, hence $\mathsf g_N(\tau-\zeta)\le C_m\mathsf v_N^me^{-r'N/64}(1+\abs\tau)^{-m}$; since $(2\pi)^{-1}\norm{\widehat B_N}_{L^1}=1$, this gives (c).
\end{proof}

\begin{proof}[Proof of the quoted-direction bound in Theorem~\ref{thm:log-modulus}]
If $\be\in E$, then $g_k=\phi_{(\be,k)}$. If $-\be\in E$, put--call parity gives $g_k=\phi_{(-\be,-k)}-\phi_{(-\be,-R)}+R-k$. 
Lemma~\ref{lem:replication-bound} gives $\varpi_{g_k}(\epsilon)\le2\epsilon$.
\end{proof}

\begin{proof}[Proof of the upper bound in \eqref{eq:log-modulus}]
If $\pm\be\notin E$, compactness of $E$ gives the missing-cone condition with $\kappa=1-\max_{\bu\in E}\abs{\bu\cdot\be}>0$. Let $f_1,f_2\in\F_\alpha(L)$ with $\sup_{\Z_E}\abs{C_{f_1}-C_{f_2}}\le\epsilon$, and let $\eta=(f_1-f_2)\dd\bx$, hence $\norm\eta_{H^\alpha}\le2L$ and $\supp\eta\subseteq B_r$. Since $\eta$ has zero mass, $(y-k)_+=\tfrac12(y-k)+\tfrac12\abs{y-k}$ and $\abs y=\pi^{-1}\int_\R(1-\cos ty)t^{-2}\dd t$ give
\[
 \int g_k\dd\eta=\frac{M_1}2-\frac1{2\pi}\int_\R\frac{\mathsf c(t)}{t^2}\dd t,\qquad
 M_1=\int\be\cdot\bx\,\eta(\mathrm d\bx),\quad \mathsf c(t)=\int\cos(t(y-k))\,\eta_{\be}(\mathrm dy),
\]
with Fubini's theorem justified by $\abs{1-\cos ty}t^{-2}\le\min(y^2,2t^{-2})$. Lemma~\ref{lem:moments} with $j=1$ and $k=0$ gives $\abs{M_1}\le3(\Lambda+1)\epsilon$.

For $\abs t\le N$, expand $\mathsf c(t)=\sum_{j\ge1}(-1)^jt^{2j}\mathsf m_{2j}/(2j)!$ with $\mathsf m_j=\int(\be\cdot\bx-k)^j\dd\eta$. Lemma~\ref{lem:moments} gives
\[
 \frac{\abs{\mathsf c(t)}}{t^2}\le\frac{3\varrho^2\epsilon}{R}\sum_{j\ge1}\frac{(\varrho\abs t)^{2j-2}}{(2j-1)!}\le\frac{3\varrho^2\epsilon}{R}e^{\varrho\abs t},
 \qquad
 \int_{\abs t\le N}\frac{\abs{\mathsf c(t)}}{t^2}\dd t\le\frac{6\varrho\epsilon}{R}e^{\varrho N}.
\]
For $\abs t>N$, $\abs{\mathsf c(t)}\le\abs{\widehat\eta(t\be)}$, and Cauchy--Schwarz with Lemma~\ref{lem:line-restriction} gives
\[
 \int_{\abs t>N}\frac{\abs{\mathsf c(t)}}{t^2}\dd t\le\left(\int\abs{\widehat\eta(t\be)}^2(1+t^2)^\alpha\dd t\right)^{\frac{1}{2}}\left(\int_{\abs t>N}t^{-4-2\alpha}\dd t\right)^{\frac{1}{2}}\le CLN^{-\alpha-\frac{3}{2}}.
\]
Taking $N=\log(1/\epsilon)/(2\varrho)$ gives \(\abs{\int g_k\dd\eta} \le C\left( \epsilon^{1/2} +L(\log(1/\epsilon))^{-(\alpha+3/2)} \right), \) which proves the upper bound in \eqref{eq:log-modulus} for sufficiently small $\epsilon$.
\end{proof}

\begin{proof}[Proof of the lower bound in \eqref{eq:log-modulus}]
Write $\bx=y\be+\bx'$ with $\bx'\perp\be$ and $p=\alpha+\tfrac32$. Choose $r'>0$ and $r_\perp>0$ with $r_\perp^2+(\abs k+r')^2<r^2$. Let $\psi_N$ be the profile of Lemma~\ref{lem:profile} on $\be^\perp\cong\R^{d-1}$, built from a nonnegative $\psi_0\in C_c^\infty$ and an $a_\perp$ for which every $\psi_N$ is supported in $\{\abs{\bx'}<r_\perp\}$; write $M_\psi=\int\psi_N=\int\psi_0>0$. Let $\mathsf K_N$ be the window of Lemma~\ref{lem:window} with this $r'$, and $\Phi_N(y)=\cos(Ny)\mathsf K_N(y)$.

\emph{Envelope.} In Lemma~\ref{lem:window}, $\abs{B_N^{(j)}(y)}\le C_jN^j(1+N\abs y)^{-1}$ for $j\le2$, and the Gaussian factor with its first two derivatives is bounded by $C_jN^{j/2}e^{-c_1Ny^2}$, $c_1=1/(64r')$. Leibniz's rule gives $\abs{\Phi_N''(y)}\le CN^2(1+N^2y^2)^{-1/2}e^{-c_1Ny^2}$ on $\supp\mathsf K_N$. Let $\bar\chi\in C_c^\infty(-r',r')$ equal one on $\supp\mathsf K_N$ and
\[
E_N(y)=C_0N^{\frac{1}{2}-\alpha}(1+N^2y^2)^{-\frac{1}{2}}e^{-c_1Ny^2}\bar\chi(y),
\]
with $C_0$ fixed so that $E_N\ge N^{-p}\abs{\Phi_N''}$. Writing $E_N(y)=C_0N^{1/2-\alpha}h_N(Ny)$ with $h_N(t)=(1+t^2)^{-1/2}e^{-c_1t^2/N}\bar\chi(t/N)$, every derivative of $h_N$ is bounded in $L^2(\R)$ uniformly in $N$, and $\norm{h_N(N\,\cdot)}_{H^\alpha}^2\le N^{2\alpha-1}\norm{h_N}_{H^\alpha}^2$ gives $\sup_N\norm{E_N}_{H^\alpha}<\infty$, while $\int E_N\le CN^{-\alpha-1/2}\log(2+N)$.

\emph{Densities.} Fix a smooth probability density $f_*$ supported in $B_r$ and set
\[
 \eta_N(\bx)=bN^{-p}\psi_N(\bx')\Phi_N''(y-k),\qquad
 m_N=bM_\psi\int E_N,
\]
\[
 f_N^\pm=(1-m_N)f_*+b\,\psi_N(\bx')E_N(y-k)\pm\eta_N.
\]
For small $b$, $m_N\le\tfrac12$ for all $N$. Both functions are nonnegative since $E_N\ge N^{-p}\abs{\Phi_N''}$ and $\psi_N\ge0$; $\eta_N$ has zero mass and zero first moments, by integration by parts in $y$, so both have mass one; their supports lie in $B_r$. Since $(1+\abs{\bxi'}^2+\tau^2)^\alpha\le(1+\abs{\bxi'}^2)^\alpha(1+\tau^2)^\alpha$, products satisfy $\norm{a(\bx')b(y)}_{H^\alpha(\R^d)}\le C\norm a_{H^\alpha(\R^{d-1})}\norm b_{H^\alpha(\R)}$. Lemma~\ref{lem:profile} bounds $\norm{\psi_N}_{H^\alpha}$; on $\{\abs{\abs\tau-N}\le N/8\}$, Lemma~\ref{lem:window}(a) gives $\int\tau^4(1+\tau^2)^\alpha\abs{\widehat\Phi_N}^2\le CN^{2\alpha+3}$, and Lemma~\ref{lem:window}(c) makes the remaining contribution bounded, so $\norm{\Phi_N''}_{H^\alpha}\le CN^p$. Hence $\norm{f_N^\pm}_{H^\alpha}\le\norm{f_*}_{H^\alpha}+Cb$, and for $L$ large and $b$ small, $f_N^\pm\in\F_\alpha(L)$ for all $N\ge1$.

\emph{Price gap.} Integrating by parts twice in $y$,
\[
 \int g_k(f_N^+-f_N^-)\dd\bx=2bN^{-p}M_\psi\Phi_N(0)=2bM_\psi N^{-p}.
\]

\emph{Quote discrepancy.} As in \eqref{eq:slice}, for $\bu\in E$ write $\bu'$ for the component of $\bu$ orthogonal to $\be$; the missing-cone condition gives $\abs{\bu'}\ge\kappa'=(1-(1-\kappa)^2)^{1/2}$. Then
\[
 \widehat{\cA\eta_N}(\bu,t)=bN^{-p}(\bu\cdot\be)^2\,\widehat\psi_N(t\bu')\,\widehat\Phi_N(t\,\bu\cdot\be)\,e^{-\mathrm ikt\,\bu\cdot\be},
\]
where $\widehat\Phi_N(\tau)=\tfrac12(\widehat{\mathsf K}_N(\tau-N)+\widehat{\mathsf K}_N(\tau+N))$.
For $\abs t\le N/2$, $\abs{t\,\bu\cdot\be\mp N}\ge N/2$, and Lemma~\ref{lem:window}(c) bounds the integral over this region by $CN^{m+1}e^{-r'N/256}$. For $\abs t>N/2$, $\abs{t\bu'}\ge\kappa'N/2$, and Lemma~\ref{lem:profile} with $c_0=\kappa'/2$ and $J=2$ gives $\abs{\widehat\psi_N(t\bu')}\le Ce^{-cN}(1+\kappa'\abs t)^{-2}$, while $\norm{\widehat\Phi_N}_\infty\le C/N$. Fourier inversion in $t$ gives
\[
 \epsilon_N:=\sup_{\Z_E}\abs{C_{f_N^+}-C_{f_N^-}}\le Ce^{-cN}.
\]
For small $\epsilon$, choose the integer $N$ with $Ce^{-cN}\le\epsilon<Ce^{-c(N-1)}$. Then $\epsilon_N\le\epsilon$ and $N\le C'\log(1/\epsilon)$, so monotonicity of $\varpi$ and the price gap give $\varpi_{g_k}^{\alpha,L}(\epsilon)\ge2bM_\psi N^{-p}\ge c(\log(1/\epsilon))^{-p}$.

The exponent $\alpha+\tfrac32$ is the Sobolev scale of the longitudinal packet $\Phi_N''$. For $\alpha<\tfrac12$ its amplitude $N^{2-p}$ exceeds any fixed baseline, and the ridge $b\,\psi_NE_N$, whose $H^\alpha$ norm stays bounded since the ridge has width $N^{-1}$, restores positivity.
\end{proof}

\begin{proof}[Proof of Corollary~\ref{cor:log-minimax}]
(i) If $g_E=c''+\cA^*m_E$, then $\omega_g(\epsilon;\mathcal C_{\alpha,L})\ge\varpi_g^{\alpha,L}(\epsilon)\ge\abs c\,\varpi_{g_k}^{\alpha,L}(\epsilon)-\norm{m_E}\epsilon$. The lower bound in \eqref{eq:minimax-functional} with $\epsilon=\varsigma/\sqrt n$ and the lower bound in \eqref{eq:log-modulus} give the assertion once $\varsigma/\sqrt n\le\epsilon_0$.

(ii) By the upper bound in \eqref{eq:minimax-functional}, it suffices to show $\omega_g(\epsilon;\mathcal C_{\alpha,L})\le C(\log(1/\epsilon))^{-(\alpha+3/2)}$ for small $\epsilon$. Let $G=C_\mu-C_\nu$ with $\mu,\nu\in\mathcal C_{\alpha,L}$ and $\norm G_{L^2(Q)}\le\epsilon$, and let $A=\sup_{\Z_E}\abs G$. By \eqref{eq:lipschitz}, $G$ is $2R$-Lipschitz in $\rho$, so $\abs G\ge A/2$ on the $\rho$-ball of radius $A/(4R)$ about a maximizer. For small $\epsilon$ this radius is at most one, and \eqref{eq:ahlfors} gives $\epsilon^2\ge\norm G_{L^2(Q)}^2\ge q_*c(A/4R)^d(A/2)^2$, hence $A\le C\epsilon^{2/(d+2)}$. Thus
\[
 \omega_g(\epsilon;\mathcal C_{\alpha,L})\le\varpi_g^{\alpha,L}(C\epsilon^{2/(d+2)})\le\abs c\,\varpi_{g_k}^{\alpha,L}(C\epsilon^{2/(d+2)})+C\norm{m_E}\epsilon^{2/(d+2)},
\]
and the upper bound in \eqref{eq:log-modulus} gives the bound. At $\epsilon=\varsigma/\sqrt n$, $\log(1/\epsilon)\asymp\log n$.
\end{proof}

\begin{lemma}[Comparison with exact quotes]\label{lem:noise-monotone}
For every $\varsigma\ge0$, $\mathfrak R_n^Q(g;\mathcal C,0)\le\mathfrak R_n^Q(g;\mathcal C,\varsigma)$.
\end{lemma}

\begin{proof}
Adding $\varsigma\xi_i'$, with $\xi_i'$ independent copies of $\xi_1$, to exact quotes produces data with the law of the experiment at level $\varsigma$. Let $\widehat p$ be an estimator for level $\varsigma$, truncated to the range of $\int g\dd\mu$ over $\mathcal C$, which does not increase its risk. Then $\widetilde p(z,v)=\E_{\xi'}\widehat p(z,v+\varsigma\xi')$ is a measurable function of exact quotes, and Jensen's inequality for the absolute loss shows that its risk at every $\mu$ is at most that of $\widehat p$.
\end{proof}

\begin{lemma}[Bayesian bound for a linear functional]\label{lem:functional-nullspace}
Let $T\in\R^{n\times M}$ be deterministic, let $\boldsymbol\Theta$ have the product density \eqref{eq:compact-prior}, and let $\boldsymbol a\in\R^M$. For every unit vector $\bv\in\ker(T)$ and every measurable $\widehat\psi:\R^n\to\R$,
\[
 \E\left[\left(\widehat\psi(T\boldsymbol\Theta)-\boldsymbol a\cdot\boldsymbol\Theta\right)^2\right]\ge\frac{(\boldsymbol a\cdot\bv)^2}{\pi^2}.
\]
\end{lemma}

\begin{proof}
It suffices to treat bounded $\widehat\psi$, since truncation at levels $\pm\sup\abs{\boldsymbol a\cdot\btheta}$ over the cube does not increase the error. As in the proof of Lemma~\ref{lem:prior-nullspace}, $\int\widehat\psi(T\btheta)\,\partial_{\bv}p(\btheta)\dd\btheta=0$ since $\widehat\psi(T\btheta)$ is constant along lines parallel to $\bv$, while $\int(\boldsymbol a\cdot\btheta)\,\partial_{\bv}p(\btheta)\dd\btheta=-\boldsymbol a\cdot\bv$. With $\partial_{\bv}p=p\,\bv\cdot\boldsymbol S$, this gives $\E[(\widehat\psi(T\boldsymbol\Theta)-\boldsymbol a\cdot\boldsymbol\Theta)\,\bv\cdot\boldsymbol S(\boldsymbol\Theta)]=\boldsymbol a\cdot\bv$. Cauchy--Schwarz and $\E(\bv\cdot\boldsymbol S)^2=\pi^2$, from \eqref{eq:prior-score}, prove the claim.
\end{proof}

\begin{lemma}[Perturbations with zero exchange price]\label{lem:absorbing}
Let $0<r_0<r_1<r$ and suppose $Q\le q^*Q_E$. Fix a nonnegative $\zeta\in C_c^\infty(\R^d)$ supported in the ball of radius $\tfrac14$ about the origin, with $\int\zeta=1$, and fix $\varepsilon_0>0$. There are $C_3$ and $n_0$, depending only on $d,R,r_0,r_1,\alpha,q^*,\zeta,\varepsilon_0$, $E$, and $\boldsymbol w$, such that for every $n\ge n_0$ the following holds. There are functions $b_1,\ldots,b_K\in C_c^\infty(B_{r_1})$ with $K\le n^{C_3}$, depending on $n$ but not on the design, such that
\begin{enumerate}[label=(\alph*)]
\item $\int b_j\dd\bx=0$ and $\int g_kb_j\dd\bx=0$ for all $j$ and all $k\in\R$;
\item $\norm{\sum_j\theta_jb_j}_\infty+\norm{\sum_j\theta_jb_j}_{H^\alpha}\le\varepsilon_0$ for every $\btheta\in[-1,1]^K$;
\item there is an event $G_n$, depending only on $Z_1,\ldots,Z_n$, with $\Pr(G_n)\ge\tfrac34$, on which every finite signed measure $\eta$ supported in $B_{r_0}$ with zero mass and zero first moments admits $\boldsymbol c\in\R^K$ with
\[
 \sum_{j=1}^Kc_jC_{b_j}(Z_i)=C_\eta(Z_i)\quad(1\le i\le n),
 \qquad
 \abs{\boldsymbol c}_2\le n^{C_3}\max_{i\le n}\abs{C_\eta(Z_i)}.
\]
\end{enumerate}
\end{lemma}

\begin{proof}
Set $w_*=\min_{\bu\in E}\bu\cdot\boldsymbol w>0$.

\emph{The family.} Let $\Delta\in(0,1)$ and $\beta>0$ be chosen below, let $\bx_1,\ldots,\bx_K$ enumerate the points of $\Delta\mathbb Z^d$ in $B_{r_1-\Delta}$, so that $K\le C\Delta^{-d}$, and put $\zeta_j(\bx)=\zeta((\bx-\bx_j)/\Delta)$ and $b_j=\beta\Delta\,\partial_{\boldsymbol w}\zeta_j$. The supports of the $\zeta_j$ are disjoint subsets of $B_{r_1}$. Integrating by parts, $\int b_j=0$ and
\[
 \int g_kb_j\dd\bx=-\beta\Delta\,(\boldsymbol w\cdot\be)\int\mathbf1\{\be\cdot\bx>k\}\zeta_j\dd\bx=0,
\]
which is (a). The functions $b_j$ have sup norm at most $\beta\norm{\nabla\zeta}_\infty$ and disjoint supports, and for every integer $m\ge0$ their sum satisfies $\norm{\sum_j\theta_jb_j}_{H^m}^2\le\sum_j\norm{b_j}_{H^m}^2\le C_m\beta^2\Delta^{-2m}$. Interpolation gives $\norm{\sum_j\theta_jb_j}_{H^\alpha}\le C\beta\Delta^{-\alpha}$, and the choice $\beta=c\,\varepsilon_0\Delta^{\alpha}$ with a small constant $c$ proves (b). For $z=(\bu,k)\in\Z_E$, integration by parts and $\partial_{\boldsymbol w}(\bu\cdot\bx-k)_+=(\bu\cdot\boldsymbol w)\mathbf1\{\bu\cdot\bx>k\}$ give
\begin{equation}\label{eq:absorbing-entries}
 C_{b_j}(z)=-\beta\Delta\,(\bu\cdot\boldsymbol w)\int\mathbf1_{H_z^+}\zeta_j\dd\bx,
 \qquad H_z^+=\{\bx:\bu\cdot\bx>k\}.
\end{equation}

\emph{The event.} Let $\gamma=R/(8q^*n)$ and $\rho_0=(8C_*q^*n^2)^{-1/d}$, where $C_*$ is such that $Q_E(B_\rho(z,t))\le C_*t^d$ for all $z\in\Z_E$ and $t>0$. Let $G_n$ be the event that no $K_i$ satisfies $r_1-\gamma\le\abs{K_i}<r_1$ and that $\rho(Z_i,Z_{i'})\ge\rho_0$ for all $i\ne i'$. Since the strike marginal of $Q$ has density at most $q^*/(2R)$, the first condition fails with probability at most $nq^*\gamma/R=\tfrac18$, and the second fails with probability at most $n^2q^*C_*\rho_0^d=\tfrac18$.

\emph{Reduction.} Write $H_i=\{\bx:U_i\cdot\bx=K_i\}$ and $H_i^+=\{\bx:U_i\cdot\bx>K_i\}$. On $G_n$, every index belongs to exactly one of $I=\{i:\abs{K_i}<r_1-\gamma\}$, $I_+=\{i:K_i\le-r_1\}$, and $I_0=\{i:K_i\ge r_1\}$. Let $\eta$ be as in (c) and $t_i=C_\eta(Z_i)$. For $n\ge n_0$ we have $r_0<r_1-\gamma$. Hence $H_i$ misses $\supp\eta$ for $i\in I_+\cup I_0$, and the vanishing mass and first moments of $\eta$ give $t_i=0$. By \eqref{eq:absorbing-entries}, $C_{b_j}(Z_i)=0$ for $i\in I_0$, and $C_{b_j}(Z_i)=-\beta\Delta(U_i\cdot\boldsymbol w)s_j$ with $s_j=\int\zeta_j=\Delta^d$ for $i\in I_+$. Hence it suffices to solve $M\boldsymbol c=(t_I,0)$, where $M$ has the rows $(C_{b_j}(Z_i))_j$ for $i\in I$ and the row $(s_j)_j$. If $\sigma=\min\{\abs{M^{\mathsf T}(\boldsymbol\lambda,\lambda'')}:\abs{\boldsymbol\lambda}^2+\lambda''^2=1\}$ is positive, the minimum-norm solution satisfies $\abs{\boldsymbol c}\le\sigma^{-1}\abs{t_I}\le\sigma^{-1}\sqrt n\max_i\abs{t_i}$. It remains to show $\sigma\ge n^{-C}$ on $G_n$.

\emph{Separated points.} Fix $i\in I$. On $G_n$, $H_i\cap B_{r_1-\gamma/2}$ is a $(d-1)$-dimensional disk of radius at least $(r_1\gamma/2)^{1/2}$, hence of area at least $A_n=c_d(r_1\gamma)^{(d-1)/2}$ for $\gamma\le r_1/2$. Let $\varrho>0$ and $i'\ne i$. Since $\bu\cdot\boldsymbol w\ge w_*$ on $E$, $\abs{U_i+U_{i'}}\ge2w_*$, hence $\sin\angle(U_i,U_{i'})=\tfrac12\abs{U_i-U_{i'}}\abs{U_i+U_{i'}}\ge w_*\abs{U_i-U_{i'}}$. If $\abs{U_i-U_{i'}}\ge\rho_0/4$, the set of $\bx\in H_i\cap B_R$ with $\abs{U_{i'}\cdot\bx-K_{i'}}<2\varrho$ is a slab of width at most $16\varrho/(w_*\rho_0)$, of area at most $C_dR^{d-2}\varrho/(w_*\rho_0)$. If $\abs{U_i-U_{i'}}<\rho_0/4$, then $\abs{K_i-K_{i'}}\ge3R\rho_0/4$ on $G_n$, and every $\bx\in H_i\cap B_R$ satisfies $\abs{U_{i'}\cdot\bx-K_{i'}}\ge\abs{K_i-K_{i'}}-R\abs{U_i-U_{i'}}\ge R\rho_0/2$, which exceeds $2\varrho$ when $\varrho\le R\rho_0/4$. Choosing $\varrho=c\min(\gamma,w_*\rho_0A_n/(nR^{d-2}))$ with a small constant $c$, the excluded area is at most $A_n/2$. There is therefore $\boldsymbol p_i\in H_i\cap B_{r_1-\gamma/2}$ with $\dist(\boldsymbol p_i,H_{i'})\ge2\varrho$ for all $i'\in I\setminus\{i\}$. For $i'\notin I$, the same holds since $H_{i'}$ misses $B_{r_1}$. Since the $2\varrho$-neighborhoods of all $H_{i'}$ cover at most $Cn\varrho$ of the volume of $B_{r_1/2}$, there is also $\boldsymbol p_0\in B_{r_1/2}$ at distance at least $2\varrho$ from every $H_{i'}$. The balls $B_i=B(\boldsymbol p_i,\varrho)$, $i\in I\cup\{0\}$, lie in $B_{r_1-\gamma/4}$, are pairwise disjoint, and $B_i$ meets no hyperplane other than $H_i$.

\emph{Lower bound on $\sigma$.} Let $\Delta=\varrho/(8\sqrt d)$ and fix $(\boldsymbol\lambda,\lambda'')$ of unit norm. By \eqref{eq:absorbing-entries}, $M^{\mathsf T}(\boldsymbol\lambda,\lambda'')=(\ip{F}{\zeta_j})_j$ with
\[
 F=-\beta\Delta\sum_{i\in I}\lambda_i(U_i\cdot\boldsymbol w)\mathbf1_{H_i^+}+\lambda''.
\]
On $B_i$, $i\in I$, the function $F$ takes two values that differ by $\beta\Delta\abs{\lambda_i}(U_i\cdot\boldsymbol w)\ge\beta\Delta w_*\abs{\lambda_i}$, one on each side of $H_i$. On $B_0$ it is a constant $v_0$, and $\abs{\lambda''-v_0}\le\beta\Delta\sqrt n\abs{\boldsymbol\lambda}$. Each side of $H_i$ in $B_i$ contains at least $c_d(\varrho/\Delta)^d$ of the supports of the $\zeta_j$, on each of which $\ip F{\zeta_j}=\Delta^dF$, and $B_0$ contains as many. Disjointness of the balls therefore gives
\[
 \abs{M^{\mathsf T}(\boldsymbol\lambda,\lambda'')}^2
 \ge c\,\varrho^d\Delta^d\left[\frac{w_*^2(\beta\Delta)^2}{2}\abs{\boldsymbol\lambda}^2+v_0^2\right]
 \ge c\,\varrho^d\Delta^d\min\left\{\frac{w_*^2(\beta\Delta)^2}{4},\frac{w_*^2}{8n}\right\},
\]
where the last step uses $v_0^2\ge\tfrac12\lambda''^2-(\beta\Delta)^2n\abs{\boldsymbol\lambda}^2$ with weight $w_*^2/(4n)$. Since $\gamma$, $\rho_0$, $A_n$, $\varrho$, $\Delta$, and $\beta$ depend only on $n$ and are fixed powers of $n$ up to constants, the family $b_1,\ldots,b_K$ is fixed before the design, $\sigma\ge n^{-C}$, and $K\le C\Delta^{-d}\le n^C$.
\end{proof}

\begin{proof}[Proof of Theorem~\ref{thm:exchange-uniform}]
\emph{Lower bound.} By Lemma~\ref{lem:noise-monotone} it suffices to treat $\varsigma=0$. Let $p=\alpha+\tfrac32$ and use the objects of the proof of the lower bound in \eqref{eq:log-modulus} with integer $N\ge1$ to be chosen: the packet $\eta_N$, the ridge $b\,\psi_NE_N$, and $m_N$. Choose $r_0<r_1<r$ so that $\supp\eta_N\subseteq B_{r_0}$ for all $N$, and choose $f_*$ positive on $B_{r_1}$, with minimum $c_*$ there. Let $b_1,\ldots,b_K$ and $G_n$ be as in Lemma~\ref{lem:absorbing} with $\varepsilon_0=c_*/4$. For $\btheta=(\theta_0,\theta_1,\ldots,\theta_K)\in[-1,1]^{K+1}$ set
\[
 f_{\btheta}=(1-m_N)f_*+b\,\psi_N(\bx')E_N(y-k)+\theta_0\eta_N+\sum_{j=1}^K\theta_jb_j.
\]
The first and last terms together are at least $c_*/4$ on $B_{r_1}$ and nonnegative elsewhere, and the middle two terms are nonnegative since $E_N\ge N^{-p}\abs{\Phi_N''}$. Hence $f_{\btheta}$ is a probability density supported in $B_r$, and its $H^\alpha$ norm is at most $\norm{f_*}_{H^\alpha}+Cb+\varepsilon_0\le L$. By Lemma~\ref{lem:absorbing}(a) and the price computation in the proof of the lower bound in \eqref{eq:log-modulus},
\[
 \int g_kf_{\btheta}\dd\bx=p_*+a_0\theta_0,\qquad a_0=bM_\psi N^{-p},
\]
for a constant $p_*$ not depending on $\btheta$. The exact quotes are $C_{f_{\boldsymbol0}}(Z_i)+(T\btheta)_i$, where $T$ has first column $t=(C_{\eta_N}(Z_i))_i$ and remaining columns $(C_{b_j}(Z_i))_i$.

Draw $\boldsymbol\Theta$ from the product prior \eqref{eq:compact-prior}, independently of the design and of any estimator randomization, and condition on both. On $G_n$, Lemma~\ref{lem:absorbing}(c) applied to $\eta_N$ yields $\boldsymbol c$ with $\sum_jc_jC_{b_j}(Z_i)=t_i$ and $\abs{\boldsymbol c}\le n^{C_3}\sup_{\Z_E}\abs{C_{\eta_N}}\le Cn^{C_3}e^{-cN}$, by the bound on $\epsilon_N$ in the proof of the lower bound in \eqref{eq:log-modulus}. Choose $N=\lceil C_4\log n\rceil$ with $C_4$ large enough that $\abs{\boldsymbol c}\le1$. Then $\bv=(1,-\boldsymbol c)/\abs{(1,-\boldsymbol c)}$ is a unit vector in $\ker(T)$ with $(a_0,0,\ldots,0)\cdot\bv\ge a_0/\sqrt2$. Clip an arbitrary estimator $\widehat p$ to $[p_*-a_0,p_*+a_0]$, which does not increase its error. Lemma~\ref{lem:functional-nullspace} gives a conditional mean squared error at least $a_0^2/(2\pi^2)$. Since the clipped error is at most $2a_0$, the conditional mean absolute error is at least $a_0/(4\pi^2)$ on $G_n$. Averaging over the design and bounding the prior average by the supremum,
\[
 \mathfrak R_n^Q\bigl(g_k;\mathcal C_{\alpha,L},0\bigr)\ge\Pr(G_n)\frac{a_0}{4\pi^2}\ge\frac{bM_\psi}{16\pi^2}\lceil C_4\log n\rceil^{-p}.
\]

\emph{Upper bound.} Let $h=n^{-1/(4d)}$, let $\{w_l\}$ be an $h$-net of $\Z_E$ in $\rho$ of cardinality at most $Ch^{-d}$, and let $\mathcal I_l=\{i:\rho(Z_i,w_l)\le h\}$. Put $\widehat C_l$ equal to the average of $Y_i$ over $\mathcal I_l$, or zero if $\mathcal I_l$ is empty. Let $\widehat f\in\F_\alpha(L)$ satisfy
\[
 \max_l\abs{C_{\widehat f}(w_l)-\widehat C_l}\le\inf_{f'\in\F_\alpha(L)}\max_l\abs{C_{f'}(w_l)-\widehat C_l}+\frac1n,
\]
chosen measurably by the countable-dense-set construction in the proof of Proposition~\ref{prop:noisy-upper}, and estimate the price by $\int g_k\widehat f$. By \eqref{eq:ahlfors} and $Q\ge q_*Q_E$, each ball $B_\rho(w_l,h)$ has $Q$-probability at least $cq_*h^d$. A Chernoff bound gives an event of probability at least $1-n^{-1}$ on which every $\abs{\mathcal I_l}\ge cq_*nh^d$. On this event, the averages $\bar\xi_l$ of the errors over $\mathcal I_l$ have conditional variances at most $(cq_*nh^d)^{-1}$, hence $\E\max_l\varsigma\abs{\bar\xi_l}\le\overline\varsigma(Ch^{-d}/(q_*nh^d))^{1/2}\le C\overline\varsigma n^{-1/4}$, and Markov's inequality gives $\max_l\varsigma\abs{\bar\xi_l}\le n^{-1/8}$ outside an event of probability at most $Cn^{-1/8}$. Since $C_f$ is $R$-Lipschitz in $\rho$, on the remaining event the true density $f$ satisfies $\max_l\abs{C_f(w_l)-\widehat C_l}\le\delta_n:=Rh+n^{-1/8}$. Hence $\max_l\abs{C_{\widehat f}(w_l)-C_f(w_l)}\le2\delta_n+1/n$, and the $2R$-Lipschitz property extends this to $\sup_{\Z_E}\abs{C_{\widehat f}-C_f}\le\epsilon_n:=2\delta_n+1/n+2Rh\le Cn^{-1/(8d)}$. Theorem~\ref{thm:log-modulus} bounds the price error by $\varpi_{g_k}^{\alpha,L}(\epsilon_n)\le C(\log(1/\epsilon_n))^{-p}\le C'(\log n)^{-p}$ for $n\ge n_0$. Off these events, the error is at most $2R$, which contributes at most $Cn^{-1/8}$.
\end{proof}

\begin{proof}[Proof of Theorem~\ref{thm:directional}]
(a) Fix $\bv\in D$. Let $N\ge1$ be an integer, $p=R/N$, and let $\mu$ and $\nu$ be supported on the segment $\{t\bv:\abs t\le R/2\}$ with densities $R^{-1}(1+\cos(2\pi t/p))$ and $R^{-1}$ in $t\in[-R/2,R/2]$. The difference $\eta=\mu_{\bv}-\nu_{\bv}$ has zero mass and zero mean, since the interval contains whole periods and the cosine is even. Its call function $\Xi(l)=\int(t-l)_+\eta(\mathrm dt)$ satisfies $\Xi''=\eta$ and vanishes for $\abs l\ge R/2$, hence $\Xi(l)=(p^2/(4\pi^2R))\left((-1)^N-\cos(2\pi l/p)\right)$ on $[-R/2,R/2]$ and $\sup\abs\Xi=\epsilon_N:=p^2/(2\pi^2R)$. For $\bu\in S^{d-1}$ and $\lambda=\bu\cdot\bv$, the call difference $C_\mu(\bu,k)-C_\nu(\bu,k)=\int(\lambda t-k)_+\eta(\mathrm dt)$ equals $\lambda\Xi(k/\lambda)$ if $\lambda>0$, equals $\abs\lambda\Xi(-k/\abs\lambda)$ if $\lambda<0$ by put--call parity and the vanishing mass and mean of $\eta$, and vanishes if $\lambda=0$. Hence $\sup_{\Z}\abs{C_\mu-C_\nu}\le\epsilon_N$, while $W_1(\mu_{\bv},\nu_{\bv})=\int\abs{\Xi'}=p/\pi^2=\pi^{-1}\sqrt{2R\epsilon_N}$. Monotonicity of $\varpi_D$ between consecutive $\epsilon_N$, whose square roots differ by the factor $N/(N+1)\ge\tfrac12$, gives the lower bound. The upper bound is Proposition~\ref{prop:positive}(ii), which is uniform over $1$-Lipschitz $G$ and $\bv\in E$.

(b) Use the perturbation \eqref{eq:oscillation} rotated to $\be$, as in the proof of Theorem~\ref{thm:missing-cone}, with integer $N\ge1$ and the bump $\zeta_N$ of Lemma~\ref{lem:profile} with $d'=d$, with common support in $B_r$; the uniform Sobolev bounds on $\zeta_N$ keep the densities $f_N^\pm$ in $\F_\alpha(L)$. The transform $\widehat{\cA H_N}(\bu,t)=\tfrac b2N^{-\alpha-2}(\bu\cdot\be)^2(\widehat\zeta_N(t\bu-N\be)+\widehat\zeta_N(t\bu+N\be))$, the bound $\abs{t\bu\pm N\be}\ge c_\kappa(\abs t+N)$ for $\bu\in E$, and Lemma~\ref{lem:profile} with $c_0=c_\kappa$ and $J=2$ give $\sup_{\Z_E}\abs{C_{f_N^+}-C_{f_N^-}}\le Ce^{-cN}$. The test claim $-N^{-1}\cos(N\be\cdot\bx)$ lies in $\mathcal G_D$, and its integral against $f_N^+-f_N^-$ is $2bN^{-\alpha-1}\int\zeta_N\cos^2(N\be\cdot\bx)\dd\bx=bN^{-\alpha-1}(\int\zeta_N+\operatorname{Re}\widehat\zeta_N(2N\be))\ge cN^{-\alpha-1}$ for large $N$. Choosing $N\asymp\log(1/\epsilon)$ gives the bound.

(c) For $\bv\in D\cap E$ use (a). For $\bv=\be$, let $f_1,f_2\in\F_\alpha(L)$ have discrepancy at most $\epsilon$ on $\Z_E$, let $\eta=(f_1-f_2)\dd\bx$, and let $\Psi(y)=\eta_{\be}((-\infty,y])$. Then $W_1(\mu_{\be},\nu_{\be})=\int_{-r}^r\abs\Psi\le\sqrt{2r}\norm\Psi_{L^2}$ and $\widehat\Psi(\tau)=\widehat\eta(\tau\be)/(\mathrm i\tau)$. Lemma~\ref{lem:moments} with $k=0$ bounds the moments $\mathsf m_j=\int(\be\cdot\bx)^j\dd\eta$ by $(3j/R)\varrho^j\epsilon$, hence $\abs{\widehat\eta(\tau\be)}\le\sum_{j\ge1}\abs\tau^j\abs{\mathsf m_j}/j!\le(3\varrho/R)\epsilon\abs\tau e^{\varrho\abs\tau}$. Hence the frequencies $\abs\tau\le N$ contribute at most $C\epsilon e^{\varrho N}$ to $\norm{\widehat\Psi}_{L^2}$, and Lemma~\ref{lem:line-restriction} bounds the frequencies $\abs\tau>N$ by $CLN^{-\alpha-1}$. Taking $N=\log(1/\epsilon)/(2\varrho)$ gives $W_1(\mu_{\be},\nu_{\be})\le C(\epsilon^{1/2}+L(\log(1/\epsilon))^{-(\alpha+1)})$. Part (b) gives the matching lower bound when $\be\in D$.
\end{proof}

\begin{proof}[Proof of Corollary~\ref{cor:indexed-finite}]
Two laws in the same admissible set have discrepancy at most $2Rh_E$ on $\Z_E$, as in Corollary~\ref{cor:finite}, which also bounds $h_E$. Apply Theorem~\ref{thm:directional}.
\end{proof}

\putbib
\end{bibunit}
\end{appendix}

\end{document}